\documentclass[11pt]{article}

\usepackage{amsmath, amsfonts,amssymb, amsthm, euscript,makeidx,xcolor,mathrsfs,latexsym}

\newcommand{\la}{\langle}
\newcommand{\ra}{\rangle}

\newcommand{\ol}{\overline}
\newcommand{\ul}{\underline}
\newcommand{\eps}{\varepsilon}

\newcommand{\ba}{\begin{array}}
\newcommand{\ea}{\end{array}}
\newcommand{\be}{\begin{equation}}
\newcommand{\ee}{\end{equation}}
\newcommand{\bea}{\begin{eqnarray}}
\newcommand{\eea}{\end{eqnarray}}
\newcommand{\beaa}{\begin{eqnarray*}}
\newcommand{\eeaa}{\end{eqnarray*}}

\newcommand{\smax}{\operatornamewithlimits{strict-max}}

\def\dbD{\mathbb{D}}
\def\dbE{\mathbb{E}}
\def\dbF{\mathbb{F}}

\def\dbL{\mathbb{L}}

\def\dbN{\mathbb{N}}
\def\dbP{\mathbb{P}}
\def\dbR{\mathbb{R}}
\def\dbS{\mathbb{S}}

\def\dbX{\mathbb{X}}

\def\a{\alpha}
\def\b{\beta}
\def\g{\gamma}
\def\d{\delta}
\def\e{\varepsilon}

\def\k{\kappa}
\def\l{\lambda}

\def\si{\sigma}

\def\f{\varphi}
\def\th{\theta}
\def\o{\omega}

\def\arm{{\rm a}}
\def\G{\Gamma}
\def\D{\Delta}

\def\L{\Lambda}

\def\O{\Omega}
\def\cA{{\cal A}}
\def\cB{{\cal B}}

\def\cF{{\cal F}}

\def\cI{{\cal I}}
\def\cJ{{\cal J}}

\def\cL{{\cal L}}
\def\cM{{\cal M}}

\def\cP{{\cal P}}

\def\cX{{\cal X}}

\def\no{\noindent}

\def\ss{\smallskip}
\def\ms{\medskip}
\def\bs{\bigskip}
\def\q{\quad}
\def\qq{\qquad}

\def\pa{\partial}
\def\cd{\cdot}
\def\cds{\cdots}

\def\tr{\hbox{\rm tr}}

\def\qed{ \hfill \vrule width.25cm height.25cm depth0cm\smallskip}

\newcommand{\basa}{\begin{assumption}}
\newcommand{\easa}{\end{assumption}}

\newcommand{\bas}{\begin{assum}}
\newcommand{\eas}{\end{assum}}

\def\limsup{\mathop{\overline{\rm lim}}}
\def\liminf{\mathop{\underline{\rm lim}}}

\def\pa{\partial}

 \def\cd{\cdot}
\def\cds{\cdots}

\def\tr{\hbox{\rm tr$\,$}}

\def\dis{\displaystyle}

\def\bx{{\bf x}}

\def\1{{\bf 1}}

\def\:{\!:\!}
\def\reff{\eqref}
\def \proof{{\noindent \bf Proof.\quad}}

\numberwithin{equation}{section}

\allowdisplaybreaks

  \renewcommand\appendix{\par
    \setcounter{section}{0}
    \setcounter{subsection}{0}
    \gdef\thesection{ Appendix \Alph{section}}}

    \newtheorem{thm}{Theorem}[section]
\newtheorem{lem}[thm]{Lemma}

\newtheorem{prop}[thm]{Proposition}
\newtheorem{rem}[thm]{Remark}
\newtheorem{eg}[thm]{Example}
\newtheorem{defn}[thm]{Definition}
\newtheorem{assum}[thm]{Assumption}

\title{Viscosity Solutions for HJB Equations on the Process Space} 

\author{Jianjun Zhou\footnote{College of Science, Northwest A$\&$F University,Yangling 712100, Shaanxi, P. R. China, zhoujianjun@nwsuaf.edu.cn. This author is partially supported by  Shaanxi Natural Science Foundation Grant 2025JC-YBMS-021.} \quad Nizar Touzi\footnote{New York University, Tandon School of Engineering, United States, nizar.touzi@nyu.edu.}  \quad  and\q Jianfeng Zhang\footnote{Department of Mathematics, University of Southern California, United States, jianfenz@usc.edu. This author is partially supported by NSF grant DMS-2205972. }}

      \date{}

\begin{document}

\maketitle

\pagestyle{plain}

\begin{abstract}

In this paper we investigate a path dependent optimal control  problem on the process space with both drift and volatility controls,  with possibly degenerate volatility. The dynamic value function is characterized by a fully nonlinear second order path dependent HJB equation on the process space, which is by nature infinite dimensional. In particular, our model  covers mean field control problems with common noise as a special case. We shall introduce a new notion of viscosity solutions and establish both the existence and the comparison principle, under merely Lipschitz/H\"{o}lder continuity assumptions. The main feature of our notion is that, besides the standard smooth part, the test function consists of an extra singular component  which allows us to handle the second order derivatives of the smooth test functions without invoking the Crandall-Ishii lemma. We shall use the doubling variable arguments, combined with the Ekeland-Borwein-Preiss variational principle in order to overcome the noncompactness of the state space. A smooth gauge-type function on the path space is crucial for our estimates.
\medskip

\no{\bf Keywords:}  Mean field control, common noise, viscosity solutions, HJB equations, Wasserstein space, path dependence, comparison principle, variational principle
\end{abstract}


\vfill\eject

\section{Introduction}
In this paper we consider a stochastic control problem whose data rely on the whole underlying state process $X$,  namely on the mapping $X:  [0, T]\times \O\to \dbR^d$ rather than on the paths $X(\o)$ for each $\o\in \O$. Consequently, its dynamic value function satisfies an infinite dimensional HJB equation on the space of processes $X$. This setting in particular covers mean field control problems, which are law invariant in terms of the state process.  In the state dependent case, the idea of lifting probability measures to random variables is due to Lions \cite{lio4}. Our framework is very general and has the following features:

\begin{itemize}
\item We allow for path dependence and thus our HJB equation is path dependent. For finite dimensional path dependent PDEs, we refer to a series of works by the authors and their collaborators: \cite{EKTZ, ETZ1, ETZ2, RTZ, zhou}, and the references therein. We also refer to Wu-Zhang \cite{WZ} and Cosso-Gozzi-Kharroubi-Pham-Rosestolato \cite{cosso3}   for path dependent mean field control problems.

\item We consider both the drift and volatility controls, with possibly degenerate volatilities. Consequently, besides infinite dimensionality, our HJB equation is a fully nonlinear degenerate second order PDE.  {There have been interesting attempts for mean field control problems with volatility controls, see e.g.  Cosso-Gozzi-Kharroubi-Pham-Rosestolato \cite{CGKPR}, Cheung-Tai-Qiu \cite{CTQ}, and Bayraktar-Ekren-Zhang \cite{BEZ}.}

\item Our model covers mean field control problems with common noise as a special case.  {We note that \cite{BEZ, CTQ}, Gangbo-Mayorga-Swiech \cite{GMS}, Mayorga-Swiech \cite{MS}, and Daudin-Jackson-Seeger \cite{DJS} also considered common noise.} 
Moreover, in the contexts of mean field games, Ahuja-Ren-Yang \cite{ARY}  {and Cardaliaguet-Souganidis  \cite{CS}} employed the lifting idea to study common noise.

\item We allow the problem to depend on the joint law of the state process and the control process. This is in the spirit of  mean field games of controls (which were called extended mean field games in the early stage), see e.g.  Gomes-Voskanyan \cite{GV}.
\end{itemize}
Initiated independently by Caines-Huang-Malhame \cite{CHM} and Lasry-Lions \cite{LL}, the theory of  mean field games and mean field controls  has received extremely strong attention in the literature. We refer to Lions's lecture \cite{lio4}, the books Bensoussan-Frehe-Yam \cite{BFY}, Carmona-Delarue \cite{CD1, CD2}, Cardaliaguet-Porretta \cite{CP}, and the references therein for a general exposition of the theory. One popular approach in the literature is to consider PDEs on the Wasserstein space of probability measures. We should note that the master  equations for mean field games  have quite different nature than the HJB equations for  mean field controls. In this paper we address the latter equations and our main focus is the comparison principle for viscosity solutions, while for mean field game master equations even classical solutions typically violate the comparison principle.

Partially due to its infinite dimensionality, such an HJB equation can rarely have a classical solution. In recent years there have been serious efforts on viscosity solutions for HJB equations arising from mean field control problems. We shall provide a literature review in the end of this introduction.  Our goal of this paper is to propose an appropriate notion of viscosity solutions and show that the dynamic value function of our control problem is the unique viscosity solution of the HJB equation.  In particular, we shall establish the comparison principle for viscosity solutions under merely Lipschitz/H\"{o}lder continuity assumptions, by using the doubling variable arguments. Our results also imply the wellposedness of the second order HJB equation on the Wasserstein space of probability measures induced from mean field control problems with common noise. To the best of our knowledge, even in the mean field control framework, our model is most general, covering most models studied in the literature, and our technical conditions are the weakest.

Notice that a notion of viscosity solutions is essentially determined by the set of test functions one chooses. In particular, the proof of the comparison principle relies heavily on this choice of test functions. Inspired by the viscosity solution theory for PDEs in infinite dimensions, see, e.g., Crandall-Lions \cite[Definition 2.1]{cran1}, Li-Yong \cite[Chapter 6 Definition 3.1]{LY}, and Fabbri-Gozzi-{Swiech} \cite[Definitions 3.32 \& 3.35]{fab1},  {as well as the viscosity solution theory for   PDEs with discontinuous time-dependence, see, e.g., Ishii \cite{Ishii}, Lions-Perthame \cite{LP}, and Nunziante \cite{Nunziante1, Nunziante2},} our test functions take the form
\bea
\label{testfunction}
\f + \phi,
\eea
where the first part $\f$ is smooth (in appropriate sense) and thus is standard; and the extra part $\phi$ is not smooth, but is absolutely continuous in time with respect to the Lebesgue measure. Our construction of $\phi$ is motivated from an observation in the constant volatility case. In this case, by a simple transformation one can convert the state process from a controlled SDE to a controlled ODE with random coefficients. Consequently, the resulting HJB equation becomes first order,\footnote{\label{1storder}
In the literature, the order often refers to the derivatives with respect to the measure variable $\mu$. Then mean field control problems with idiosyncratic noise
 are related to first order equations while those with common noise are related to second order equations.
  Here by first order we mean the derivatives with respect to $x$ are of the first order. In particular, 
  we view the Lions derivative $\pa_\mu = \pa_x {\d\over \d \mu}$ as first order in $x$, where ${\d\over \d\mu}$ is the linear functional derivative. 
   So our first order equation corresponds to mean field control problems with neither idiosyncratic noise nor common noise. We shall take this convention 
   throughout the paper.}  whose viscosity solution is a lot easier to study. By applying the inverse transformation on the test functions for the latter equation,
   we obtain a candidate test function $\phi$ for our original dynamic value function, which turns out to be absolutely continuous in time but  
   is in general not smooth. See Subsection \ref{sect-firstorder} for details. 
   For the general case with volatility controls, we tailor the construction of $\phi$ which in essence cancels the diffusion term. 
    Indeed, technically $\phi$ is used to cancel some terms appearing in the doubling variable arguments, which involve the second order 
     derivatives\footnote{\label{2ndorder} As in Footnote \ref{1storder}, here we refer to the second order derivatives with respect to $x$, 
      which are due to the volatility of the state process, or say due to the presence of idiosyncratic noise or common noise.}
of $\f$ and are otherwise hard to estimate.  To the best of our knowledge, this type of test functions is new in the literature of viscosity solutions for mean field control problems.

Another important consequence of introducing the singular component $\phi$ of test functions is that we can establish the comparison principle without using the Crandall-Ishii lemma, even though we are  using the doubling variable arguments for a  second order equation. Indeed, as we just explained, in the constant volatility case, the introduction of $\phi$ allows us to convert the HJB equation into a first order one, which does not require the Crandall-Ishii lemma. For standard second order HJB equations, the Crandall-Ishii lemma is used exactly to handle the second order derivatives of the test functions. So our approach provides an alternative solution to this important issue. Moreover, our general framework covers the standard finite dimensional (path dependent) HJB equations, as well as HJB equations on the Wasserstein space of probability measures, so our results imply the comparison principle for viscosity solutions of those equations as well, without using the Crandall-Ishii lemma. However, we should point out that our notion of viscosity solutions is not equivalent to the  ``standard" ones for those equations. So our results do not imply directly the wellposedness results in the literature.

Unfortunately, even in the mean field framework, this function $\phi$ is typically not law invariant, which prohibits us from defining viscosity solutions intrinsically on the Wasserstein space. Thus we are required to consider functions on the process space, which as a by product enables us to cover the common noise case for free. Another advantage for working directly on the  process space is that, in the state dependent case, the square distance of two random variables is a smooth functional, but the square $2$-Wasserstein distance of two probability measures is not differentiable. While being more involved, our path dependent case benefits from this feature as well.  We should note that Soner-Yan \cite{SY1, SY2} introduced a nice norm on the Wasserstein space by using the Fourier transform, whose square is smooth. See also \cite{CGKPR}, and Bayraktar-Ekren-Zhang \cite{BEZ1}, Daudin-Seeger \cite{DS} for serious efforts to get around of this difficulty. However, it is not clear how to extend these ideas to serve for our purpose in the general case, especially in the path dependent setting.

One drawback of  the process space is its lack of local compactness, which is crucial for the proof of the comparison principle.  To circumvent this difficulty, we shall use the Ekeland-Borwein-Preiss variational principle, see Aubin-Ekeland \cite{AE} and Borwein-Preiss \cite{bor1}. Roughly speaking, to optimize a continuous function on a non-compact space, one may construct an approximate function by using the so called gauge type function such that the approximate function has a strict optimal argument. To serve as a test function for our purpose, we require the gauge type  function on the process space to be smooth with desired estimates for its derivatives. This is achieved by utilizing the smooth gauge type function on continuous paths constructed by Zhou \cite{zhou} for path dependent PDEs.

To prove the comparison principle, we shall first double the spatial variable and then double the temporal variable. This procedure is standard for parabolic equations, see, e.g., Crandall-Ishii \cite[Lemma 8]{cran-ishii}  and Crandall-Ishii-Lions \cite[Theorem 8.3]{cran2}. However, due to a subtle adaptedness requirement of $\phi$, we need a third step of approximation to derive the desired contradiction.

The rest of the paper is organized as follows. First, we conclude this introduction by providing a brief literature review on viscosity solutions for mean field control problems. In Section \ref{sect-model}  we introduce our control problem and establish some basic properties of its dynamic value function.  In Section \ref{sect-test} we introduce smooth functions $\f$ on the process space, which leads to the target HJB equation. In Section \ref{sect-phi} we construct the singular component $\phi$ and derive some crucial estimates. In Section \ref{sect-HJB} we propose our notion of viscosity solutions and present the main results. In Section  \ref{sect-MFC0} we illustrate how our general model covers the mean field control problem with common noise as a special case. Sections \ref{sect-comparison1}, \ref{sect-comparison2} and \ref{sect-comparison3} are devoted to the proof of the comparison principle for viscosity solutions. Finally in the Appendix we complete some technical proofs.

\ms
\no{\bf Some literature review on viscosity solutions for mean field control problems.}  We first remark that these equations are by nature infinite dimensional, and thus the related works are intrinsically connected to the viscosity solution theory for infinite dimensional HJB equations, for which we refer to Lions \cite{lio1, lio2, lio3} and the books Li-Yong \cite{LY} and Fabbri-Gozzi-{Swiech} \cite{fab1}. For first order (in the sense of Footnote \ref{1storder}) HJB equations on the Wasserstein space, arising from mean field control problems with deterministic controls and other related problems, we refer to the works Bertucci \cite{Bertucci}, Cardaliaguet-Quincampoix \cite{CQ}, Conforti-Kraaij-Tonon \cite{CKT},  Feng-Katsoulakis \cite{FK},  Gangbo-Nguyen-Tudorascu \cite{GNT}, Gangbo-Tudorascu \cite{GT}, and Jimenez-Marigonda-Quincampoix \cite{JMQ1, JMQ2}. These equations involve the Lions derivative $\pa_\mu V(t,\mu, x)$, but not the higher order derivative $ \pa_{x\mu} V(t, \mu, x)$, and in the finite dimensional case they correspond to the standard first order HJB equations.

We are mainly interested in mean field control problems with idiosyncratic noise and/or common noise, where the state process is a controlled diffusion, with or without volatility controls. These equations involve $\pa_{x\mu} V(t, \mu, x)$ and/or $\pa_{\mu\mu} V$, and in the finite dimensional case correspond to the standard  second order HJB equations. In the path dependent setting, Wu-Zhang \cite{WZ} proposed a notion of viscosity solutions by restricting the viscosity neighborhood of some point $(t, \mu)$  to certain compact set, and established the partial comparison principle. However, the full comparison principle requires certain perturbed equation to have a classical solution, which is verified only in some special cases.
In a state dependent setting with volatility controls, the work  Cosso-Gozzi-Kharroubi-Pham-Rosestolato \cite{CGKPR} removed the compactness requirement by applying the Ekeland-Borwein-Preiss variational principle. Based on a so called Gaussian-smooothed $2$-Wasserstein distance,  { \cite{CGKPR} constructed a finite dimensional smooth approximation of the value function, which allows one to
compare viscosity semi-solutions with the value function. Following this approach but by modifying the definition of viscosity solutions, the recent paper Cheung-Tai-Qiu \cite{CTQ} established the comparison principle rigorously. The work Bayraktar-Cheung-Ekren-Qiu-Tai-Zhang \cite{BCEQTZ} continued along this line.} The works Burzoni-Ignazio-Reppen-Soner \cite{BIRS} and Soner-Yan \cite{SY1} derived the comparison result by means of the doubling variable arguments.
 They did not invoke the Crandall-Ishii lemma either, but due to a completely different nature than ours, see Remark \ref{rem-Ishii2} below for the detailed explanation. In particular, they introduced a smooth metric on the Wasserstein space by using the Fourier transform. However, these works require certain uniform Lipschitz continuity on the controls. When the volatility is a positive constant and hence the equation is semilinear, this serious constraint was removed in the recent work Soner-Yan \cite{SY2} and the comparison principle was established. More recently, by also using the Fourier-Wasserstein metric, Bayraktar-Ekren-Zhang \cite{BEZ} and Bayraktar-Ekren-He-Zhang \cite{BEHZ} established a Crandall-Ishii lemma for functions on the Wasserstein space of probability measures, and considered a mean field control problem with partial information and common noise. In the case with common noise  but without idiosyncratic noises, Gangbo-Mayorga-Swiech \cite{GMS}  and Mayorga-Swiech \cite{MS} studied the so called $L$-viscosity solution by lifting the equation to the Hilbert space.  The work Daudin-Seeger \cite{DS} studied semilinear HJB equations, by applying the doubling variable arguments with a further entropy penalization. Another recent work Daudin-Jackson-Seeger \cite{DJS} studied semilinear equations with common noise and non-convex Hamiltonian, which allows to consider zero-sum game problem in the mean field setting. They established the  comparison principle by exploiting the idea of \cite{BEZ} together with some delicate regularity estimates.

 { We remark that, when restricting to the mean field control setting, our model covers most models mentioned in the previous paragraph, with the exceptions that \cite{WZ} considers general parabolic equations, \cite{DJS} considers zero-sum game problems, and \cite{BEZ, BEHZ} consider problems with partial information. We believe  our approach can be extended to cover more general cases. Moreover, we require only Lipschitz or even H\"{o}lder continuity on the data, under $W_2$ for the measure variable $\mu$. To the best of our knowledge, even when restricting to the settings in those works, our technical assumptions are the weakest, for example:
\begin{itemize}

\item The works \cite{BIRS,  DS, DJS, GMS, MS, SY1, SY2} require the volatility coefficient (for idiosyncratic noise or for common noise) to be uniformly non-degenerate, or even to be a constant;

\item The works \cite{BCEQTZ, CTQ, CGKPR} do not allow the volatility coefficient to depend on the measure variable $\mu$; and the works \cite{BEHZ, BEZ} require a strong technical assumption, which is verified only when all the coefficients are independent of $\mu$;

\item The works \cite{BCEQTZ, CTQ, DJS, GMS, MS} require the common noise coefficient to be independent of control, or even to be a constant; and \cite{BEHZ, BEZ} verified their crucial condition only in the case that the common noise coefficient depends only on the control, but not on the state or its law;

\item The works \cite{BCEQTZ,  BEHZ, BEZ, BIRS,  CTQ,  CGKPR, DS, DJS, SY1, SY2} require the coefficients to be $W_1$-Lipschitz continuous in $\mu$, while \cite{GMS, MS} require $W_p$-Lipschitz continuity  for some $p<2$.
\end{itemize}
}

We should also mention the following works concerning potential mean field games, where the mean field control problem is involved automatically: Bensoussan-Graber-Yam \cite{BGY1, BGY2}, Bensoussan-Tai-Yam \cite{BTY}, Carmona-Cormier-Soner \cite{CCS}, Cecchin-Delarue \cite{CD}, and Gangbo-Meszaros \cite{GM}.  {In particular, \cite{BGY1, BGY2, BTY} worked directly on the Hilbert space of random variables. Following the approach in \cite{WZ}, Talbi-Touzi-Zhang \cite{TTZ1, TTZ2} established the complete wellposedness for a mean field optimal stopping problem. By using a finite dimensional projection and modifying the standard Crandall-Ishii lemma, Soner-Tissot Daguette-Zhang \cite{STDZ} proved the comparison principle for an HJB equation arising from controlled occupied process, which involves a special type of path dependence.} Moreover, the works Cox-Kallblad-Larsson-Svaluto-Ferro \cite{CKLS} on controlled measure-valued martingales and Feng-Swiech \cite{FS} where the controlled dynamics involves a mixture of a Hamiltonian flow and a  gradient flow   are also closely related.

\ms

\no{\bf Some notations.} We shall denote $x \cd x' := \sum_{i=1}^n x_i x'_i$ for $x, x'\in \dbR^n$, and $M \!:\! M' := \tr(M^\top M')$  for $M, M'\in \dbR^{m\times n}$, with $M^\top$ the transpose of $M$. Moreover, $|x|^2 := x\cd x$, $|M|^2 := M \!:\! M$. The set of $d\times d$-symmetric matrices is denoted by $\dbS^d$.

 {Throughout the paper, we fix a finite time horizon $[0, T]$ and a filtered probability space $(\O, \cF, \dbF, \dbP)$, with $\dbF=\{\cF_t\}_{0\le t\le T}$.} We assume $\cF_0$ is rich enough to support any probability measure on $\dbR^d$, and $\cF_t = \cF_0 \vee \cF^B_t$, where $B$ is a $d-$dimensional Brownian motion on $(\Omega,\cF,\dbP)$. Moreover, we fix a sub-filtration $\dbF^0=\{\cF^0_t\}_{0\le t\le T} \subset \dbF$, and denote $\dbE^0_t := \dbE[\cd|\cF^0_t]$.

For a Euclidian space $E$ and $p\ge 1$, we denote by $\dbL^p(\cF_{t}; E)$ the space of $\cF_{t}$-measurable $E$-valued random variables $\xi_*$ such that $\|\xi_*\|_p^p:= \dbE[|\xi_*|^p] <\infty$; and  $\dbL^p(\dbF_{[t, T]}; E)$ the space of $\dbF$-progressively measurable $E$-valued processes $\xi$ on $[t, T]$ such that $\dbE\big[\int_{t}^T |\xi_s|^pds\big] <\infty$.

Denote $\dbX:= C([0, T]; \dbR^d)$, equipped with the uniform norm $|\cd|_\infty$.  For any $p\ge 1$, let  $\cX_p$  denote the set of $\dbF$-progressively measurable continuous processes $\xi$ with $\|\xi\|_p^p :=\dbE[|\xi|^p_\infty] < \infty$, equipped with the norm $\|\cd\|_p$. Note that $\xi(\o) \in \dbX$ for all $\o\in \O$. Let $\ol\cX_{[t, T]}^p := [t, T]\times \cX_p$ denote the time and state space.

In order to distinguish the dependence on the whole process, or more precisely on the (deterministic) mapping on $[0, T]\times \O$,  from that on the realized paths of the process,  we introduce the notation $\ul\xi=\xi$ to emphasize the dependence on the whole process. That is, for a function $\f$ on $[0, T] \times \dbX \times \cX_p$,  we shall write $\f_t( \bx, \ul\xi)$ instead of $\f_t( \bx,\xi)$. In particular, this allows us to express the following without confusion:
\beaa
\f_t(\xi, \ul\xi)(\o) = \f_t(\xi(\o), \ul\xi),
\eeaa
and, when not involving $\bx$, the value $\f_t( \ul\xi)$ is deterministic. When it is more convenient, especially when the functions are state dependent, we may also use the notation $\f(t, \bx, \ul\xi) = \f_t(\bx, \ul\xi)$.
Moreover, throughout the paper we shall always assume all involved path dependent functions  $\f$ are adapted in the sense:
\bea
\label{adapted}
\f_t( \bx, \ul\xi) = \f_t\big( \bx_{\cd\wedge t}, \ul\xi_{\cd\wedge t}\big),
&\mbox{for all}&
(t,\bx,\xi)\in[0,T]\times\dbX\times\cX_p.
\eea

\section{Formulation of the process dependent control problem}
\label{sect-model}
\setcounter{equation}{0}

We will consider open loop controls taking values in $A$, a  { (possibly unbounded)} domain in a Euclidian space.\footnote{Although closed loop controls are expected to induce the same value function under appropriate regularity conditions, we refrain from considering this case for technical simplicity.}  Denote $\cA_{[t, T]}:= \dbL^2(\dbF_{[t, T]}; A)$,  $\cA_t:= \dbL^2(\cF_t; A)$, and consider the data:
\bea
\label{coefficients}
\left.\ba{c}
(b, \si): [0, T]\times \O\times \dbX \times  A \times \cX_2\times \cA_{T}
\longrightarrow
(\dbR^d, \dbR^{d\times d}),\\
 f: [0, T]\times \cX_2\times \cA_{T}
 \longrightarrow
 \dbR,
 \qq
 g: \cX_2 \longrightarrow \dbR.
\ea\right.
\eea
As usual we omit the variable  $\o\in \O$ inside $b$ and $\si$, and we recall the convention \reff{adapted}.
For any $(t,\xi)\in \ol \cX^2_{[0, T]}$ and $\a\in \cA_{[t,T]}$, consider the following path dependent SDE on $[t, T]$:
\bea\label{state1}
\left.\ba{c}
\dis X_s^{t,\xi,\a}=\xi_s,\q s\in [0, t];\ms\\
\dis dX^{t,\xi,\a}_s= b_s\big(X^{t,\xi,\a},\a_s, \ul X^{t,\xi,\a}, \ul \a_s\big)ds+ \si_s\big(X^{t,\xi,\a},\a_s, \ul X^{t,\xi,\a}, \ul\a_s\big)dB_s,\q s\in [t, T].
\ea\right.
\eea
Here $X$ is path dependent, while $\a$ involves only the current state, and the notation $\ul\a_s$ refers to the whole random variable $\a_s$. The value function of our control problem is:
\bea\label{value1}
\left.\ba{c}
\dis V_t( \ul\xi):=\inf_{\a\in \cA_{[t,T]}}J_t( \ul\xi,\ul\a), \q (t,\xi)\in \ol\cX_{[0, T]}^2,\\
\dis \mbox{where}\q J_t( \ul\xi, \ul\a) : = g({\ul X^{t, \xi,\a}}) +  \int_{t}^{T} f_s( \ul X^{t,\xi,\a}, \ul\a_s)ds.
 \ea\right.
 \eea

Throughout the paper, the following assumptions will always be in force.
\begin{assum}
\label{assum-standing}
{\rm (i)} The coefficients $b, \si, f$ are progressively measurable in all variables and adapted in the sense of \reff{adapted};  in particular, $b, \si$ are $\dbF$-progressively measurable; and  $h_t({\bf 0}, \a_t, \ul {\bf 0}, \ul\a_t)$, for  $h=b,\si$, $f_t(\ul {\bf 0}, \ul\a_t)$ and $g(\ul {\bf 0})$ are bounded \footnote{ This boundedness requirement is just for technical convenience. In particular, it can be relaxed so as to cover the linear quadratic case.} by a constant $C_0$.

\no {\rm (ii)} $b, \si$ are uniformly Lipschitz continuous in $(\bx, \xi)$ with a Lipschitz constant $L$:
      \bea
   &\dis |h_t(\bx,a, \ul\xi,\ul \a_t)-h_t( \bx',a,\ul\xi,\ul \a_t)| \le L|\bx_{\cd\wedge t} - \bx'_{\cd\wedge t}|_\infty,\q h=b, \si;\ms\nonumber\\
   \label{LipschitzCond}
    &\dis   |h_t(\bx,a, \ul\xi,\ul \a_t)-h_t( \bx,a,\ul\xi',\ul \a_t)| \le L \Big(\dbE^0_t\big[|\xi_{\cd\wedge t}-\xi'_{\cd\wedge t}|_\infty^2\big]\Big)^{\frac12}, ~\dbP\mbox{-a.s.},\q h=b,\si;
    \eea
for all $t\in [0, T]$, $\bx,\bx'\in \dbX$, $\xi, \xi'\in \cX_2$, $a\in A$, and $\alpha\in\cA_{[0,T]}$.

\no {\rm (iii)} {$f, g$ are uniformly H\"{o}lder-$\beta$ continuous in $\xi$ for some $0<\beta\leq 1$:}
      \beaa
   \left.\ba{c}
   { \dis |f_t( \ul\xi, \ul \a_t)-f_t( \ul\xi', \ul\a_t)| \le L  \|\xi_{\cd\wedge t}-\xi'_{\cd\wedge t}\|^{\beta}_2,\q |g(\ul\xi)-g(\ul\xi')| \le L  \|\xi-\xi'\|^{\beta}_2,}
   \ea\right.
    \eeaa
for all $t\in [0, T]$, $\xi, \xi'\in \cX_2$, and $\alpha\in\cA_{[0,T]}$.
\end{assum}

Here ${\bf 0}$ denotes the zero path  and $\ul{\bf 0}\in \cX_2$  the  zero stochastic process.  Notice that the volatility $\si$ is possibly degenerate,  so there is no loss of generality in taking it as a square matrix and considering $X$ and $B$ with the same dimension.

\begin{rem}
\label{rem-LipschitzCond}
 { The Lipschitz continuity in \reff{LipschitzCond} is under the conditional expectation, which is stronger than the following standard Lipschitz condition:
\bea
 \label{Lipschitz}
   |h_t(\bx,a, \ul\xi,\ul \a_t)-h_t( \bx,a,\ul\xi',\ul \a_t)| \le L \|\xi_{\cd\wedge t}-\xi'_{\cd\wedge t}\|_2,\q h=b,\si.
    \eea
    This is mainly to deal with the common noise case in Section \ref{sect-MFC}. For mean field control problems without common noise, it is sufficient to replace \reff{LipschitzCond} with \reff{Lipschitz}.

    We also note that, in the mean field control setting with common noise, \reff{LipschitzCond} is implied by the standard Lipschitz conditions for the data on the Wasserstein space, see Assumption \ref{assum-MFC} and Proposition \ref{prop-MFC1} below.}
    \end{rem}

\begin{lem}
\label{lem-DPP}
Let Assumption \ref{assum-standing} hold true.

\no{\rm (i)} For any $(t, \xi)\in \cX_{[0, T]}^2$ and $\a\in\cA_{[t,T]}$, SDE \reff{state1} admits a unique strong  solution $X^{t,\xi,\a}$. Moreover, for any $p\ge 2$, there exists a constant $C_p$, depending only on $p$, $T$, $d$, and the constants $L$, $C_0$ in Assumption \ref{assum-standing}, such that, for any $\xi'\in \cX_p$,
\bea
\label{Xreg}
\left.\ba{c}
\dis \|X^{t,\xi,\a}\|_p \leq C_p\big(1+\|\xi_{\cdot\wedge t}\|_p\big);\q \|X^{t,\xi,\a}-X^{t,\xi',\a}\|_p \leq C_p\|\xi_{\cdot\wedge t}-\xi_{\cdot\wedge t}'\|_p;\ms\\
\dis \|X_{\cdot\wedge s}^{t,\xi,\a}-X_{\cdot\wedge s'}^{t,\xi,\a}\|_p\leq C_p\big(1+\|\xi_{\cdot\wedge t}\|_p\big)|s-s'|^{\frac{1}{2}}, \q s,s'\in [t,T].
\ea\right.
\eea
\no{\rm (ii)} The functions $J$ and  $V$ are adapted in $\xi$ in the sense of \reff{adapted}, and $V$ satisfies the dynamic programming principle: for all $(t,\xi)\in \overline{\cX}^2_{[0,T]}$ and $\d \le T-t$,
\bea
\label{DPP}
\dis V_t( \ul\xi)
=
\inf_{\a\in \cA_{[t, T]}}
\Big\{V_{t+\d}(\ul X^{t,\xi, \a})
         + \int^{t+\delta}_{t} f _s( \ul X^{t,\xi,\a},\ul\a_s)ds\Big\}.
\eea
Moreover, denoting $\D t:= t-t'$, $\D\xi := \xi-\xi'$, we have:
\bea
\nonumber
&\dis |J_t( \ul\xi, \ul\a)|\leq C\big(1+{\|\xi_{\cdot\wedge t}\|^{\beta}_2}\big);\q |J_t(\ul\xi, \ul\a)-J_t( \ul\xi', \ul\a)|\leq C{\|\D\xi_{\cdot\wedge t}\|^{\beta}_2};
\\
\label{Vreg}
&\dis\!\!\!\!\!\!\!\!  |V_t(\ul\xi)|\leq C\big(1+{\|\xi_{\cdot\wedge t}\|^{\beta}_2}\big);~
 |V_t(\ul\xi)-V_{t'}(\ul\xi')|\leq C\Big[{{\|\xi_{\cdot\wedge t}-\xi'_{\cdot\wedge t'}\|^{\beta}_2} + \big(1+ \|\xi_{\cdot\wedge t}\|^{\beta}_2\big) |\D t|^{\frac{\beta}{2}}}\Big].
\eea
\end{lem}
The DPP \reff{DPP} follows from similar arguments as in \cite[Theorem 3.4]{cosso3}, and we shall sketch a proof in the Appendix for the reader's convenience. All the involved estimates are rather standard, in particular the conditional $\dbL^2$-type regularity of $b, \si$ with respect to $\ul X$ in Assumption \ref{assum-standing} (ii) does not induce any difficulty. We thus omit those proofs.

\section{The HJB equation and classical solutions}
\label{sect-test}

Following the standard control theory, the DPP \reff{DPP} induces an HJB equation for $V$. For this purpose, we first introduce derivatives for functions on $\ol \cX^p_{[0, T]}$.

\subsection{Smooth functions on $\ol \cX^p_{[0, T]}$}
For a generic metric space $E$ and for $t \in [0, T)$, $p\ge 2$, let $C^0(\cX_{[t, T]}^p; E)$ denote the space of adapted and continuous functions $\f: \cX_{[t, T]}^p\to E$. In particular, $C^0(\cX_{[t, T]}^p):=C^0(\cX_{[t, T]}^p; \dbR)$.
Throughout this paper, $\cM^p_{[t, T]}$ denotes the space of It\^o processes $X\in \cX_p$ such that
\bea
\label{cMp}
dX_s = \b_s ds + \g_s dB_s, \q s\in [t,T],\q\mbox{for some}~ (\b,\g)\in \dbL^p(\dbF_{[t, T]}; \dbR^d\times\dbR^{d\times d}).
\eea
Here we abuse the notation $\b$ with the H\"{o}lder continuity of $f$ and $g$.

 \begin{defn}
 \label{defn-C12cXp}
For $\hat t<T$, $p\ge 2$,  we say $\f\in C^{1,2}(\ol\cX^p_{[\hat t, T]})$ if  $\f\in C^0(\ol\cX^p_{[\hat t, T]})$ and there exist
\beaa
 \partial_t\f\in C^0\big(\ol\cX^p_{[\hat t, T]}\big),
 \q\partial_{X}\f\in C^0\big(\ol\cX^p_{[\hat t, T]}; \dbL^{\frac{p}{p-1}}(\cF;{\mathbb{R}}^d)\big),
 \q\partial_{\bx X}\f \in C^0\big(\ol\cX^p_{[\hat t, T]}; \dbL^{\frac{p}{p-2}}(\cF;\dbS^d)\big),
 \eeaa
 such that  $\partial_t\f_t( \ul \xi)$, $\partial_X\f_t( \ul \xi)$ and $\partial_{\bx X}\f_t(\ul\xi)$ are  $\cF_t$-measurable for all $(t, \xi)\in \ol\cX^p_{[\hat t, T]}$, and
 for any 
 $X\in \cM^p_{[\hat t, T]}$, the following functional It\^{o} formula holds:
\bea
\label{Ito}
{d\over dt}\f_t( \ul X)=\partial_t\f_t( \ul X)+ \dbE\Big[\partial_{X}\f_t(\ul X)\cd  \b_t+ {1\over 2} \partial_{\bx X}\f_t(\ul X) :  \g_t\g_t^\top \Big], \q t\in [\hat t, T].
\eea
\end{defn}

The following simple lemma is crucial for the above definition. The proof is similar to the justification of \cite[Definition 2.8]{ETZ1}, and is postponed to the Appendix.

\begin{lem}
\label{lem-dupire}
For any $\f\in C^{1,2}(\ol\cX^p_{[\hat t, T]})$, the derivatives $\pa_t \f, \pa_X \f, \pa_{\bx X}\f$ are unique.
\end{lem}

\begin{rem}
\label{rem-dupire}
{\rm (i)} In the state dependent case:  $\f_t(\ul\xi) = \psi(t, \ul  \xi_t)$ for some $\psi$ defined on $[0, T]\times \dbL^2(\cF_T)$,  we see that $\pa_t \f = \pa_t \psi$ is standard and  $\pa_X \f = \pa_X \psi$ is the Fr\'echet derivative:
\beaa
&\dis \pa_t \psi(t, \ul Y) = \lim_{\d\downarrow 0} {1\over \d} [\psi(t+\d, \ul Y)-\psi(t, \ul Y)],\q Y\in \dbL^2(\cF_T);\\
&\dis \psi(t, \ul Y+\ul Z)  = \psi(t, \ul Y) + \dbE\Big[ \pa_X \psi(t, \ul Y+\ul Z) Z\Big] + o(\|Z\|_2), \q Y, Z\in \dbL^2(\cF_T).
\eeaa
Moreover, $\pa_{xX}\f = \pa_{xX}\psi$ can be determined by:
\beaa
\pa_{xX}\psi(t, \ul X_t) = {d\over ds} \Big\la \pa_{X}\psi(t, \ul X_t + \ul B_\cd-\ul B_t), B_\cd\Big\ra_s \Big|_{s=t}
&\mbox{for all}&
X\in\cM^p_{[t,T]}.
\eeaa

\no {\rm (ii)} In the law invariant case: $\f_t( \ul \xi) = \check \f_t( \dbP_\xi)$ for some $\check \f$ smooth as in \cite[Theorem 2.7]{WZ},  where { $\dbP$ is the fixed probability measure on the measurable space $(\O, \cF)$,  and $\dbP_\xi\in \cP_2(\dbX)$ is the law of the process $\xi$ under $\dbP$}, we can easily see that $\f\in C^{1,2}(\ol\cX^2_{[0, T]})$ and
\beaa
\pa_t \f_t( \ul \xi) = \pa_t \check \f_t( \dbP_\xi),\q \pa_X \f_t( \ul \xi) = \pa_\mu \check \f_t( \dbP_\xi, \xi),\q \pa_{\bx X} \f_t( \ul \xi) = \pa_{\o\mu} \check \f_t( \dbP_\xi, \xi),
\eeaa
with  $\pa_\mu \check\f$ the path-dependent Lions derivative  as in \cite{WZ}, and $\pa_\o$ the Dupire vertical derivative with respect to the paths of $\xi$.

\no{\rm (iii)} For the general case, one may define the path derivatives first and then prove the It\^{o} formula \reff{Ito}, as in \cite{Dupire, cotn1} for functions on $[0, T]\times \dbX$ and in \cite{WZ, cosso3} for functions on $[0, T]\times \cP_2(\dbX)$. However, for the viscosity theory later, we will need only the  It\^{o} formula, so we follow the approach of \cite{ETZ1} by using \reff{Ito} directly to define the path derivatives.
\end{rem}

\begin{eg}
\label{eg-smooth} All smooth functions involved in this paper are of the form
\beaa
\f_t( \ul\xi) := \dbE\big[\bar\f_t\big( \xi - \hat\xi_{\cd\wedge \hat t}\big)\big],\q (t, \xi) \in \ol\cX^p_{[\hat t, T]},
&\mbox{for some fixed}&
(\hat t, \hat \xi)\in \ol\cX^p_{[0, T]},
\eeaa
where $p\ge 2$ and $\bar\f: [0, T]\times \dbX\to \dbR$ is continuous, adapted, and admits continuous pathwise Dupire's derivatives $(\pa_t \bar\f,  \pa_{\bx}\bar\f, \pa_{\bx\bx}\bar\f): [0, T]\times \dbX\longrightarrow (\dbR, \dbR^d, \dbS^d)$ such that
\bea
\label{growth}
\! |\pa_{t}\bar\f_t( \bx)| \le C\big[1+|\bx_{\cd\wedge t}|_\infty^{p}\big], |\pa_{\bx}\bar\f_t( \bx)| \le C\big[1+|\bx_{\cd\wedge t}|_\infty^{p-1}\big], |\pa_{\bx\bx}\bar\f_t( \bx)| \le C\big[1+|\bx_{\cd\wedge t}|_\infty^{p-2}\big].
\eea
Then $\f\in C^{1,2}(\ol\cX^p_{[\hat t, T]})$, and denoting $\xi' := \xi -  \hat\xi_{\cd\wedge \hat t}$, we have
\bea
\label{derivative}
\left.\ba{c}
\dis \pa_t \f_t( \ul\xi) := \dbE\big[\pa_t\bar\f_t( \xi')\big],\q  \pa_X\f_t( \ul\xi) := \pa_{\bx}\bar\f_t( \xi'),\q \pa_{\bx X}\f_t( \ul\xi) := \pa_{\bx\bx}\bar\f_t( \xi').
\ea\right.
\eea
\end{eg}
\proof For any 
$X\in \cM^p_{[\hat t, T]}$, denoting $X' := X - \hat\xi_{\cd\wedge \hat t}$, by the standard functional It\^{o} formula (cf. \cite{Dupire, cotn1}) we have
\beaa
d \bar\f_t(X')
=
\pa_t  \bar\f_t( X') dt
+ \pa_{\bx} \bar\f_t(X') \!\cd\! d X_t
+ {1\over 2} \pa_{\bx\bx} \bar\f_t(X')\!:\!d \la X\ra_t,
~t\ge \hat t.
\eeaa
 {Taking the expectation on both sides and comparing it with \reff{Ito},} we obtain \reff{derivative} immediately. In particular, the integrability of $\pa_X\f$ and $\pa_{\bx X}\f$ follows directly from \reff{growth} and the $\dbL^p$-integrability of $X$ and $\hat\xi$,  and it follows from the Burkholder-Davis-Gundy inequality that $\pa_{\bx} \bar\f_t(X') \!\cd\! \g_td B_t$ is a  true martingale.
\qed

\subsection{The HJB equation on the  process space  $\ol \cX^2_{[0, T]}$}
We now consider the following HJB equation on $\ol\cX^2_{[0,T]}$:
 \bea
 \label{HJB}
&\dis  \cL U_t(\ul\xi)
:=
\partial_t U_t(\ul\xi)
+\inf_{\a_t \in \cA_t} H_t\big( \ul\xi,\partial_{X}U_t(\ul\xi),\partial_{\bx X}U_t(\ul\xi), \ul\a_t
                                      \big)= 0,~  (t,\xi)\in\ol \cX^2_{[0,T)};\nonumber\\
& \dis \mbox{where} \
 H_t\big( \ul\xi,\ul Z, \ul \G, \ul\a_*\big)
 :=\dbE\Big[ \big(b_t(\cd)\cd Z + \frac{1}{2}\sigma\si^\top_t(\cd) : \G
                       \big)(\xi,  \a_*, \ul\xi,\ul\a_*)+f_t(\ul\xi,\ul\a_*)\Big],\  \ms\\
& \dis  (t,\xi)\in \ol\cX^2_{[0, T]},~ Z\in \dbL^2(\cF_t, \dbR^d),~ \G \in \dbL^\infty(\cF_t, \dbS^d),~  \a_*\in \cA_T.\nonumber
\eea

\begin{rem}
\label{rem-HJBinvariant}
When $b, \si$ do not depend on $(\ul X, \ul \a)$ and $g(\ul X) = \dbE\big[\bar g(X)\big]$, $f_t( \ul X, \ul \a_t) = \dbE\big[\bar f_t( X, \a_t)\big]$, for some deterministic function $\bar f, \bar g$, one can easily show that
\bea
\label{V=v}
V_t( \ul \xi) = \dbE\big[v_t( \xi)\big],
\eea
where $v$ is associated with the standard path dependent HJB equation on $[0, T]\times \dbX$:
\bea
 \label{HJBstandard}
\left.\ba{c}
\dis  \partial_t v_t(\bx)+\inf_{a\in A} \Big[b_t(\bx,  a)\cd \partial_{\bx} v_t(\bx) + \frac{1}{2}\si\si^\top_t(\bx,  a) : \partial_{\bx\bx}v_t(\bx)+\bar f_t( \bx, a)\Big]=0,
 \ea\right.
\eea
and $\pa_t v, \pa_{\bx} v, \pa_{\bx\bx} v$ are Dupire's path derivatives.
Moreover, in the state dependent case, \reff{HJBstandard} reduces to the standard HJB equation:
\bea
 \label{HJBstate}
\left.\ba{c}
\dis  \partial_t v(t, x)+\inf_{a\in A} \Big[b(t, x,  a)\cd \partial_x v(t, x) + \frac{1}{2}\si\si^\top(t, x,  a) : \partial_{xx}v(t, x)+\bar f(t, x, a)\Big]=0.
 \ea\right.
\eea
\end{rem}

As usual, we start with classical solutions.

\begin{defn}
\label{defn-classical} (Classical solution)
 A functional $U\in C^{1,2}(\ol\cX^2_{[0, T]})$   is called a classical solution (resp. subsolution, supersolution) to the equation \reff{HJB} if
 \beaa
 -\cL U_t(\ul\xi) = ~(\mbox{resp.} \le,~ \ge)~ 0,\q \forall (t,\xi)\in \ol\cX^2_{[0,T)}.
 \eeaa
 \end{defn}

In the contexts of classical solutions, we also need the regularity in $t$ which is not included in Assumption \ref{assum-standing}. However, this regularity will not be needed for the viscosity solution.

\begin{assum}
\label{assum-regt}
The functions $h=b, \si, f$ are  locally uniformly continuous in $(t, a, \a_*)\in [0, T]\times A\times \cA_T$ in the following sense: for any $R>0$, there exists a modulus of continuity function $\rho_R$ such that, for all $\bx\in \dbX$,  $\|\xi_{\cd\wedge t}\|_2 \le R$, and all $ t\le t'$, $a, a'\in A$, $\a_*, \a'_*\in \cA_T$,
\beaa
\big| h_t(\bx, a, \ul\xi, \ul\a_*) - h_{t'}(\bx_{\cd\wedge t},a', \ul\xi_{\cd\wedge t}, \ul\a'_* )\big| \le \rho_R\big(|t-t'| + |a-a'| + \|\a_*-\a'_*\|_2\big).
\eeaa
\end{assum}

\begin{thm}
\label{thm-classical}
Let Assumptions \ref{assum-standing} and \ref{assum-regt} hold, and assume that the value function $V\in C^{1,2}(\ol\cX^2_{[0, T]})$. Then $V$ is the unique classical solution of the HJB equation \reff{HJB}.
\end{thm}
Combining the DPP \reff{DPP} and the It\^o formula \reff{Ito}, it is rather standard to verify that $V$ satisfies \reff{HJB}. The uniqueness is not hard either, however, since it will be a consequence of the uniqueness of the viscosity solution later, we omit the proof.

\section{The singular component of test functions}
\label{sect-phi}
In general one can hardly expect \reff{HJB} to have a classical solution. Our goal of this paper is to propose  an appropriate notion of viscosity solutions which allows for a complete characterization of the value function through the corresponding HJB equation. The main feature of our approach is that, besides the standard smooth functions $\f\in C^{1,2}(\ol\cX^2_{[0, T]})$, our test functions contain an additional component which is singular in certain sense. To motivate this, we first consider a simple  setting with constant volatility.

\subsection{A motivating case}
\label{sect-firstorder}
In this subsection we study heuristically the state dependent case with  constant volatility $\si \equiv I_{d\times d}$. In this case, it is natural to consider the following change of variables to convert the SDE into a random ODE:
\bea
\label{tildeX}
\left.\ba{c}
\dis \tilde b_s(x, \o, a, \ul {X_s}, \ul \a_s) := b_s(x + B_s(\o), \o, a, \ul X_s+\ul B_s, \ul\a_s),\ms\\
\dis \tilde f_s(\ul X_s, \ul \a_s) := f_s(\ul X_s+\ul B_s, \ul \a_s),\q \tilde g(\ul X_T) := g(\ul X_T+\ul B_T).
\ea\right.
\eea
Then, recalling \reff{state1} and \reff{value1}, we have
\bea
\label{tildeV}
\left.\ba{c}
\dis V_t(\ul\xi_t) = \tilde V_t(\ul\xi_t-\ul B_t),
\q\mbox{where}\q
\tilde V_t(\ul{\tilde\xi}_t)
:=
\inf_{\a\in \cA_{[t,T]}}
\tilde g({\ul {\tilde X}_T^{t,\tilde\xi_t,\a}}) +  \int_{t}^{T}\tilde f_s( \ul {\tilde X}_s^{t,\tilde\xi_t,\a}, \ul\a_s)ds,
\\
\mbox{and}\q
\dis \tilde X^{t,\tilde\xi_t,\a}_s= \tilde\xi_t + \int_t^s \tilde b_r\big(\tilde X^{t,  \tilde\xi_t,\a}_r,\a_r, \ul {\tilde X}_r^{t, \tilde\xi_t,\a}, \ul \a_s\big)ds,~ s\in [t, T].
\ea\right.
\eea
This is a deterministic control problem, and thus \reff{HJB} becomes a first order equation:
\bea
 \label{HJBtilde}
&\dis  \tilde \cL \tilde U_t(\ul{\xi}_t):= \partial_t \tilde U_t(\ul{\xi}_t)+\inf_{\a_t \in \cA_t} \dbE\Big[ \tilde b_t(\xi_t,  \a_t, \ul\xi_t,\ul\a_t)\cd \pa_X \tilde U_t(\ul\xi_t) +\tilde f_t(\ul\xi_t,\ul\a_t)\Big]= 0.
\eea
Since \reff{HJBtilde} does not involve $\pa_{xX} \tilde U$, it suffices to consider $\tilde U\in C^{1,1}([0, T] \times \dbL^2(\cF_T))$. That is, there exist appropriate $\pa_t \tilde U$ and $\pa_X\tilde U$ satisfying the chain rule:
\bea
\label{chainrule}
d\tilde U_t( \ul X_t)=\partial_t\f_t( \ul X_t)dt+ \dbE\Big[\partial_{X}\tilde U_t(\ul X_t)\cd \b_t dt\Big],~\mbox{where}~ dX_t = \b_t dt,
~\mbox{and}~\beta\in\dbL^p(\dbF,\dbR^d).
\eea

We now investigate the comparison principle for appropriately defined viscosity solution for the first order equation \reff{HJBtilde}. Let $\tilde U_1, \tilde U_2$ be a viscosity subsolution and supersolution, respectively.  {As standard for viscosity solutions,} 
we consider the doubling variable approach:
\bea
\label{tildepsi}
\tilde \psi(t, \ul \xi_t, s, \ul \zeta_s) := \tilde U_1(t, \ul\xi_t) - \tilde U_2(s, \ul\zeta_s) - n\Big[|t-s|^2 + \dbE\big[|\xi_t- \zeta_s|^2\big]\Big] - \cds,
\eea
where ``$\cds$'' indicates appropriate further penalty functions so as to guarantee the existence of a maximizer, denoted as
$(\hat t, \hat s, \hat \xi_{\hat t}, \hat \zeta_{\hat s})$. Then we shall consider a test function of $\tilde U_1$ in the form (omitting possibly additional terms):
\bea
\label{tildef}
\tilde \f_1(t, \ul \xi_t) := n\Big[|t-\hat s|^2 + \dbE\big[|\xi_t- \hat\zeta_{\hat s}|^2\big]\Big].
\eea
Indeed, by this argument one can easily prove rigorously the comparison principle for \reff{HJBtilde}.

We now turn back to \reff{HJB}. Recall \reff{tildeV}, it is natural to consider $U_1(t, \ul \xi_t) := \tilde U_1(t, \ul\xi_t-\ul B_t)$ as a viscosity subsolution of \reff{HJB},  with a test function induced by \eqref{tildef}:
\bea
\label{tildetest}
\f_1(t, \ul \xi_t) = \tilde \f_1(t, \ul\xi_t-\ul B_t) = n\Big[|t-\hat s|^2 + \dbE\big[|\xi_t - B_t - \hat\zeta_{\hat s}|^2\big]\Big].
\eea
The above $\f_1$ has two features:
\begin{itemize}
\item[{\rm (i)}] $\f_1$ is not law invariant, namely $\dbP_{\xi_t} = \dbP_{\xi'_t}$ does not imply $\f_1(t, \ul\xi_t) = \f_1(t, \ul \xi_t')$. This is another motivation for us to consider functions on the process space $\ol\cX^2_{[0, T]}$ directly, instead of on the Wasserstein space of probability measures.
\item[{\rm (ii)}]  $\f_1$ is not in $C^{1,2}(\ol\cX^2_{[0, T]})$ in the sense of Definition \ref{defn-C12cXp}. Indeed, assume $t\ge \hat t \vee \hat s$ and $dX_s = \b_s ds + \g_s dB_s$, $s\ge t$, then
\bea
\label{f1abs}
{d\over ds} \f_1(s, \ul X_s)
\!=\!
2n (s-\hat s) + 2n \dbE\big[(X_s -B_s - \hat \zeta_{\hat s}) \!\cd\! \b_s + {1\over 2} (\g_s - I_{d\times d})\!:\!(\g_s - I_{d\times d})^\top\big].
\eea
Recalling \reff{Ito}, we see that $\pa_{\bx X}\f_1$ does not exist and thus $\f_1 \notin C^{1,2}(\ol \cX^2_{[0, T]})$. However, from \reff{f1abs} it is clear that $s\mapsto \f_1(s, \ul X_s)$ is absolutely continuous. This motivates us to consider test functions which are absolutely continuous in $t$ but are not in $C^{1,2}(\ol \cX^2_{[0, T]})$.
\end{itemize}

\subsection{The singular component of test functions in the general case}
We now introduce the singular component of test functions, denoted as $\phi$, which is absolutely continuous in $t$ but is not in $C^{1,2}(\ol \cX^2_{[0, T]})$ in general. This part $\phi$ is new in the  literature of mean field control problems and is the main feature of our notion of viscosity solutions due to its crucial role in our proof of the comparison principle.

 For $\tilde t\in [0, T]$ and $p\ge 2$, let $\Xi^p_{\tilde t}$ denote the set of maps $(\tilde b, \tilde \si, \tilde f)$ where $(\tilde{b},\tilde{\sigma}): [\tilde t, T]\times \O\times  A \times\cA_{T}
 \to {(\dbR^d, \dbR^{d\times d})}$ are $\dbF_{[\tilde t, T]}$-progressively measurable, $\tilde f: [\tilde t, T]\times\cA_{T} \to \dbR$ is progressively measurable, and there exists a random variable $\G_*\in \dbL^p(\cF_T)$ such that
 \bea
 \label{tildebsifbound}
\sup_{[\tilde t, T]\times A\times \cA_T} [|\tilde b|+|\tilde \si|]\le \G_*, ~\dbP\mbox{-a.s.},\qq    \sup_{[\tilde t, T]\times\cA_T}|\tilde f|<\infty.
\eea
   Given $(\tilde t, \tilde\xi)\in \ol\cX^p_{[0,T]}$ and $(\tilde{b},\tilde{\sigma}, \tilde f)\in \Xi^p_{\tilde t}$, we introduce the maps defined for all $(t,\xi)\in \ol\cX^p_{[\tilde t,T]}$:
 \bea
\label{cIaxi}
\left.\ba{c}
 \dis \cI^{\a}_t(\xi)
 :=
 \cI^{\tilde{b},\tilde{\si},\tilde t, \tilde\xi,\a}_t(\xi)
 :=
 \xi_t
 - \tilde\xi_{\tilde t}
- \int_{\tilde t}^t \tilde{b}^{\a}_s ds -  \int_{\tilde t}^t \tilde{\si}^{\a}_s dB_s,\q
F^{\a}_t
:=
F^{\tilde f,\tilde t,\a}_t  := \int_{\tilde t}^t \tilde f^{\a}_s ds,\ms\\
\dis\mbox{where}\q \tilde{b}^{\a}_s:=\tilde b_s(\alpha_s,\ul\alpha_s)$, $\tilde{\si}^{\a}_s:=\tilde\si_s(\alpha_s,\ul\alpha_s)$, $\tilde{f}^{\a}_s:=\tilde f_s(\ul\alpha_s).
\ea\right.
\eea
We now introduce the class of singular test functions inspired from \eqref{tildef} which allows us to handle the general setting of non-constant and even controlled volatilities.

 \begin{defn}
 \label{defn-C+cXp}
 For $\tilde t\in [0, T]$ and $p\ge 2$, we denote $C^+(\ol\cX^p_{[\tilde t,T]})$ the set of maps of the form:
 \bea
 \label{C+cXp}
 \dis \phi_t( \ul\xi)
 :=
 \inf_{\a\in \cA_{[\tilde t, T]}}\Big\{
 k \dbE\Big[ \big|\cI^{\a}_t(\xi) \big|^p
                + \big|\cI^{\a}_{t'}(\xi')
                 \big|^p
          \Big]
+ \int_{t'}^t \tilde f^{\a}_s ds\Big\},
\!\!&\mbox{for all}& \!\!
(t, \xi) \in\ol \cX^p_{[\tilde t, T]},
 \eea
for some $\tilde \xi\in \cX_p$, $(t',\xi')\in \ol\cX^p_{[\tilde t, T]}$, $(\tilde{b},\tilde{\sigma}, \tilde f)\in \Xi^p_{\tilde t}$,  and some constant $k\geq0$.

Moreover, let $C^-(\ol\cX^p_{[\tilde t, T]})$ denote the set of $\phi$ such that $-\phi \in C^+(\ol\cX^p_{[\tilde t, T]})$.
\end{defn}

\begin{rem}
\label{rem-C+cXp}
{\rm (i)} Due to the involvement of $\tilde \xi, \xi'$ in \reff{C+cXp}, $\phi$ is not law invariant. This explains partially that we need to work on $\cX_p$ instead of the Wasserstein space $\cP_p(\dbX)$.

\no {\rm (ii)} All the results in this paper will remain true, after obvious modifications, if we restrict the $(\tilde b, \tilde \si, \tilde f)$ in \reff{C+cXp} to $(b, \si)(\hat \xi_{\cd\wedge \tilde t}, \cd, \ul{\hat\xi}_{\cd\wedge \tilde t}, \cd)$ and $f(\ul{\hat\xi}_{\cd\wedge \tilde t}, \cd)$ for some $\hat\xi\in \cX_p$, by using the true data $(b, \si, f)$. We allow for the flexibility on $(\tilde b, \tilde \si, \tilde f)$ so that our definition of viscosity solutions is model independent.

\no {\rm (iii)} When the volatility is not controlled:  $\si = \si_t( \xi,\ul\xi)$, we may consider the simpler set $C^+(\ol\cX^p_{[\tilde{t},T]})$ consisting of functions
\bea
\label{phidrift}
\phi_t(\xi) := \dbE\big[\big|\xi_t - \tilde\xi_{\tilde{t}} - \tilde\si_{\tilde{t}} (B_t - B_{\tilde{t}})\big|^p\big], \q\mbox{for some} ~\tilde \xi_{\tilde{t}}\in \dbL^p(\cF_{\tilde{t}}; \dbR^d),~\tilde\si_{\tilde{t}}\in \dbL^p(\cF_{\tilde{t}}; \dbS^d),
\eea
or  even restricting further to $\tilde\si_{\tilde{t}} =  \si_{\tilde t}(\hat \xi_{\cd\wedge \tilde t},  \ul{\hat\xi}_{\cd\wedge \tilde t})$  for some $\hat\xi\in \cX_p$, by using the true data $\si$ as in (ii). In these cases the related arguments can be simplified significantly.

It is clear that the term $\tilde \si_{\tilde t}$ is designed to cancel the impact of the diffusion term $\si_t(X, \ul {X})dB_t$ in the related estimates. We note that, in the drift control case, there is no need to involve $\tilde b, \tilde f$ in \reff{phidrift}. However, in the volatility control case, the estimates become much more involved. Then, besides a more carefully designed $\tilde \si^\a$, the terms $\tilde b^\a, \tilde f^\a$ in \reff{cIaxi} are also crucial for the estimates later.

\no {\rm (iv)} For the estimates later, we will often use the following equivalent formulation of $\phi$:
\bea
 \label{C+cXp2}
 \dis \phi_t( \ul\xi)
 :=
 \inf_{\a\in \cA_{[\tilde t, T]}}\Big[
 k \dbE\big[ \big|\cI^{\a}_t(\xi) \big|^p \big] + F^\a_t + \k(\a)\Big],~\mbox{where}~ \k(\a):= k \dbE\big[ \big|\cI^{\a}_{t'}(\xi') \big|^p \big] - F^\a_{t'}.
 \eea
\end{rem}

The next result states that the function $t\mapsto\phi_t$ is absolutely continuous with respect to the Lebesgue measure, and provides some crucial estimates for the comparison principle of viscosity solutions. We observe that the   {time derivative $\dot\phi$} of $\phi$ is not continuous in general.
\begin{prop}
\label{prop-C+cXp}
Let $\phi\in C^+(\ol\cX^p_{[\tilde t,T]})$ as in \eqref{C+cXp} with corresponding $(t',\xi')$ and $p\ge 2$. Then, for any $X\in \cM^p_{[\tilde t, T]}$, the mapping $t\mapsto \phi_t( \ul X)$ is absolutely continuous, with  time derivative $\dot{\phi}$ satisfying:

\no {\rm (i)} when $t\ge t'$,
\begin{eqnarray}
\label{pa+phi-est2}
\int_t^{t+\d} \dot\phi_s( \ul X)ds
&\le&
 \inf_{\a\in \cA_{[t, T]}}
\int_t^{t+\d} \big[{p(p-1)\over 2}kI^p_s( \a)+ \tilde{f}^{\a}_s\big]ds,
\end{eqnarray}

 \vspace{-8mm}$$
 \mbox{where}~
 I^p_s(\a)
 :=
 \big\|\b_s - \tilde{b}^{\a}_s\big\|_p
 \sup_{ \a'\in \cA_{[\tilde t, T]}}
 \big\|\cI^{\a'}_s(X)\big\|_p^{p-1}
 + \big\|\gamma_s - \tilde{\si}^{\tilde\a}_s\big\|_p^2
 \sup_{ \a'\in \cA_{[\tilde t, T]}}
 \big\|\cI^{\a'}_s(X)\big\|_p^{p-2};
 $$

\no {\rm (ii)} alternatively, when $t< t'$, for any $\d>0$, there exists $\a^\d\in \cA_{[\tilde t, T]}\subset \cA_{[t, T]}
$, which may depend on $X_{\cd\wedge t}$ but is independent of $(\b_s, \g_s)_{s\ge t}$, such that
\begin{eqnarray}
\int_t^{t+\d}\dot{\phi}_s( \ul X)ds
&\le&
 \int_t^{t+\d} \Big({p(p-1)\over 2}kI^p_s(\a^\d) +  \tilde{f}^{\a^\d}_s\Big)ds+\d^2.
\label{pa+phi-est1}
\end{eqnarray}
\end{prop}
The proof is postponed to  { the Appendix}.

\section{Viscosity solutions of the HJB  equation}
\label{sect-HJB}

In this section we propose a notion of viscosity solutions for \reff{HJB}. 
 For any $p\geq2$, $U: \ol\cX^2_{[0, T]}\to \dbR$ and $(t, \xi) \in\ol \cX^{p}_{[0,T]}$, denote
\bea
\nonumber
{\mathfrak{F}^+_{p}}U_t( \ul\xi)
&:=&
\Big\{(\f, \phi)\in C^{1,2}(\ol\cX^p_{[t,T]})\times C^+\ol\cX^p_{[\tilde t,T]})\ \mbox{for some}\ \tilde{t}\in [0, t]:
\\
&&\hspace{10mm}
 \big[U-(\varphi+\phi)\big]_t(\ul\xi)
=
\sup_{(s,\zeta)\in \ol \cX^p_{[t,T]}} \big[U-(\varphi+\phi)\big]_s(\ul\zeta)\Big\};
\label{cAtest+}
\\
{\mathfrak{F}^-_{p}}U_t(\ul\xi)
&:=&
\Big{\{}(\f, \phi)\in C^{1,2}(\ol\cX^p_{[t,T]})\times C^-(\ol\cX^p_{[\tilde t,T]})\ \mbox{for some}\ \tilde{t}\in [0, t]:\nonumber\\
&&\hspace{10mm}
 \big[U - (\varphi+\phi)\big]_t(\ul\xi)
=
\inf_{(s,\zeta)\in \ol\cX^p_{[t,T]}} \big[U-(\varphi+\phi)\big]_s(\ul\zeta)\Big\}.
\label{cAtest-}
\eea

Due to the use of the singular test functions $\phi$, we need to introduce the frozen state process defined for all $(t,\xi)\in \ol\cX^2_{[0, T]}$ and $\a\in \cA_{[t, T]}$ by:
\bea
\label{barXtxi}
\left.\ba{c}
\dis
\bar X^{t, \xi, \a}_s := \xi_s,~ s\in [0, t];\q \bar X^{t, \xi, \a}_s
:=
\xi_{t}
+
\int_{t}^s b^{t,\xi, \a}_r dr +
\int_{t}^s \si^{t,\xi,\a}_r dB_r, ~s\in [t, T],
\\
\mbox{where}~
 h^{t,\xi, \a}_s
 := h_s(\xi_{\cdot\wedge t}, \a_s, \ul {\xi}_{\cdot\wedge t},\ul \a_s),
 ~\mbox{for}~
 h = b, \si, \q\mbox{and}~ f^{t,\xi, \a}_s
 := f_s(\ul {\xi}_{\cdot\wedge t},\ul \a_s),\q s\ge t.
\ea\right.
 \eea
We notice the slight difference between $\bar X$ and the $X$ in  \reff{state1}. In particular, in $h^{t, \xi, \a}_s$ the state process is frozen while the control part $\a_s$ is evolving in $s$.

\begin{defn}
\label{defn-viscosity}
 {\rm (i)} For any $p\geq 2$, $U\in USC(\ol\cX^2_{[0, T]})$ is a viscosity {$p$-subsolution} of HJB equation (\ref{HJB}) if, for all {$(t,\xi)\in \ol\cX^{p}_{[0, T)}$ and $(\varphi, \phi)\in \mathfrak{F}^+_{p}U_t(\ul\xi)$,}
 \bea
 \label{sub}
\!\!\!\!\!\! \partial_t \f_t(\ul\xi)+\liminf_{\d\to 0}\inf_{\a \in \cA_{[t, T]}} {1\over \d}\int_t^{t+\d} \big[ H_s( \ul\xi_{\cd\wedge t}, \partial_X\f_t(\ul\xi),\partial_{\bx X}\f_t(\ul\xi), \ul \a_s)  + \dot \phi_s( \ul {\bar X}^{t,\xi, \a})\big]ds \ge 0.
  \eea

 \no{\rm (ii)} For any $p\ge 2$, $U\in LSC(\ol\cX^2_{[0, T]})$ is a viscosity {$p$-supersolution} of  HJB equation (\ref{HJB}) if, for all {$(t,\xi)\in \ol\cX^{p}_{[0, T)}$ and $(\varphi, \phi)\in  \mathfrak{F}^-_{p}U_t(\ul\xi)$,}
\bea
 \label{super}
\!\!\!\!\!\! \partial_t \f_t(\xi)+\limsup_{\d\to 0}\inf_{\a \in \cA_{[t, T]}} {1\over \d}\int_t^{t+\d}\big[H_s( \ul\xi_{\cd\wedge t}, \partial_X\f_t(\ul\xi),\partial_{\bx X}\f_t(\ul\xi), \ul \a_s)  + \dot \phi_s( \ul {\bar X}^{t,\xi, \a})\big]ds \le 0.
  \eea
  \no{\rm (iii)} {$U$ is a viscosity subsolution (resp., supersolution) of HJB equation (\ref{HJB}) if there exists $p_0\ge 2$ such that $U$ is a viscosity $p$-subsolution (resp., $p$-supersolution) of HJB equation (\ref{HJB}) for all $p\geq p_0$.}

 \no{\rm (iv)}  $U\in C^0(\ol\cX^2_{[0, T]})$ is a viscosity solution of  HJB equation (\ref{HJB}) if it is both a viscosity subsolution and a viscosity supersolution of \reff{HJB}.
\end{defn}

\begin{rem}
\label{rem-viscosity1}
{\rm (i)} While in quite different forms, our idea of introducing the singular component $\phi$ is inspired by works on viscosity solutions for PDEs in infinite dimensions, see, e.g., Crandall-Lions \cite[Definition 2.1]{cran1}, Li-Yong \cite[Chapter 6 Definition 3.1]{LY}, and Fabbri-Gozzi-{Swiech} \cite[Definitions 3.32 $\&$ 3.35]{fab1}.  { A similar idea has also been used to study viscosity solutions for parabolic PDEs whose coefficients are discontinuous in time, see, e.g., Ishii \cite{Ishii}, Lions-Perthame \cite{LP}, and Nunziante \cite{Nunziante1, Nunziante2}.}

\no {\rm (ii)} The standard definition of viscosity solutions amounts to setting $\phi=0$. So a viscosity solution in our sense is always a viscosity solution in the standard sense. Consequently, by introducing the component $\phi$, the existence of viscosity solutions becomes slightly harder, but as we will see it significantly helps for the comparison principle of viscosity solutions.

\no {\rm (iii)} We take the integral form in the left side of \reff{sub} and \reff{super} because $\dot\phi$ is discontinuous, in general. If it were continuous,  then these expressions would reduce to a simpler and more standard form under Assumption \ref{assum-regt}:
\beaa
\partial_t \f_t(\ul\xi)+\inf_{\a \in \cA_{[t, T]}} \big[H_t( \ul\xi, \partial_X\f_t(\ul\xi),\partial_{\bx X}\f_t(\ul\xi), \ul \a_t)  + \dot \phi_t( \ul {\bar X}^{t,\xi, \a})\big].
  \eeaa
\end{rem}

\begin{rem}
\label{rem-viscosity2}
{\rm (i)} { For $2\le p_1<p_2$, $\mathfrak{F}^+_{p_1}U_t( \ul\xi)$ and $\mathfrak{F}^+_{p_2}U_t( \ul\xi)$ do not contain each other. Consequently, the viscosity $p_1$-subsolution property and the viscosity $p_2$-subsolution property do not imply each other.}

\no {\rm (ii)} In light of \reff{Upsilon0} below, it is convenient to choose $p$ as an even integer.
 We shall prove the comparison  principle only on $\ol\cX^{p}_{[0, T]}$ with some large even integer $p$. Since  $\cX_{p}$ is dense in $\cX_2$, then, when $U$ is continuous under $\|\cd\|_2$, we obtain the comparison principle on $\ol\cX^2_{[0, T]}$.

\no {\rm (iii)} We should note that $\cX_p$ is not compact under $\|\cdot\|_p$. We shall circumvent this difficulty by using the Ekeland-Borwein-Preiss variational principle.
\end{rem}

We first show that our notion of viscosity solutions is consistent with that of classical solutions.

\begin{prop}\label{prop-consistent}
 Let Assumptions \ref{assum-standing} and \ref{assum-regt} hold true and $U \in C^{1,2}(\ol\cX^2_{[0, T]})$. Then $U$ is a viscosity subsolution (resp. supersolution) of \reff{HJB} if and only if it is a classical subsolution (resp. supersolution) of \reff{HJB}.
\end{prop}

The proof of this result is rather standard, and thus is postponed to the Appendix.
Our main result of the paper is the following characterization of the value function by means of the corresponding HJB equation \eqref{HJB}.

\begin{thm}\label{thm-exist}
Under Assumption \ref{assum-standing}, the value function $V$ is the unique viscosity solution of the equation (\ref{HJB}) with terminal condition $V_T = g$ in the class of functions satisfying \reff{Vreg}.
\end{thm}

The viscosity property of the value function $V$ follows from standard arguments and will be reported below. The uniqueness is as usual more challenging, and is a consequence of the following comparison result.

\begin{thm}\label{thm-comparison}
Let Assumption \ref{assum-standing} hold true, and  $U_1\in USC(\ol\cX^2_{[0, T]})$, $U_2\in LSC(\ol\cX^2_{[0, T]})$ be a viscosity subsolution and supersolution of HJB equation  \reff{HJB},  respectively.  Assume further that there exists a modulus of continuity function $\rho_R$ for each $R>0$ such that

\no{\rm (i)} one of $h=U_1$ or $U_2$ satisfies the following estimate slightly weaker than \reff{Vreg},
\bea
\label{Vreg1}
|h_t(\ul\xi)-h_{t'}(\ul\xi')|\leq C\Big[\|\D \xi_{\cdot\wedge t}\|^{\beta}_2 + \big(1+ \|\xi_{\cdot\wedge t}\|_2\big) |\D t|^{\frac{\beta}{2}}\Big].
\eea

\no{\rm (ii)} and the other one satisfies, for $h=U_1$ or $-U_2$ and for any $R>0$:
\bea
\label{Ureg}
\dis |h_t(\ul \xi)|\le C\big(1+\|\xi_{\cd\wedge t}\|_2\big),
~~h_t(\ul\xi) - h_s(\ul\xi_{\cd\wedge t}) \le \rho_R(s-t),
~\mbox{for all}~
t<s, \|\xi_{\cd\wedge t}\|_2\le R.
\eea
Then $U_1(T,\cd) \le U_2(T, \cd)$, on $\cX_2$, implies that $U_1\le U_2$ on $\ol\cX^2_{[0, T]}$.
\end{thm}

We defer the proof of the last theorem to Sections \ref{sect-comparison1}, \ref{sect-comparison2}, and \ref{sect-comparison3}. This theorem also implies the following comparison result immediately.

\begin{thm}\label{thm-comparison2}
Let Assumption \ref{assum-standing} hold true, and  $U_1\in USC(\ol\cX^2_{[0, T]})$, $U_2\in LSC(\ol\cX^2_{[0, T]})$ be a viscosity subsolution and supersolution of \reff{HJB},  respectively.  Assume $U_1, U_2$ satisfy the estimates \reff{Ureg}, and  $U_1(T,\cd) \le g\le U_2(T, \cd)$ on $\cX_2$.  Then $U_1\le U_2$ on $\ol\cX^2_{[0, T]}$.
\end{thm}
\proof Note that the value function $V$ is a viscosity solution of \reff{HJB} with terminal condition $V_T = g$ and satisfies \reff{Vreg}, and hence \reff{Vreg1}. Then, by applying Theorem \ref{thm-comparison} on $U_1$ and $V$ we have $U_1 \le V$, and by applying Theorem \ref{thm-comparison} on $V$ and $U_2$ we have $V \le U_2$. Thus $U_1 \le U_2$.
\qed

\begin{rem}
\label{rem-Ishii}
{\rm (i)} We remark that we will not use the Crandall-Ishii lemma in the proof of Theorem \ref{thm-comparison}, despite that our control problem involves the diffusion term. This is not completely surprising because, as we explained in Subsection \ref{sect-firstorder}, the singular component $\phi$ of the test function is essentially involved to cancel the diffusion term.

\no {\rm (ii)} In the setting of Remark \ref{rem-HJBinvariant}, our HJB equation \reff{HJB} reduces to standard equations \reff{HJBstandard} or \reff{HJBstate}. Defining $v$ as a viscosity solution of \reff{HJBstandard} or \reff{HJBstate} if the function $V$ defined in \reff{V=v} is a viscosity solution of the equation \reff{HJB}, then we obtain the comparison principle for the fully nonlinear second order (path dependent) PDE by using the doubling variable arguments but without using the  Crandall-Ishii lemma.

\no {\rm (iii)} However, we shall emphasize that the above definition of viscosity solutions is not equivalent to the standard notion of viscosity solutions for HJB equations, e.g. in \cite{cran2}. So we are not claiming that we can avoid the Crandall-Ishii lemma  {in the doubling variable arguments} for the standard viscosity solutions of fully nonlinear second order HJB equations.

\no {\rm (iv)} We also refer to Remark \ref{rem-Ishii2} below for a highly related comment.
\end{rem}

\ms
\no {\bf Proof of Theorem \ref{thm-exist}.} As the uniqueness is implied by the comparison result of Theorem \ref{thm-comparison}, we only focus on the existence part.  {Recall Definition \ref{defn-viscosity} and we fix arbitrary $p\ge 4$.}

(i) We first prove the viscosity $p$-subsolution property. Fix  $(t, \xi) \in \ol\cX^p_{[0, T)}$ and $(\varphi, \phi)\in \mathfrak{F}^+_{p}V_t( \ul\xi)$. For any $\a\in \cA_{[t, T]}$, recall \reff{state1} and denote $X^\a:= X^{t,\xi,\a}$. By DPP \reff{DPP} and \reff{cAtest+} we have for any $\d>0$,
\beaa
\dis  0 &\le& V_{t+\d}( \ul X^{\a}) - V_t( \ul\xi)+ \int^{t+\delta}_{t} f _s(\ul X^{\a},  \ul\a_s)ds\\
\dis &\le& [\f+\phi]_{t+\d}( \ul X^{\a}) - [\f+ \phi]_t( \ul\xi) + \int^{t+\delta}_{t}f _s(\ul X^{\a},  \ul \a_s)ds.
 \eeaa
Applying the It\^o formula \reff{Ito} on $\f_s(  \ul X^{\a})$ we have
 \beaa
\dis 0 &\le & {1\over \d} \int_t^{t+\d} \dbE\Big[\pa_t\f_s( \ul X^{\a}) + \pa_X\f_s(\ul X^\a) \cd b_s( X^\a, \a_s, \ul X^\a, \ul \a_s) \\
&& \hspace{18mm}
+ {1\over 2} \pa_{\bx X}\f_s( \ul X^{\a}) : \si\si^\top_s( X^\a, \a_s, \ul X^\a, \ul \a_s)+f _s(\ul X^{\a},  \ul \a_s) + \dot\phi_s( \ul X^\a)\Big]ds.
 \eeaa
 Recall \reff{barXtxi} and denote $\bar X^\a:= \bar X^{t,\xi,\a}$.
 We claim that
 \bea
 \label{existence-claim}
 \sup_{\a\in \cA_{[t, T]}} \D^\a_\d  = o(\d),\q \mbox{where}\q \D^\a_\d := \Big| \int_t^{t+\d}  \big[ \dot\phi_s( \ul X^\a) - \dot\phi_s( \ul{\bar X}^\a)\big]ds\Big|.
 \eea
 Together with the regularity of $\f$ and $b, \si, f$, and recalling the {$h^{t,\xi,\a}$} in \reff{barXtxi}, this implies:
\beaa
&&\dis 0 \le \inf_{\a\in \cA_{[t, T]}}{1\over \d} \int_t^{t+\d}
\dbE\Big[\pa_t\f_t( \ul\xi) + \pa_X\f_t( \ul\xi) \cd b^{t, \xi, \a}_{s}\\
&&\dis\hspace{40mm}  + {1\over 2} \pa_{\bx X}\f_t(\ul\xi) : (\si\si^\top)^{t,\xi, \a}_{s} +f^{t, \xi, \a}_{s}  + \dot\phi_s( \ul{\bar X}^\a)\Big]ds + o(1).
\eeaa
Send $\d\to 0$, we obtain \reff{sub}, namely the viscosity subsolution property of $V$.

We now prove \reff{existence-claim}. Recall \reff{cIaxi} and \reff{C+cXp2} with corresponding $\tilde t, \tilde \xi, t', \xi'$, $\tilde b, \tilde \si, \tilde f$ and $p\geq 4$. By \reff{tildebsifbound} we denote
 \bea
 \label{cp}
c_p := \Big(\dbE[|\G_*|^p]\Big)^{1\over p} + \sup_{[\tilde t, T]\times\cA_T}|\tilde f|<\infty.
\eea
 Note that $X^\a_{\cd\wedge t} = \xi_{\cd\wedge t} = \ol X^\a_{\cd\wedge t}$, then $\phi_t(\ul X^\a) =  \phi_t( \ul{\bar X}^\a)$. Thus, by \reff{C+cXp2},
\bea
\label{phit+d}
\left.\ba{c}
\dis \D^\a_\d = \Big| \phi_{t+\d}( \ul X^\a) - \phi_{t+\d}( \ul{\bar X}^\a)\Big| \le k \sup_{\tilde\a\in \cA_{[\hat{t}, T]}} \D^{\a,\tilde \a}_\d, \\\dis\mbox{where}\q \D^{\a,\tilde \a}_\d:=\Big|\dbE\big[\big|\cI^{\tilde \a}_{t+\d}( X^\a)\big|^p-\big|\cI^{\tilde \a}_{t+\d}( \bar X^\a)\big|^p\big]\Big|.
\ea\right.
\eea
Here we use $\tilde \a$ in \reff{C+cXp2} so as to distinguish with the $\a$ in  \reff{existence-claim}.
Fix $\a, \tilde \a$ and denote
\beaa
h^{\a,\tilde \a}_s := h_s( X^\a, \a_s, \ul X^\a, \ul\a_s) - \tilde{h}^{\tilde \a}_s,\q \ol h^{\tilde \a}_s := h^{t,\xi, \tilde\a}_s - \tilde{h}^{\tilde \a}_s,\q h= b, \si.
\eeaa
Then by  \reff{state1} we have
\beaa
\!\!\!\!\!\!
\dbE\Big[ \big|\cI^{\tilde \a}_{t+\d}( X^\a)\big|^p \Big]
&\!\!\!\!\!\!=&\!\!\!
\dbE\Big[  \big|\cI^{\tilde \a}_t(\xi) + \int_t^{t+\d} b^{\a, \tilde \a}_s  ds + \int_t^{t+\d} \si^{\a, \tilde \a}_s dB_s\big|^p \Big]
\\
&\!\!\!\!\!\!=&\!\!\!\dbE\Big[|\cI^{\tilde \a}_t(\xi) |^p + \int_t^{t+\d} \big[p|\cI^{\tilde \a}_t(\xi) |^{p-2}\cI^{\tilde \a}_t(\xi) \cdot b^{\a, \tilde \a}_s
+{p\over 2}|\cI^{\tilde \a}_t(\xi) |^{p-2}|\si^{\a, \tilde \a}_s|^2  \\
&&\hspace{32mm}+ {p(p-2)\over 2} |\cI^{\tilde \a}_t(\xi) |^{p-4}|\cI^{\tilde \a}_t(\xi) \si^{\a, \tilde \a}_s|^2\big] ds\Big] + o(\d),
\eeaa
where $o(\d)$ is uniform in $\a, \tilde \a$. Similarly,  by \reff{barXtxi} we have
\beaa
\!\!\dbE\Big[ \big|\cI^{\tilde \a}_{t+\d}( \bar X^\a)\big|^p \Big]
&\!\!\!\!\!=& \!\!\!\!\!
\dbE\Big[|\cI^{\tilde \a}_t(\xi) |^p \!+\!\! \int_t^{t+\d}\!\! \big\{p|\cI^{\tilde \a}_t(\xi) |^{p-2}\cI^{\tilde \a}_t(\xi) \!\cdot\! \ol b^{\tilde \a}_s   \\
&&\hspace{15mm}
+{p\over 2}|\cI^{\tilde \a}_t(\xi) |^{p-2}|\ol \si^{\tilde \a}_s|^2 + {p(p\!-\!2)\over 2} |\cI^{\tilde \a}_t(\xi) |^{p-4}|\cI^{\tilde \a}_t(\xi) \ol\si^{\tilde \a}_s|^2\big\} ds\Big]
\!+\! o(\d).
\eeaa
 Then
 \beaa
 \D^{\a,\tilde \a}_\d
 &\le&
 C\dbE\Big[ \int_t^{t+\d} \big[|\cI^{\tilde \a}_t(\xi) |^{p-1}|b^{\a, \tilde \a}_s \!-\! \ol b^{\tilde \a}_s| + |\cI^{\tilde \a}_t(\xi) |^{p-2}|\si^{\a, \tilde \a}_s \!+\! \ol \si^{\tilde \a}_s||\si^{\a, \tilde \a}_s \!-\! \ol \si^{\tilde \a}_s|\big]ds\Big] + o(\d)
 \\
 &\le& C\int_t^{t+\d} \Big[\|\cI^{\tilde \a}_t(\xi) \|_p^{p-1} \|b^{\a, \tilde \a}_s \!-\! \ol b^{\tilde \a}_s\|_p \!+\! \|\cI^{\tilde \a}_t(\xi) \|_p^{p-2}\|\si^{\a, \tilde \a}_s \!+\! \ol \si^{\tilde \a}_s\|_p\|\si^{\a, \tilde \a}_s \!-\! \ol \si^{\tilde \a}_s\|_p\Big]ds + o(\d).
 \eeaa
Recall \reff{cp}. Note that, for $s\in [t, t+\d]$, by Assumption \ref{assum-standing} and standard SDE estimates,
 \beaa
\|\cI^{\tilde \a}_t(\xi) \|_p
\!\!\!&\le& \!\!\!
C\Big(\|\xi_t\|_p + \|\tilde{\xi}_{\tilde t}\|_p + c_p\Big);
\\
\|\si^{\a, \tilde \a}_s + \ol \si^{\tilde \a}_s\|_p
\!\!\!&\le& \!\!\!
C\Big(1+\|X^\a_{\cd\wedge s}\|_p + \|\bar X^\a_{\cd\wedge s}\|_p+c_p\Big)\le C\Big(1+\|\xi_{\cd\wedge t}\|_p +c_p\Big);
\\
\!\!\!\dis  \|b^{\a, \tilde \a}_s \!-\! \ol b^{\tilde \a}_s\|_p^p
&\!\!\!\!\!\!\!=& \!\!\!\!\!\!\!
\dbE\Big[\Big|b_s( X^\a\!, \a_s, \ul X^\a\!, \ul\a_s)
                      \!-\!
                      b_s(X^\a_{\cd\wedge t}, \a_s, \ul X^\a_{\cd\wedge t}, \ul\a_s)\Big|^p\Big]
\le
 C\dbE\Big[|X^\a_{\cd\wedge s}\!-\! X^\a_{\cd\wedge t}|_\infty^p\Big] \!=\! o(1).
\eeaa
Similarly, $\|\si^{\a, \tilde \a}_s - \ol \si^{\tilde \a}_s\|_p \le o(1)$. Then, we derive  \reff{existence-claim} from \reff{phit+d}:
 \beaa
\D^\a_\d \le Ck\Big(1+\|\xi_{\cd\wedge t}\|_p + \|\tilde{\xi}_{\tilde t}\|_p + c_p\Big)^{p-1}o(\d) + o(\d) = o(\d).
\eeaa

 (ii) We next prove the viscosity $p$-supersolution property. Fix  $(t, \xi) \in \ol\cX^p_{[0, T)}$ and $(\varphi, \phi)\in \mathfrak{F}^-_{p}V_t( \ul\xi)$. For any $\d>0$, by DPP \reff{DPP} there exists $\a^\d\in \cA_{[t, T]}$ such that
  \beaa
   0&\ge&  V_{t+\d}( \ul X^{\a^\d}) - V_t(\ul \xi) + \int^{t+\delta}_{t}f _s( \ul X^{\a^\d}, \ul \a^\d_s)ds-\d^2\\
\dis &\ge&  [\f+\phi]_{t+\d}( \ul X^{\a^\d}) - [\f+ \phi]_t( \ul\xi)+ \int^{t+\delta}_{t}f_s( \ul X^{\a^\d},\a^\d_s)ds-\d^2.
 \eeaa
 Now following similar arguments as in (i) we can show
 \beaa
  \partial_t \f_t(\ul\xi)+\limsup_{\d\to 0}{1\over \d}\int_t^{t+\d}\Big[ H_s( \ul\xi_{\cd\wedge t}, \partial_X\f_t(\ul\xi),\partial_{\bx X}\f_t(\ul\xi), \ul \a^\d_s)  + \dot \phi_s( \ul{\bar X}^{\a^\d})\Big]ds \le 0.
  \eeaa
  This  implies \reff{super}, and hence $V$ is a viscosity supersolution.
 \qed

 \section{Connection with HJB equations on the Wasserstein space}\label{sect-MFC0}

\subsection{The mean field control problem with common noise}
\label{sect-MFC}

In this section, we explain that our general setting in this paper covers the mean field control problem with common noise as a special case.
For a generic Polish space $(E,d_E)$, we denote by $\cP_p(E)$ the set of  probability measures $\mu$ on the Borel field $\cB(E)$ with $\int_{E} d_E^p(\bx,\bx_0) \mu(d\bx) <\infty$, for some (and hence for all) $\bx_0\in E$, equipped with the $p$-Wasserstein distance $W_p$. In this section, we decompose $ B^\top= \big((B^1)^\top, (B^0)^\top\big)$, where $B^0$ denotes the common noise, and set $\dbF^0 := \dbF^{B^0}$. Given $\xi\in \cX_2$, let $\dbP_\xi$ denote the law of $\xi$ under $\dbP$, and $\dbP_{\xi|\cF^0_t}$ the conditional law of $\xi$, conditional on $\cF^0_t$. Moreover, for the canonical process ${\bf X}$ on $\dbX$, define $\mu_{\cd\wedge t} := \mu \circ {\bf X}_{\cd\wedge t}^{-1}$, the law of the stopped process.
Let $\cX_p^{\perp}(\dbF^0_{t})$ denote the set of $\xi\in \cX_p$ such that $\xi_{\cd\wedge t}$ is independent of $\cF^0_t$, and for given $(t, \mu)\in [0, T]\times \cP_2(\dbX)$,
\beaa
\cX_p^\perp(\dbF^0; t, \mu) := \Big\{\xi\in \cX_p^\perp(\dbF^0_t): \dbP_{\xi_{\cd\wedge t}} = \mu_{\cd\wedge t}\Big\}.
\eeaa
Our mean field control problem with common noise involves the following data:
\bea
\label{MFCdata}
(\check{b}, \check{\si}, \check f): [0, T]\times \dbX \times  A \times \cP_2(\dbX\times A)
\longrightarrow (\dbR^d, \dbR^{d\times d}, \dbR),\q
\check{g}: \dbX\times \cP_2(\dbX) \longrightarrow \dbR.
\eea
We emphasize that $\check b, \check \si$ are deterministic here.
 Given $(t,\xi)\in \ol\cX^2_{[0, T]}$, define
 \bea
 \label{MFCV}
  \check V_t(\ul \xi) := \inf_{\a\in \cA_{[t, T]}} \check J_t(\ul \xi, \ul \a),
 \q
 \check J_t(\ul \xi, \ul \a)
 :=
 \dbE\Big[\check g\big(X, \dbP_{ X|\cF^0_T}\big)
                +
                \int_t^T \check f_s\big(X,\a_s, \dbP_{ (X, \a_s)|\cF^0_s}\big)ds
        \Big],~~
\eea
where $ X :=  X^{t,\xi,\a}$ is the controlled state process defined by the SDE:
$$
 X_{\cd\wedge t}=\xi_{\cd\wedge t},
 ~\mbox{and} ~
d X_s= \check b_s\big( X,\a_s, \dbP_{ (X, \a_s)|\cF^0_s}\big)ds+ \check \si_s\big( X,\a_s, \dbP_{ (X, \a_s)|\cF^0_s}\big)dB_s,\q s\in [t, T].
$$
We note that here, by writing $\check \si = (\check \si^1, \check \si^0)$, we also allow the control to act on the diffusion coefficient $\check \si^0$ of the common noise.

Similarly to Assumptions \ref{assum-standing} and \ref{assum-regt}, we assume the following.
\begin{assum}
\label{assum-MFC}
{\rm (i)} The data $\check h= \check b, \check \si, \check f$ are  progressively measurable in all variables and adapted in the spirit of \reff{adapted}: $\check h_t(\bx, a, \dbP_{(\xi, \a_*)}) = \check h_t(\bx_{\cd\wedge t}, a, \dbP_{(\xi_{\cd\wedge t}, \a_*)})$; and  $\check h_t({\bf 0}, a, \dbP_{({\bf 0}, \a_*)})$ and $\check g({\bf 0}, \d_{\bf 0})$ are bounded by a constant $C_0$.

\no {\rm (ii)} $\check h= \check b, \check \si$ are uniformly Lipschitz continuous in $(\bx, \mu)\in \dbX\times \cP_2(\dbX)$:
      \beaa
\big|\check h_t(\bx,a, \dbP_{(\xi, \a_*)})-\check h_t( \bx',a,\dbP_{(\xi', \a_*)})\big| \le L\Big[|\bx_{\cd\wedge t} - \bx'_{\cd\wedge t}|_\infty + W_2(\dbP_{\xi_{\cd\wedge t}}, \dbP_{\xi'_{\cd\wedge t}}) \Big].
       \eeaa
\no {\rm (iii)} {$ \check f, \check g$ are uniformly H\"{o}lder-$\b$ continuous in $(\bx, \mu)\in \dbX\times \cP_2(\dbX)$ for some $0<\beta\leq 1$:}
       \beaa
&\dis  {\big|\check f_t(\bx,a, \dbP_{(\xi, \a_*)})-\check f_t( \bx',a,\dbP_{(\xi', \a_*)})\big| \le L\Big[|\bx_{\cd\wedge t} - \bx'_{\cd\wedge t}|^{\beta}_\infty + W^{\beta}_2(\dbP_{\xi_{\cd\wedge t}}, \dbP_{\xi'_{\cd\wedge t}}) \Big];}\\
&\dis  {\big|\check g(\bx,\mu)-\check g( \bx',\mu')\big| \le L\Big[|\bx_{\cd\wedge t} - \bx'_{\cd\wedge t}|^{\beta}_\infty + W^{\beta}_2(\mu, \mu') \Big].}
       \eeaa
       \end{assum}

  \begin{assum}
\label{assum-MFC2}
For any $R>0$, there exists a modulus of continuity function $\rho_R$ such that, for $\check h=  \check b, \check \si, \check f$, for all $t\in [0, T]$, $\bx\in \dbX$,  $a, a'\in A$, $\a_*, \a'_*\in \cA_T$, and all $\|\xi_{\cd\wedge t}\|_2\le R$,
\beaa
\big|\check h_t(\bx, a, \dbP_{(\xi, \a_*)} ) -\check h_t(\bx, a', \dbP_{(\xi, \a'_*)} )\big| \le \rho_R(|a-a'| + \|\a_*-\a'_*\|_2).
\eeaa
\end{assum}

Now given $\check b, \check \si, \check f, \check g$, define
\bea
\label{check-data}
\left.\ba{c}
\dis h_t(\o, \bx, a, \ul \xi, \ul\a_*)   := \check h_t \big(\bx, a, \dbP_{(\xi, \a_*)|\cF^0_t}(\o)\big),\q h = b, \si;\ms\\
\dis  f_t(\ul \xi, \ul \a_*)   := \dbE\big[\check f_t \big(\xi, \a_*, \dbP_{(\xi, \a_*)|\cF^0_t}\big)\big],\q  g(\ul \xi)   := \dbE\big[\check g \big(\xi, \dbP_{\xi|\cF^0_T}\big)\big].
\ea\right.
\eea
We remark that, given $(\xi, \a_*)\in \cX_2\times \cA_T$, we know the joint law $\dbP_{(\xi, \a_*, B^0)}$ and hence the conditional law $\dbP_{(\xi, \a_*)|\cF^0_t} = \dbP_{(\xi, \a_*)|\cF^{B^0}_t}$, so the above functions are well defined.
The following result verifies that  \reff{MFCV} is a special case of the problem \reff{value1}.

\begin{prop}
\label{prop-MFC1}
 Let $\check b, \check \si, \check f, \check g$ satisfy Assumption \ref{assum-MFC}, then  $b, \si, f, g$ defined by \reff{check-data} satisfy Assumption \ref{assum-standing}, and $V(t, \ul\xi) = \check V(t, \ul \xi)$ for all $(t, \xi)\in \ol\cX^2_{[0, T]}$.

 Moreover, if $\check b, \check \si, \check f$ satisfy Assumption \ref{assum-MFC2}, then  $b, \si, f$ satisfy Assumption \ref{assum-regt}.
\end{prop}

\proof We shall only verify the properties of $b$ in Assumption \ref{assum-standing}. Those for $\si$ are similar and those for $f$ and $g$ are slightly easier. In particular, the implication of Assumption \ref{assum-regt} from Assumption \ref{assum-MFC2} is also similar.
First, the properties in Assumption \ref{assum-standing} (i) are obvious. Next, by omitting the dependence on $(\o,a)$, we have
\beaa
&&\dis\!\!\!\!\!\! \big|b_t(\bx,\ul \xi, \ul\a_*)-b_t(\bx',\ul \xi, \ul\a_*)\big|
=
\big|\check b_t \big(\bx, \dbP_{(\xi, \a_*)|\cF^0_t}\big)
        -\check b_t \big(\bx', \dbP_{(\xi,\a_*)|\cF^0_t}\big)\big|
\le L|\bx_{\cd\wedge t}-\bx'_{\cd\wedge t}|_\infty;
\\
&&\dis \!\!\!\!\!\! \big|b_t(\bx,\ul \xi, \ul\a_*)-b_t(\bx,\ul \xi', \ul\a_*)\big|
=
\big|\check b_t \big(\bx, \dbP_{(\xi, \a_*)|\cF^0_t}\big)
        -\check b_t \big(\bx, \dbP_{(\xi',\a_*)|\cF^0_t}\big)\big|\\
&&\q \le
LW_2(  \dbP_{(\xi, \a_*)|\cF^0_t},   \dbP_{(\xi', \a_*)|\cF^0_t})
\le L\Big(\dbE^0_t [|\xi_{\cd\wedge t}-\xi'_{\cd\wedge t}|_\infty^2]\Big)^{1\over 2}.
\eeaa
Moreover, by \reff{check-data} it is clear that $X^{t, \xi, \a} = \check X^{t, \xi, \a}$. Then one can easily get $J_t(\ul\xi, \ul \a_t) = \check J_t(\ul \xi, \ul \a_t)$, and hence $V_t(\xi) = \check V_t(\xi)$.
\qed

\begin{prop}
\label{prop-MFC2}
 Let Assumption \ref{assum-MFC} hold.

\no {\rm (i)}  $\check V$ is conditional law invariant in the sense that $
\check V_t(\ul\xi)=\check V_t(\ul\xi')$ for all $\xi, \xi'\in \ol\cX^2_{[0, T]}$ satisfying $\mathbb{P}_{\xi_{\cd\wedge t}|\cF^0_t}=\mathbb{P}_{\xi'_{\cd\wedge t}|\cF^0_t}$, a.s.
Consequently, by abusing the notation $\check V$, we may define
\bea
\label{checkV}
\check V_t(\mu) := \check V_t(\ul\xi),\q\mbox{for all}~(t, \mu) \in [0, T]\times \cP_2(\dbX),\q \mbox{where}~ \xi\in \cX_2^\perp(\dbF^0; t, \mu).
\eea
\no {\rm (ii)} For any $(t, \xi)\in \ol\cX^2_{[0,T]}$, we have $V_t(\ul\xi) = \dbE\big[\check V_t(\dbP_{\xi|\cF^0_t})\big]$.

 \ms
\no {\rm (iii)} For any $0\le t<t+\d \le T$,  $\xi\in \cX_2^\perp(\dbF^0; t, \mu)$,  we have
\bea
\label{DPP55}
\left.\ba{c}
\dis \check V_t( \mu) =  \inf_{\a\in \cA_{[t, T]}} \dbE\Big[\check V_{t+\d}( \dbP_{X^{t,\xi, \a}|\cF^0_{t+\d}}) + \int^{t+\delta}_{t}\check f _s(X^{t,\xi,\a}, \a_s, \mathbb{P}_{(X^{t,\xi,\a},\a_s)|\cF^0_s})ds\Big].
\ea\right.
\eea
\no {\rm (iv)} For the constant $C$ in Lemma \ref{lem-DPP} {\rm (ii)}, we have
\bea
\label{Vreg55}
\left.\ba{c}
  |\check V_t( \mu)| \le C\Big(1+ W^{\beta}_2(\mu_{\cd\wedge t}, \d_{\bf 0})\Big);\ms\\
\dis |\check V_t(\mu)- \check V_{t'}(\mu')|\leq C\Big(W^{\beta}_2(\mu_{\cd\wedge t}, \mu'_{\cd\wedge t'})  + \big(1+ W^{\beta}_2(\mu_{\cd\wedge t}, \d_{\bf 0})\big) |t-t'|^{\beta\over2}\Big).
\ea\right.
\eea
 \end{prop}
 This result is not surprising, for example, \cite[Theorem 3.6]{cosso3} proved the law invariance property when there is no common noise. However, we believe the precise form under such a generality here is new in the literature.  For completeness we sketch a proof in the Appendix.

To derive the PDE from the DPP \reff{DPP55}, as usual we need the appropriate It\^{o} formula. The following result can be proved by combining \cite[Theorem 4.17]{CD2} and \cite[Theorem 2.7]{WZ}, see also \cite{BLPR, FT24}. The precise meaning of the space $C_b^{1,2,2}([0, T]\times \dbX \times \cP_2(\dbX))$ and that of the derivatives, as well as the proof are again postponed to the Appendix.

\begin{prop}
\label{prop-MFCIto}
 Assume $U\in C_b^{1,2,2}([0, T]\times \dbX \times \cP_2(\dbX))$.  For $i=1,2$,  consider
 $$
 d X^i_t := b^i_t dt + \si^{i,1}_t dB^1_t + \si^{i,0}_t d B^0_t,\q\mbox{and introduce the conditional law}~\mu_{\cd\wedge t}:=\dbP_{X^2_{\cd\wedge t}|\cF^{0}_t},
$$
 where $b^i\in \dbL^2({\dbF}; \dbR^d)$ and $\si^i =(\si^{i,1}, \si^{i,0})\in \dbL^2({\dbF}; \dbR^{d\times d})$. Then
\bea
\label{Ito-common}
&d U_t(X^1, \mu) =  \Big[\pa_t U_t + \pa_\bx U_t \cdot b^1_t + \frac{1}{2} \pa_{\bx\bx} U_t : \si_t^1 (\si_t^1)^\top\Big](X^1, \mu) dt +\pa_\bx U_t(X^1,\mu)\cd\si^1_tdB_t\nonumber\\
& + \dbE_{\cF_t}\big[\pa_\mu U_t(X^1, \mu, \tilde X^2_t) \cd \tilde \si^{2,0}_t \big]dB_t^0
\nonumber\\
&+ \dbE_{\cF_t}\Big[\pa_\mu U_t(\cd,\tilde X^2_t) \cd\tilde b^{2}_t +\frac{1}{2} \pa_{\tilde x}\pa_\mu U_t(\cd, \tilde X^2_t) : \tilde \si_t^2 (\tilde \si_t^2)^\top \\
&\q +\pa_\bx\pa_\mu U(\cd,\tilde X^2_t): \si^{1,0}_t (\tilde \si_t^{2,0})^\top+\frac{1}{2}\pa_{\mu\mu}U_t(\cd,\tilde X^2_t,\bar X^2_t) : \tilde\si_t^{2,0}(\bar \si_t^{2,0})^\top\Big](X^1,\mu) dt,\nonumber
\eea
where $(\tilde X^2, \tilde b^2, \tilde \si^2)$ and $(\bar X^2, \bar b^2, \bar \si^2)$ are conditionally independent copies of $(X^2,  b^2,  \si^2)$, conditional on $\dbF^0$; and by extending the filtered probability space $(\O, \dbF, \dbP)$ in a natural way to include the conditionally independent copies, $\dbE_{\cF_t}$ is the conditional expectation, conditional on $\cF_t = \cF_0 \vee \cF^B_t$.
\end{prop}

Now apply Proposition \ref{prop-MFCIto} on the DPP \reff{DPP55}, we obtain the following second order path dependent HJB equation on the Wasserstein space: for $(t,\mu)\in [0, T]\times \cP_2(\dbX)$,
\bea
\label{HJB-MFC}
 \check \cL \check U_t(\mu)
 :=
 \pa_t \check U_t(\mu)
 + \inf_{\a_t\in \cA_t^\perp(\cF^0_t)}
    \check H_t\big(\mu, \pa_\mu \check U_t(\mu,\cd),
                                    \pa_{\tilde x \mu} \check U_t(\mu, \cd),
                                    \pa_{\mu\mu} \check U_t(\mu, \cd,\cd),
                                    \ul \a_t\big)
=
0,
\eea
where, for $\a_t\in \cA_t^\perp(\cF^0_t)$ independent of $\cF^0_t$,  $\xi \in \cX_2^\perp(\dbF^0; t, \mu)$,  and letting $(\tilde \xi, \tilde \a_t)$, $(\bar\xi, \bar \a_t)$ be independent copies of $(\xi, \a_t)$, the Hamiltonian term $\check H$ is defined as follows:
\bea
\label{MFC-H1}
\check H(\cds)
\!\!\!&:=&\!\!\!
\dbE\Big[\pa_\mu \check U_t(\mu,\tilde \xi_t) \cd\check b_t(\tilde \xi, \tilde \a_t,   \dbP_{(\xi,\a_t)}) +\frac{1}{2} \pa_{\tilde x\mu} \check U_t(\mu, \tilde \xi_t) : \check \si_t \check \si_t^\top(\tilde \xi, \tilde \a_t, \dbP_{(\xi, \a_t)})
\\
\!\!\!&&\!\!\!
+\frac{1}{2}\pa_{\mu\mu}\check U_t(\mu,\tilde \xi_t,\bar \xi_t)
                  \!:\! \check \si_t^{0}(\tilde \xi, \tilde \a_t,  \dbP_{(\xi,\a_t)})
                        (\check \si_t^{0})^\top(\bar \xi, \bar \a_t,  \dbP_{(\xi, \a_t)})
    + \check f_t(\xi, \a_t,  \dbP_{(\xi,\a_t)})\Big].
\nonumber
\eea

\begin{defn}
\label{defn-MFC-viscosity}
We say $\check U\in [0, T]\times \cP_2(\dbX)\to \dbR$ is a viscosity solution (resp. subsolution, supersolution) of the equation \reff{HJB-MFC}-\reff{MFC-H1} if $U_t(\ul\xi) := \dbE[\check U_t(\dbP_{\xi|\cF^0_t})]$ is a viscosity solution (resp. subsolution, supersolution) of the equation \reff{HJB}, where $b, \si, f, g$ are defined by \reff{check-data}.
\end{defn}

The following result is a direct consequence of Theorems \ref{thm-exist} and \ref{thm-comparison}.
\begin{thm}
Let Assumption  \ref{assum-MFC} hold.

\no {\rm (i)} The $\check V$ defined by \reff{checkV} is the unique viscosity solution, in the sense of Definition \ref{defn-MFC-viscosity}, to the equation \reff{HJB-MFC}-\reff{MFC-H1}  with terminal condition $\check g$.

\no {\rm (ii)} Let $\check U_1, \check U_2$ be a viscosity subsolution and viscosity supersolution, respectively,  to the equation \reff{HJB-MFC}-\reff{MFC-H1}, in the sense of Definition \ref{defn-MFC-viscosity}.  Assume, for $\check h  =\check U_1$ and $-\check U_2$, for any $R>0$, there exists a modulus of continuity function $\rho_R$ such that
\beaa
 |\check h_t(\mu)|\le C\big(1+W_2(\mu_{\cd\wedge t}, \d_{\bf 0})\big),\q \check h_t(\mu) - \check h_{t'}(\mu_{\cd\wedge t}) \le \rho_R(t'-t),
 \eeaa
for all $t<t'$ and $W_2(\mu_{\cd\wedge t}, \d_{\bf 0})\le R$. If $\check U_1(T,\mu)\le \int_{\dbX} \check g(\bx, \mu) \mu(d\bx)\le \check U_2(T, \mu)$ for all $\mu\in \cP_2(\dbX)$,  then $\check U_1\le \check U_2$  on $[0, T]\times \cP_2(\dbX)$.
\end{thm}
We remark again that the above notion of viscosity solutions is not equivalent to the definitions in the literature as mentioned in Introduction. Therefore, our results do not imply the uniqueness or comparison principle for the viscosity solutions in those publications.

\begin{rem}
\label{rem-common}
(i) Given $\check U\in  C_b^{1,2}([0, T] \times \cP_2(\dbX))$, similarly to the $\f_1$ in \reff{tildetest}, in general $U_t(\ul\xi) := \dbE[\check U_t(\dbP_{\xi|\cF^0_t})]$ may not be in $C^{1,2}(\cX^2_{[0, T]})$ in the sense of Definition \ref{defn-C12cXp}. So, while we are identifying the viscosity solutions for \reff{HJB-MFC}-\reff{MFC-H1} and  \reff{HJB}, in general  we cannot identify their classical solutions. The situation is  different when there is no common noise, see Proposition \ref{prop-MFCclassical} below.

\no {(ii) It is worth pointing out that, in the common noise case, \reff{Ito-common} and hence the PDE \reff{HJB-MFC}-\reff{MFC-H1} involve the second order derivative $\pa_{\mu\mu} \check U$, however, the corresponding PDE \reff{HJB} in the process space does not involve the second order derivative $\pa_{XX} U$, which is never used in this paper. {It is not clear if this type of first order nature in terms of $\pa_X$ helps in our approach, in particular, our construction of the singular component $\phi$ is mainly to help for the estimates related to the second order derivative $\pa_{\bx\bx}$ (or $\pa_{\bx X}$). We refer to Remark \ref{rem-Ishii2} below for some highly related comments,} and we also refer to Footnotes \ref{1storder} and \ref{2ndorder} for the conventions of orders.  }
\end{rem}

\subsection{Mean field control without common noise}

When there is no common noise, recall Remark \ref{rem-dupire} (ii) and define the classical semisolutions of \reff{HJB-MFC}-\reff{MFC-H1} in the obvious manner. The following result is obvious.

\begin{prop}
\label{prop-MFCclassical}
Let Assumptions  \ref{assum-MFC} and \ref{assum-MFC2} hold. Assume there is no common noise, and $s\in [t, T]\mapsto \check h_s(\bx_{\cd\wedge t}, a, \dbP_{(X_{\cd\wedge t}, \a_*)})$ is uniformly continuous  for $h=b, \si, f$, uniformly in $(\bx, a, X, \a_*)$. Let $\check U\in C^{1,2}([0, T]\times \cP_2(\dbX))$ and denote $U_t(\ul\xi) := \check U(t, \dbP_\xi)$. Then $U\in C^{1,2}(\ol\cX^2_{[0, T]})$, and $U$ is a classical solution (resp. subsolution, supersolution) of \reff{HJB}-\reff{check-data} if and only if $\check U$ is a classical solution  (resp. subsolution, supersolution)  of \reff{HJB-MFC}-\reff{MFC-H1}.
\end{prop}

\begin{rem}
\label{rem-Ishii2}
In the state dependent setting without  common noise, \cite{BIRS, SY1, SY2} proposed different notions of viscosity solutions for HJB equations on the Wasserstein space and established the comparison principle by doubling variable arguments, also without invoking the Crandall-Ishii lemma.  However, their mechanism is completely different from ours.

\ms
\no{\rm (i)} Our HJB equation \reff{HJB} is a second order equation due to the term $\pa_{\bx X} V$, see Footnote \ref{2ndorder}. As in \reff{tildepsi} we consider the following penalization in the doubling variable arguments:
\beaa
\psi(t, \ul\xi, s, \ul\zeta):= U_1(t, \ul\xi) - U_2(s, \ul\eta) - n \dbE[|\xi-\eta|^2] - \cds
\eeaa
Here we consider the state dependent case for simplicity and $\xi, \zeta$ denote random variables instead of processes. Assume the above has optimal arguments $(t_n, \xi_n, s_n, \zeta_n)$. By standard arguments one can easily show $\lim_{n\to \infty}n \dbE[|\xi_n-\eta_n|^2]=0$. Then naturally we would construct test functions (omitting the other terms):
\beaa
\f_1( \ul\xi) := n \dbE[|\xi - \zeta_n|^2],\q \f_2(\ul\zeta):= n \dbE[|\zeta - \xi_n|^2].
\eeaa
Ignoring the possible adaptedness issue, by \reff{Ito} we have
\beaa
\pa_X \f_1(\ul\xi) = 2n (\xi- \zeta_n),~ \pa_{xX} \f_1(\ul\xi) = 2n;\q \pa_X \f_2(\ul\zeta) = 2n (\zeta- \xi_n),~ \pa_{xX} \f_2(\ul\zeta) = 2n.
\eeaa
This implies that,
\bea
\label{Ishii-cancel1}
\pa_X \f_1(\ul\xi_n) + \pa_X \f_2(\ul\zeta_n)=0,\q\mbox{but}\q \pa_{xX} \f_1(\ul\xi_n) =  \pa_{xX} \f_2(\ul\zeta_n).
\eea
Then we can have the desired cancellation for the first order derivatives $\pa_X \f$, for example,
\beaa
\sup_\a \dbE\Big[ b(\xi_n, \a)  \pa_X \f_1(\ul\xi_n) +  b(\zeta_n, \a)  \pa_X \f_2(\ul\zeta_n)\Big]&=&\sup_\a \dbE\Big[ \big(b(\xi_n, \a)   -  b(\zeta_n, \a)\big)  2n (\xi_n- \zeta_n)\Big]
\\
&\le& Cn \dbE[|\xi_n-\eta_n|^2] \longrightarrow 0.
\eeaa
However, there won't be such type of cancellation for the second order derivatives $\pa_{xX}\f$. This is the same issue the second order derivatives $\pa_{xx}\f$ face for standard HJB equations, and in that case the Crandall-Ishii lemma is used exactly to overcome this difficulty. We alternatively introduced the singular component $\phi$, instead of the Crandall-Ishii lemma,  to get around of this difficulty.

\ms
\no {\rm (ii)}  {The HJB equation \reff{HJB-MFC}-\reff{MFC-H1} on the Wasserstein space has a special structure and the corresponding second order term $\pa_{x\mu} \f$ can be magically cancelled,}
as  \cite{BIRS, SY1, SY2}  observed, even though we are talking about the same value function: $V(t, \ul\xi) = \check V(t, \dbP_\xi)$ by Proposition \ref{prop-MFC2} (ii).  To see this, we use the setting in \cite{SY2} and assume state dependence and $d=1$ for simplicity. Introduce a smooth distance function on $\cP_2(\dbR)$ by using the Fourier transform: for some constant $k$ and $i=\sqrt{-1}$:
\beaa
\rho^2(\mu, \nu)  := \int_{\dbR} {|\cF_{\mu-\nu}(z)|^2\over 1+|z|^k} dz,\q\mbox{where}\q \cF_{\mu-\nu}(z) := \int_\dbR e^{-iz\l} (\mu-\nu)(d\l).
\eeaa
Consider the following penalization in the doubling variable arguments:
\beaa
\check\psi(t, \mu, s, \nu):= \check U_1(t, \mu) - \check U_2(s, \nu) - n \rho^2(\mu, \nu) - \cds
\eeaa
and assume it has optimal arguments $(t_n, \mu_n, s_n, \nu_n)$. Then we have  $\dis\lim_{n\to \infty}n  \rho^2(\mu_n, \nu_n)=0$, and we shall similarly construct test functions (again omitting the other terms):
\beaa
\check \f_1(\mu) := n \rho^2(\mu, \nu_n),\q \check\f_2(\nu):= n \rho^2(\mu_n, \nu).
\eeaa
By direct calculation, we have: denoting by $Re$ the real part,
{\small\beaa
\pa_\mu \check \f_1(\mu, x) = 2n \int_{\dbR} {Re\Big(iz e^{izx}\cF_{\mu-\nu_n}(z)\Big)\over 1+|z|^k} dz ,\q \pa_{x\mu} \check \f_1(\mu, x) = -2n \int_{\dbR} {Re\Big(|z|^2 e^{izx}\cF_{\mu-\nu_n}(z)\Big)\over 1+|z|^k} dz;\\
\pa_\mu \check \f_2(\nu, x) = 2n \int_{\dbR} {Re\Big(iz e^{izx}\cF_{\nu-\mu_n}(z)\Big)\over 1+|z|^k} dz ,\q \pa_{x\mu} \check \f_2(\nu, x) = -2n \int_{\dbR} {Re\Big(|z|^2 e^{izx}\cF_{\nu-\mu_n}(z)\Big)\over 1+|z|^k} dz.
\eeaa}
This implies that, unlike \reff{Ishii-cancel1} for the second equality,
\bea
\label{Ishii-cancel2}
\pa_\mu \check \f_1(\mu_n, x) +\pa_\mu \check \f_2(\nu_n, x)=0,\q\mbox{and}\q \pa_{x\mu} \check \f_1(\mu_n, x) + \pa_{x\mu} \check \f_2(\nu_n, x)=0.
\eea
So one may expect the desired cancellation for the terms involving $\pa_{x\mu}\check\f$, and in this sense the HJB equation can be viewed as a first order equation (in a different sense than Footnote \ref{1storder}). However, we should note that the infinite dimensionality, in the sense that $\pa_{x\mu} \check\f(\mu, \cd)$ is a function of $x$, will cause additional difficulty which is not present in the standard finite dimensional HJB equations.

\ms
\no {\rm (iii)} When there is common noise, the HJB equation \reff{HJB-MFC}-\reff{MFC-H1} involves $\pa_{\mu\mu} \check U$. Note that
\beaa
\pa_{\mu\mu} \check \f_1(\mu, x, \tilde x) = 2n \int_{\dbR} {Re\Big(|z|^2 e^{iz(x-\tilde x)}\Big)\over 1+|z|^k} dz =\pa_{\mu\mu} \check \f_2(\nu, x, \tilde x).
\eeaa
 Then we encounter the same difficulty caused by $\pa_{xx}\f$ in the standard case. We remark that \cite{BEZ} introduced a Crandall-Ishii lemma in this setting to overcome this difficulty.
\end{rem}

\section{The comparison result:  Doubling variable in space}
\label{sect-comparison1}

Note that the space $\cX_p$ is not compact. To prove the comparison principle without compactness, we shall use the following variation of the Ekeland-Borwein-Preiss variational principle, see \cite{bor1}.  For this purpose,  we first recall the definition of gauge-type functions. Let $(E, d)$ be a generic complete metric space.
\begin{defn}\label{defn-gaupe}
 A function $\Upsilon \in C^0( E^2; [0,\infty))$  is called a gauge-type function on  $(E, d)$ if
 $$
 \Upsilon(x,x)=0~\mbox{for all}~x\in E,
 ~~\mbox{and}~~
 \lim_{\Upsilon(x_1, x_2)\to 0}d(x_1, x_2)=0.
 $$
\end{defn}

 \begin{lem}\label{lem-variation} 
 Let $\psi: E \rightarrow \dbR$ be upper semicontinuous and  bounded from above, and $\Upsilon$ a gauge-type function on  $(E, d)$. For any $\e>0$ and $x_0 \in E$ satisfying $\psi(x_0) \ge \sup_E\psi-\varepsilon$,  there exist $\hat x\in E$ and  $\{x_i\}_{i\ge 1} \subset E$ such that, denoting $\Psi := \psi - \sum_{i=0}^{\infty} 2^{-i}\Upsilon(\cd,x_i) \le \psi$,
\beaa
\Upsilon(\hat x, x_i)\leq \e\, 2^{-i},
&\mbox{for all}~i\ge 0,~\mbox{and}&
\Psi(\hat x)
=\smax_E\,\Psi
\;\ge\; \psi(x_0)
\;\ge\; \sup_E\psi-\varepsilon.
\eeaa
\end{lem}

\subsection{Smooth gauge-type functions}
Recall Definition \ref{defn-viscosity} (iii), for uniqueness of viscosity solutions it suffices to consider large even integer $p$. For some estimates later, see Footnote \ref{footnote-p} below, from now on we shall always assume $p\ge 6$. Consider the following complete metric space (omitting the dependence of $\L_0, d_0$ on $p$ for notational simplicity):
\bea
\label{Lamda0}
\left.\ba{c}
 \L_0 := \big\{ \theta=(t, \xi_{\cd\wedge t}): (t, \xi) \in \ol\cX^p_{[0, T]}\big\},\q
  d_0\big(\ul\theta, \ul\theta'\big) := |t-t'| +\|\xi_{\cd\wedge t} - \xi'_{\cd\wedge t'}\|_p.
 \ea\right.
 \eea
 The metric $d_0$ itself is of course a gauge-type function. However, for fixed $\theta'$, the mapping $\theta\mapsto d_0\big(\ul\theta,\ul\theta'\big)$ is not smooth in the sense of Definition \ref{defn-C12cXp}.
In order to construct smooth test functions later, we need smooth gauge-type functions. For this purpose, we introduce
\bea
\label{Upsilon0}
\left.\ba{c}
\dis \Upsilon_t( \bx) :=  \frac{(|\bx_{\cdot\wedge t}|_\infty^p-|\bx_t|^{{p}})^3}{|\bx_{\cdot\wedge t}|_\infty^{{{2p}}}}\1_{\{|\bx_{\cd\wedge t}|_\infty\neq 0\}} + 3 |\bx_t|^{{p}},\q (t,\bx) \in [0, T]\times \dbX;\ms\\
\dis \Upsilon_0(\ul\theta,\ul\theta')
:=
\dbE\big[ \Upsilon_{t\vee t'}(\xi_{\cd\wedge t} - \xi'_{\cd\wedge t'})\big],
~\mbox{and}~
\bar \Upsilon_0(\ul\theta,\ul\theta')
:=
\Upsilon_0(\ul\theta,\ul\theta') + |t-t'|^2,
~\theta,\theta'\in \L_0.
\ea\right.
\eea
We emphasize that the map $\Upsilon$ is defined on deterministic paths, rather than on processes, consequently, their path derivatives should be understood in the sense of Dupire \cite{Dupire}. We also recall from \cite{zhou} the useful bounds of $\Upsilon$ and the triangle-like inequality
\bea\label{triangle}
\big|\bx_{\cdot\wedge t}\big|_\infty^{{p}}\leq \Upsilon_t(\bx) \leq 3|\bx_{\cdot\wedge t}|_\infty^{{p}},
\ \Upsilon_t(\bx+\bx')
\le
2^{{{p-1}}}[\Upsilon_t(\bx)+ \Upsilon_t(\bx')],
\ t\in [0, T], \bx,\bx' \in\dbX.~~
\eea
Moreover, since $p\ge 6$, $\Upsilon$ is differentiable in $t$ and twice continuously differentiable in $\bx$ with
\bea\label{lem-Upsilon}
 \pa_t \Upsilon_t( \bx) = 0,
 &|\partial_{\bx}\Upsilon_t(\bx)|\leq {{3p}}|\bx_t|^{{p-1}},&
\mbox{and}\q
|\partial_{\bx\bx}\Upsilon_t(\bx)|\leq {{3p(3p-1)}}|\bx_t|^{{p-2}}.
\eea
Combining this with Example \ref{eg-smooth}, we have the following result.
 \begin{lem}\label{lem-Gauge}
{\rm (i)} $\bar \Upsilon_0$ is a gauge-type function on $(\L_0, d_0)$.

\no {\rm (ii)} For any $\theta'=(t',\xi')\in \ol\cX^{{p}}_{[0, T)}$, the mapping $\theta\in  \ol\cX^{{p}}_{[t', T]} \longmapsto \Upsilon_0(\theta, \theta')$ is in $C^{1,2}(\ol\cX^{{p}}_{[t', T]})$ such that, for any $\hat t\ge t'$, $X\in \cM^{{p}}_{[\hat t, T]}$ as in \reff{cMp}, and $t\in [\hat t, T]$,
we have
\bea\label{statesop3333}
d\Upsilon_0((t, \ul X), \ul{\theta}')
=
\dbE\Big[\partial_{\bx}{\Upsilon}_t(X-\xi'_{\cdot\wedge t'})\!\cd\! dX_t
                +\frac{1}{2}\partial_{\bx\bx}{\Upsilon}_t( X-\xi'_{\cdot\wedge t'})\!:\!d\langle X\rangle_t\Big].
\eea
\end{lem}

\begin{rem}
\label{rem-Gauge}
In the state dependent case $U_t( \ul\xi) = U_t( \ul\xi_t)$, we may consider a much simpler smooth gauge-type function: $\bar\Upsilon_0\big((t, \ul\xi_t), (t', \ul\xi'_{t'})\big) := |t-t'|^2 +\|\xi_t - \xi'_{t'}\|_2^2$. The related calculations will be simplified significantly. However, our approach still requires the singular component $\phi$ in the test functions, so the main arguments for the comparison principle will remain the same. See also Remark \ref{rem-viscosity2} concerning $\phi$ in the case without volatility control.
\end{rem}

We next extend the space $(\L_0, d_0)$ by doubling the spatial variable:
\bea
\label{Lamda1}
\left.\ba{c}
 \dis \L_1(\hat{t}) := \big\{ \l:=(t, \xi_{\cd\wedge t}, \zeta_{\cd\wedge t}): t \in [\hat{t}, T]~\mbox{and}~ \xi, \zeta\in \cX_{{p}}\big\},
 ~\mbox{for all}~
 \hat{t}\in [0,T], \ms\\
 \dis  d_1\big(\ul\l, ~\ul\l'\big) := d_0\big((t, \ul\xi), (t', \ul\xi')\big)+ d_0\big((t, \ul\zeta), (t', \ul\zeta')\big),\ms\\
 \dis \Upsilon_1(\ul\l, \ul\l') := \Upsilon_0\big((t, \ul\xi), (t', \ul\xi')\big) + \Upsilon_0\big((t, \ul\zeta), (t',\ul \zeta')\big),
 \q
 \dis \bar\Upsilon_1(\ul\l, \ul\l') := \Upsilon_1(\ul\l, \ul\l') + |t-t'|^2.
  \ea\right.
 \eea
Similarly, we extend further the space  by doubling the temporal and spatial variables:
\bea
\label{Lamda2}
\left.\ba{c}
\dis \L_2(\hat{t},\hat{s})
:=
\big\{
 \iota:=(\theta,\eta):
\theta=(t,\xi_{\cd\wedge t})\in \ol\cX^{{p}}_{[\hat{t}, T]},
\eta=(s, \zeta_{\cd\wedge s}) \in \ol\cX^{{p}}_{[\hat{s}, T]}
\big\},
~~
 \hat{t},\hat s\in [0,T],
 \ms\\
 \dis  d_2\big(\ul\iota, \ul\iota'\big)
 :=
 d_0(\ul\theta, \ul\theta')
 + d_0(\ul\eta,\ul\eta'),\ms\\
  \dis \Upsilon_2(\ul\iota, \ul\iota') := \Upsilon_0(\ul\theta, \ul\theta') + \Upsilon_0(\ul\eta,\ul\eta'),
 ~\mbox{and}~
 \dis \bar\Upsilon_2(\ul\iota, \ul\iota') := \bar\Upsilon_0(\ul\theta, \ul\theta') + \bar\Upsilon_0(\ul\eta,\ul\eta').
  \ea\right.
 \eea
Clearly $\bar \Upsilon_1$ and $\bar \Upsilon_2$ are gauge-type functions on $(\L_1(\hat{t}), d_1)$ and $(\L_2(\hat{t},\hat{s}), d_2)$, respectively.

 In this  and the next sections,  we prove the comparison principle by using the doubling variable arguments. In $(\L_1(\hat{t}), d_1)$ we double the spatial variable only, while in $(\L_2(\hat{t},\hat{s}), d_2)$ we double both the spatial and temporal variables. From now on we consider the setting in Theorem \ref{thm-comparison} and assume all the conditions there hold true. {Recall Definition \ref{defn-viscosity} (iii), we may assume $U_1 (\mbox{resp.}, U_2)$ is a viscosity $p$-subsolution (resp., $p$-supersolution) of HJB equation (\ref{HJB}) for some even integer $p$ satisfying:\footnote {\label{footnote-p}Since $p$ is even and $\b\le 1$, this implies $p\ge 6$.}
 \bea
 \label{pbeta}
 p > {4\over \b}.
 \eea
 }
 Moreover, we assume without loss of generality that $U_1$ satisfies \reff{Vreg1} and denote:
 \beaa
  m:= \sup_{\th\in {\ol\cX^{p}_{[0, T]}}} [U_1-U_2](\th).
 \eeaa
We shall assume to contrary that $m>0$, and work toward a contradiction.

\subsection{Doubling variable in space}
 For any $n\ge 1$ and $\e>0$, introduce the function $\psi^{n,\e}$ defined for $\lambda=(t,\xi,\zeta)\in\L_1(0)$ by:
\begin{eqnarray}
\label{psine}
\psi^{n,\e}(\ul\l)
= \psi^{n,\e}_t(\ul\xi,\ul\zeta)
\!\!&:=&\!\!
U_1(t,\ul\xi)-U_2(t,\ul\zeta)- n\dbE\big[\Upsilon_t(\xi-\zeta)\big]- {n^{q}}\,\mathbb{E}\big[|\xi_{t}-\zeta_{t}|^{{p}}\big]\ms\ss
\\
&&\!\!\!\!\!\!\! -\e\Big(1- {t\over 2T}\Big)\mathbb{E}\big[ \Upsilon_t(\zeta)\big],~  {\mbox{with}~ q:=1+\frac{1}{2}\left[\frac{2}{p-2}+\frac{\beta(p-2)}{2p-2\beta}\right].}
\notag
 \end{eqnarray}
 { By \reff{pbeta}, one can verify straightforwardly that
\bea
\label{pbeta2}
{p\over {p-2}}-{\beta\over {p-\beta}}< q<1+{\beta p\over {2p-2\beta}}.
\eea
These will be crucial for the estimate \reff{pbeta3} below. }

 \begin{prop}
 \label{prop-psine}
Let all the conditions in Theorem \ref{thm-comparison} hold true with the above convention on $p$, including \reff{pbeta}, and $U_1$ satisfy \reff{Vreg1}. Assume to the contrary that $m>0$. Then there exists $0<\e_{m}\le 1$ such that, for all $0<\e \leq \e_{m}$, the following hold.

\ms
\no{\rm (i)} There exist  $\hat\l^{n,\e}=(\hat t^{n, \e}, \hat \xi^{n,\e}, \hat \zeta^{n,\e}) \in \L_1(0)$ and $\Psi^{n,\e} \in C^0(\L_1(\hat t^{n,\e}))$ s.t. $\Psi^{n,\e}\le \psi^{n,\e}$ and
\bea
\label{Psinemax}
 \Psi^{n,\varepsilon}(\ul{\hat \l}^{n, \e})
=\max_{\l\in \L_1(\hat t^{n,\e})} \,\Psi^{n,\varepsilon}(\ul\l)
\;\ge\; \sup_{\l\in \L_1(0)}\psi^{n, \e}(\ul \l)-\frac1n
\;\ge\;
\frac{m}{2}-\frac1n.
\eea

\no{\rm (ii)} Moreover, $[\psi^{n,\e} - \Psi^{n,\e}](\ul\l)=\pi^{n,\e}_0(t)+ \pi^{n,\e}_1(t, \ul{\xi}) + \pi^{n,\e}_2(t, \ul{\zeta})$, for some nonnegative functions $\pi^{n,\e}_0\in C^1([\hat t^{n,\e}, T])$,  $\pi^{n,\e}_1, \pi^{n,\e}_2\in C^{1,2}(\ol\cX^{{p}}_{[\hat t^{n,\e}, T]})$, such that: at $\hat\l^{n,\e}$,
\bea
\label{pineest}
\left.\ba{c}
\dis 0\le \pi^{n,\e}_j\le Cn^{-1}, ~j=0,1,2;\q |\pa_t \pi^{n,\e}_0| \le Cn^{-\frac12}; \ms\\
\dis \pa_t \pi^{n,\e}_j =0,\q \|\pa_X \pi^{n,\e}_j\|_{{{p}\over {p-1}}} \le C\textcolor{black}{n^{-\frac{p-1}{p}}},\q \|\pa_{\bx X}  \pi^{n,\e}_j\|_{\textcolor{black}{p\over p-2}} \le C \textcolor{black}{n^{-\frac{p-2}{p}}}, \q j=1,2.
\ea\right.
\eea

\no{\rm (iii)} There exists $C_\e>0$, which may depend on $\e$ but not on $n$, such that
\begin{eqnarray}
\label{CTe}
\mathbb{E}\Big[|\hat \xi^{n,\e}_{\cdot\wedge\hat t^{n,\e}}|_\infty^{\textcolor{black}{p}}
   +|\hat \zeta^{n,\e}_{\cdot\wedge\hat t^{n,\e}}|_\infty^{\textcolor{black}{p}}\Big]
&\le&
   C_\e,~\mbox{and}\ms
  \\
 n\mathbb{E}\Big[|\hat\xi^{n,\e}_{\cdot\wedge\hat t^{n,\e}}
                               -\hat \zeta^{n,\e}_{\cdot\wedge\hat t^{n,\e}}|_\infty^{\textcolor{black}{p}}\Big]
   +\textcolor{black}{n^{q}}\mathbb{E}\Big[|\hat\xi^{n,\e}_{\hat t^{n,\e}}
                                                     -\hat\zeta^{n,\e}_{\hat t^{n,\e}}|^{\textcolor{black}{p}}\Big]
   &\le&
    C \textcolor{black}{n^{-{\beta\over {p-\beta}}}},~\mbox{for all}~ n\ge 1.\q
\label{xietane}
\end{eqnarray}
{\rm (iv)} $\overline{T}_{\e, m}:=\limsup_{n\to\infty} \hat t^{n,\e} < T$. Consequently, there exists $n_{\e, m}$, which may depend on $\e$ but not on $n$, such that $\hat t^{n,\e}\le T_{\e, m} := {1\over 2}\big(\overline{T}_{\e, m} + T\big)<T$, for all $n \ge n_{\e, m}$.
\end{prop}

 \proof (i) As $m>0$, there exists $\th_0=(t_0, \xi^0)\in \L_0$ s.t. $(U_1-U_2)(\ul\th_0) \ge {2m\over 3}$. Then
 $$
 \sup_{\L_1(0)}\psi^{n,\e}
 \ge
 \psi^{n,\e}(t_0, \ul\xi^0, \ul\xi^0)
 =
 (U_1-U_2)(\ul\th_0) - \e\Big(1-{ t_0\over 2T}\Big)\dbE[\Upsilon(\th_0)]
 \ge
 {2m\over 3} - 3\e \dbE[|\xi^0_{\cd\wedge t_0}|_\infty^{\textcolor{black}{p}}],
 $$
 where the last inequality is due to \reff{triangle}.  By setting $\e_{m} := {m\over 1+18 \dbE[|\xi^0_{\cd\wedge t_0}|_\infty^{\textcolor{black}{p}}]}$, we obtain the last inequality in \reff{Psinemax}. Next, by \reff{triangle} we have
 \beaa
\mathbb{E}\big[\Upsilon_t(\zeta)+\Upsilon_t(\xi-\zeta)\big] \ge {1\over {\textcolor{black}{2^{p-1}+1}}}\mathbb{E}\big[\Upsilon_t(\xi)+\Upsilon_t(\zeta)\big] \ge {1\over {\textcolor{black}{2^{p-1}+1}}}[\|\xi_{\cd\wedge t}\|_{\textcolor{black}{p}}^{\textcolor{black}{p}} + \|\zeta_{\cd\wedge t}\|_{\textcolor{black}{p}}^{\textcolor{black}{p}}].
 \eeaa
Note that $n\ge 1 \ge \e$ and $1-{t\over 2T} \ge {1\over 2}$, then
\beaa
n\dbE\big[\Upsilon_t(\xi-\zeta)\big]+\e\Big(1- {t\over 2T}\Big)\mathbb{E}\big[ \Upsilon_t(\zeta)\big] &\ge& {\e\over 2}\mathbb{E}\big[\Upsilon_t(\zeta)+\Upsilon_t(\xi-\zeta)\big] \\
&\ge& {\e\over {\textcolor{black}{2^{p+1}}}}\big[\|\xi_{\cd\wedge t}\|_{\textcolor{black}{p}}^{\textcolor{black}{p}} + \|\zeta_{\cd\wedge t}\|_{\textcolor{black}{p}}^{\textcolor{black}{p}}\big].
\eeaa
Thus, since $U_1, U_2$ satisfy the estimates  \reff{Ureg},
 \bea
\label{psinelower1}
  \psi^{n,\e}(\ul\l) \le C\big[1+\|\xi_{\cd\wedge t}\|_2 + \|\zeta_{\cd\wedge t}\|_2\big] - {\e\over {\textcolor{black}{2^{p+1}}}} \big[\|\xi_{\cd\wedge t}\|_{\textcolor{black}{p}}^{\textcolor{black}{p}} + \|\zeta_{\cd\wedge t}\|_{\textcolor{black}{p}}^{\textcolor{black}{p}}\big].
  \eea
  This clearly implies $\sup_{\L_1(0)} \psi^{n,\e} <\infty$, and thus there exists  $\l^{n,\e}_0 = (t^{n,\e}_0, \xi^{n,\e}_0, \zeta^{n,\e}_0)\in \L_1(0)$ such that $\psi^{n,\e}(\ul\l^{n,\e}_0) \ge \sup_{\L_1(0)} \psi^{n,\e}-{1\over n}$. Now apply Lemma \ref{lem-variation} on $(\L_1(t^{n, \e}_0), d_1)$, there exist $\hat\l^{n,\e} = (\hat t^{n,\e}, \hat \xi^{n,\e}, \zeta^{n,\e})\in \L_1(t^{n, \e}_0)$ and $\l^{n,\e}_i = (t^{n,\e}_i, \xi^{n,\e}_i, \zeta^{n,\e}_i)\in \L_1(t^{n, \e}_0)$, $i\ge 1$,  such that

\ms
$\bullet$ $t^{n,\e}_i \uparrow \hat t^{n,\e}$ as $i\to \infty$, where the monotonicity of $(t^{n,\e}_i)_i$ is due to  \cite[Lemma 2.14]{zhou};

\ms
$\bullet$ $\ol \Upsilon_1(\hat \l^{n,\e}, \l^{n,\e}_i) \le {1\over n 2^i}$, $i\ge 0$, and $\dis \Psi^{n,\e} := \psi^{n,\e} -  \sum_{i=0}^\infty 2^{-i}\ol\Upsilon_1(\ul\l,  \l^{n,\e}_i) \le \psi^{n,\e}$ satisfies \reff{Psinemax}.

 (ii)  For notational convenience, in the rest of this proof we omit the superscript $^{(n,\e)}$, e.g. $\hat \l = \hat\l^{n,\e}$, $\l_i = \l^{n,\e}_i$, $\psi = \psi^{n,\e}$, $\Psi = \Psi^{n,\e}$. By the definition of $\Psi$ and \reff{Lamda1}, we have $[\psi - \Psi](\ul\l)=\pi_0(t)+ \pi_1(t, \ul{\xi}) + \pi_2(t, \ul{\zeta})$ with nonnegative maps
 $\pi_0(t) := \sum_{i=0}^\infty 2^{-i} (t-t_i)^2$,
 \bea
 \label{pi12}
 \pi_1(t, \ul\xi) := \sum_{i=0}^\infty 2^{-i} \Upsilon_0((t,\ul\xi), (t_i, \ul\xi_i)),
 ~\mbox{and}~~
 \pi_2(t, \ul\zeta) := \sum_{i=0}^\infty 2^{-i}\Upsilon_0((t,\ul\zeta), (t_i, \ul\zeta_i)).
 \eea
The regularity of $\pi^{n,\e}_j$ follow from Lemma \ref{lem-Gauge}, and it only remains to verify \reff{pineest}. First,
\beaa
 &&\dis \pi_0(t)+ \pi_1(t, \ul{\xi}) + \pi_2(t, \ul{\zeta}) = \sum_{i=0}^\infty 2^{-i} \ol \Upsilon_1(\ul{\hat \l}, \ul\l_i) \le \sum_{i=0}^\infty {1\over n 2^{2i}} \le {2\over n},
 \\
&&\mbox{and}~~\dis \pa_t \pi_0(\hat t) = \sum_{i=0}^\infty 2^{-i} (\hat t-t_i) \le  \sum_{i=0}^\infty {1\over 2^{i}} \Big({1\over n2^i}\Big)^{1\over 2} \le {C\over \sqrt{n}}.
\eeaa
Next, by  Lemmas \ref{lem-Gauge} and the estimates \eqref{lem-Upsilon} we have $\pa_t \pi_1 =0$, and
$$
\dis \|\pa_X \pi_1(\hat t,  \ul {\hat \xi}) \|_{\textcolor{black}{p\over p-1}}
\le
\sum_{i=0}^\infty 2^{-i}\|\pa_\bx \Upsilon_{\hat t}(\hat \xi_{\cd\wedge \hat t}- (\xi_i)_{\cd\wedge t_i})\|_{\textcolor{black}{p\over p-1}}
\le
C\sum_{i=0}^\infty 2^{-i}\||\hat \xi_{\hat t}- (\xi_i)_{t_i}|^{\textcolor{black}{p-1}} \|_{\textcolor{black}{p\over p-1}}
$$
Then, it follows from the lower bound of $\Upsilon$ in \eqref{triangle} that
\beaa
\|\pa_X \pi_1(\hat t,  \ul {\hat \xi}) \|_{\textcolor{black}{p\over p-1}}
&\le&
C\sum_{i=0}^\infty 2^{-i}\|\hat \xi_{\hat t}- (\xi_i)_{t_i}\|_{\textcolor{black}{p}}^{\textcolor{black}{ p-1}}
\le
C\sum_{i=0}^\infty 2^{-i} \|\Upsilon_{\hat t}(\hat \xi_{\cd\wedge \hat t}- (\xi_i)_{\cd\wedge t_i})\|_1^{\textcolor{black}{p-1\over p}}
\\
&=&
C\sum_{i=0}^\infty 2^{-i} \Big(\Upsilon_0((\hat t, \hat \xi), (t_i, \xi_i))\Big)^{\textcolor{black}{p-1\over p}}
\le
C\sum_{i=0}^\infty 2^{-i}\Big({1\over n 2^i}\Big)^{\textcolor{black}{p-1\over p}} \le C\, n^{-{\textcolor{black}{p-1\over p}}}.
\eeaa
Similarly, as $\dis|\pa_{\bx X} \pi_1(\hat t,  \ul {\hat \xi})| \le  \sum_{i=0}^\infty 2^{-i}|\pa_{\bx\bx} \Upsilon_{\hat t}(\hat \xi_{\cd\wedge \hat t}- (\xi_i)_{\cd\wedge t_i})|$,
we may get $\|\pa_{\bx X} \pi_1(\hat t,  \ul {\hat \xi}) \|_{\textcolor{black}{p\over p-2}} \le C\, n^{-{\textcolor{black}{p-2\over p}}}$. The estimates for $\pi_2$ can be proved similarly.

 \ms

 (iii) First, by \reff{Psinemax} we have
\bea
\label{psinelower2}
\psi(\ul{\hat\l}) \ge \Psi(\ul{\hat\l}) \ge {m\over 2} - {1\over n} \ge -1.
\eea
Then, noting that $\|\cd\|_2 \le \|\cd\|_{\textcolor{black}{p}}$, it follows from \reff{psinelower1} that
$$
-1 \le C\big[1+\|\hat\xi_{\cd\wedge \hat t}\|_{\textcolor{black}{p}} + \|\hat\zeta_{\cd\wedge \hat t}\|_{\textcolor{black}{p}} \big] - {\e\over {\textcolor{black}{2^{p+1}}} } \big[\|\hat\xi_{\cd\wedge \hat t}\|_{\textcolor{black}{p}}^{\textcolor{black}{p}}  + \|\hat \zeta_{\cd\wedge \hat t}\|_{\textcolor{black}{p}}^{\textcolor{black}{p}}\big].
$$
This implies  \reff{CTe}  with $C_\e = C\e^{-{\textcolor{black}{p\over p-1}}}$.

Next, as $\Psi(\ul{\hat\l}) \ge \Psi(\hat t, \ul{\hat\zeta}, \ul{\hat\zeta}) $ by the optimality of $\hat\l$ in \reff{Psinemax},  by  \reff{psine} and (ii) we have
 \bea
 \label{U1difference}
 \left.\ba{c}
\dis 0\le
 U_1(\hat t, \ul{\hat\xi})- U_1(\hat t, \ul{\hat\zeta})
 - n\dbE\big[\Upsilon_{\hat t}(\ul{\hat\xi}- \ul{\hat\zeta})\big]- \textcolor{black}{n^{q}}\mathbb{E}\big[|\hat\xi_{\hat t}-\hat\zeta_{\hat t}|^{\textcolor{black}{p}}\big]
 \ms\\
 \dis  + \sum_{i\ge 0}2^{-i}
 \dbE\Big[\Upsilon_{\hat t}(\hat \zeta - (\xi_i)_{\cdot\wedge t_i})
               -\Upsilon_{\hat t}(\hat \xi- (\xi_i)_{\cdot\wedge t_i})
        \Big].
        \ea\right.
 \eea
Since $\Upsilon \ge 0$, we have
  \beaa
&\dis n\dbE\big[\Upsilon_{\hat t}(\ul{\hat\xi}- \ul{\hat\zeta})\big]+\textcolor{black}{n^{q}}\mathbb{E}\big[|\hat\xi_{\hat t}-\hat\zeta_{\hat t}|^{\textcolor{black}{p}}\big] \le K_1 + K_2,\q \mbox{where}\ms\\
&\dis
K_1:= \big|U_1(\hat t, \ul{\hat\xi})- U_1(\hat t, \ul{\hat\zeta})\big|,\q
K_2:= \sum_{i\ge 0}2^{-i} \dbE\big[\Upsilon_{\hat t}(\hat \zeta_{\cdot\wedge\hat t} - (\xi_i)_{\cdot\wedge t_i})\big].
 \eeaa
 Since $U_1$  satisfies  \reff{Vreg1}, by  \reff{CTe} we have
 \beaa
 K_1\le
 C\|\hat \xi_{\cdot\wedge\hat t}-\hat \zeta_{\cdot\wedge\hat t}\|_{\textcolor{black}{p}}^{\textcolor{black}{\beta}}.
\eeaa
Moreover, recall from (i) that $\ol \Upsilon_1(\hat \l, \l_i) \le {1\over n 2^i}$. Then, by  \eqref{triangle}  we have
\beaa
 K_2
\le
 C\sum_{i=0}^\infty2^{-i}\dbE\Big[\Upsilon_{\hat t}(\hat \zeta-\hat \xi)+ \Upsilon_{\hat t}( \hat \xi - (\xi_i)_{\cdot\wedge t_i})  \Big] \le C\|\hat \xi_{\cdot\wedge\hat t}-\hat \zeta_{\cdot\wedge\hat t}\|_{\textcolor{black}{p}}^{\textcolor{black}{p}} + {C\over n}.
 \eeaa
Recall from \eqref{triangle} that $\Upsilon_t(\bx) \ge |\bx_{\cd\wedge t}|_\infty^{\textcolor{black}{p}}$. Then
$$
n \|\hat \xi_{\cdot\wedge\hat t}-\hat \zeta_{\cdot\wedge\hat t}\|_{\textcolor{black}{p}}^{\textcolor{black}{p}}
+ \textcolor{black}{n^{q}}\|\hat\xi_{\hat t}-\hat\zeta_{\hat t}\|_{\textcolor{black}{p}}^{\textcolor{black}{p}}
\le
 C \|\hat \xi_{\cdot\wedge\hat t}-\hat \zeta_{\cdot\wedge\hat t}\|_{\textcolor{black}{p}}^{\textcolor{black}{\beta}}+C  \|\hat \xi_{\cdot\wedge\hat t}-\hat \zeta_{\cdot\wedge\hat t}\|_{\textcolor{black}{p}}^{\textcolor{black}{p}}  +{C\over n}.
$$
Note that \textcolor{black}{$ C x^{\beta} \le {1\over 2} n x^{p} +  C n^{-{\beta\over {{p}-\beta}}}$}. Then, for $n$ large we have
\beaa
n \|\hat \xi_{\cdot\wedge\hat t}-\hat \zeta_{\cdot\wedge\hat t}\|_{\textcolor{black}{p}}^{\textcolor{black}{p}}  + \textcolor{black}{n^{q}}\|\hat\xi_{\hat t}-\hat\zeta_{\hat t}\|_{\textcolor{black}{p}}^{\textcolor{black}{p}}  \le  C n^{-{\textcolor{black}{\beta\over p-\beta}}}.
\eeaa
 This implies \reff{xietane} immediately.

\ms

(iv) First, for $n \ge {6\over m}$, by \reff{psinelower2} we have $\psi(\ul{\hat\l})\ge {m\over 2}-{1\over n}\ge {m\over 3}$.
 Then, by the second inequality in \reff{Ureg} and since $U_1(T, \cd) \le U_2(T, \cd)$, we derive from \reff{psine}, \reff{CTe},  \reff{xietane}    that
  \beaa
 &\dis {m\over 3}
 \le
 \psi(\ul{\hat\l})
 \le
 U_1(\hat t, \ul{\hat \xi})-U_2(\hat t, \ul{\hat \zeta})
 \le
 \big[U_1(\hat t,\ul{\hat \xi}_{\cdot\wedge\hat t})
        -U_2(\hat t, \ul{\hat \zeta}_{\cdot\wedge\hat t})\big]
\!-\!  \big[U_1(T, \ul{\hat \zeta}_{\cdot\wedge\hat t})-U_2(T, \ul{\hat \zeta}_{\cdot\wedge\hat t})\big]
\\
&\dis \le
 2\rho_{C_\e}(T-\hat t)+C\|\hat \xi_{\cdot\wedge\hat t}-\hat \zeta_{\cdot\wedge\hat t}\|^{\textcolor{black}{\beta}}_2\le 2\rho_{C_\e}(T-\hat t)+ C_\e \textcolor{black}{n^{-\frac{\beta}{p-\beta}}}.
\eeaa
The required result follows from taking $\liminf_{n\to\infty}$ in the last inequality.
 \qed

\section{The comparison result: a crucial estimate}
\label{sect-comparison2}

To prove the comparison principle, we need another proposition. Recall the notations $\e_{m}$, $C_\eps$ and $n_{\e, m}$ defined in Proposition \ref{prop-psine}.

\begin{prop}
\label{prop-comparison1}
Let all the conditions in Proposition \ref{prop-psine}  hold true. Assume further that
 \bea
\label{U2c}
\mbox{$U_2$ is a viscosity supersolution of  $\cL U_2(t, \ul\xi) \le -c_1$ for some $c_1 >0$.}
\eea
Then there {exist constants $C>0$ and $\tilde{q}>0$}, independent of $n, \e, m, c_1$, such that
 \bea
 \label{comparison1}
c_1
\le
\e (C -{1\over 2T})
    \big\|\hat\zeta^{n,\e}_{\cdot\wedge\hat t^{n,\e}}\big\|_{\textcolor{black}{p}}^{\textcolor{black}{p}}
     + C\e
+C_\e {n^{-\tilde{q}}},
~\mbox{for all}~
\e\le \e_{m},~n\ge n_{\e, m}.
\eea
 \end{prop}

We note that $\psi^{n,\e}$ in \eqref{psine} doubles the spatial variable only.
The rather lengthy proof of this proposition requires additional doubling variables in time and space.

\subsection{Doubling variable in time}

Let $\e_m,n_{\e,m}$, $\hat \l = \hat\l^{n,\e}=(\hat t^{n,\e}, \hat\xi^{n,\e}, \hat \zeta^{n,\e})$  be as  in Proposition \ref{prop-psine}, and set:
\begin{align*}
 \hat \arm=\hat\arm^{n,\e} := {1\over 2}[\hat\xi^{n,\e} + \hat \zeta^{n,\e}],
\q \hat\th=\hat\th^{n,\e}:=(\hat t^{n,\e}, \hat\xi^{n,\e}),
   \q\hat\eta=\hat\eta^{n,\e}:=(\hat t^{n,\e}, \hat\zeta^{n,\e}),
\end{align*}
for $0<\e\leq\e_m$ and $n\ge n_{\e,m}$, where we omit the superscripts $^{n,\e}$ when the contexts are clear.
For each $N\ge 1$ and $\iota =(\th,\eta)\in \L_2(\hat{t})$, with $\th=(t,\xi)$ and $\eta=(s,\zeta)$, define
 \bea
\psi^{N}_{n,\e}(\ul\iota)
&:=&
U_1({\ul\th})-U_2({\ul\eta})
- \big[\f^{\hat\th}_1(\ul\th,\ul{\hat\arm})
         + \f^{\hat\eta}_2(\ul\eta,\ul{\hat\arm}) \big] - \pi^{n,\e}_0(t) - N|s-t|^2
\nonumber
\\&& \hspace{24mm}
-\sup_{\a\in \cA_{[\hat t, T]}}
                 \cJ_-^{\hat \th,\hat\arm}(\theta,\alpha)
-\inf_{\a\in \cA_{[\hat t, T]}} \cJ_+^{\hat \th,\hat\arm}(\eta,\alpha),
 \label{psineN}
\eea
where, recalling the notations in \reff{cIaxi} and \reff{barXtxi}, the penalties $\f^{\hat\th}_j$ and $\cJ_\pm$ are given by
\bea
\label{fne}
\left.\ba{c}
\dis  \f^{\hat\th}_1(\ul\th,\ul{\hat\arm})
:=
\bar\Upsilon_0\big(\ul\th,  \ul{\hat\th}\big)
+\textcolor{black}{2^{p-1}} n\Upsilon_0(\ul\th,(\hat t, \ul{\hat \arm}))
+\pi^{n,\e}_1(\ul\th);
\ms\\
\dis \f^{\hat\eta}_2(\ul\eta,\ul{\hat\arm})
:=
\bar\Upsilon_0\big(\ul\eta,  \ul{\hat\eta}\big)
+\textcolor{black}{2^{p-1}} n\Upsilon_0(\ul\eta,(\hat t, \ul{\hat \arm}))
+\e\Big(1- {s\over 2T}\Big)\dbE[\Upsilon_s(\zeta)]
+\pi^{n,\e}_2(\ul\eta);\ms\\
\dis \cJ_{\pm}^{\hat\th,\hat\arm}(\th,\alpha)
:=
\textcolor{black}{2^{p-1}n^{q}}\mathbb{E}\big[|\hat \cI^{\a}_t( \xi)|^{\textcolor{black}{p}}  \big]
                         \pm \hat F^{\a}_t,\q  \hat \cI^\a:= \cI^{b^{\hat\th},\si^{\hat\th},\hat t,\hat \arm,\a}, ~~\hat F^\a:= F^{f^{\hat\th}, \hat t, \a}.
\ea\right.
\eea
We emphasize that the term $\cJ_+^{\hat \th,\hat\arm}(\eta,\a)$ in \eqref{psineN} also relies on $\hat\th$, and not on $\hat\eta$.

 \begin{prop}
 \label{prop-psineN1}
Let all the conditions in  Proposition \ref{prop-psine}  hold true. Then for any $0<\e\le \e_m$, $n\ge n_{\e,m}$, the following hold.

\ms
\no{\rm (i)} There exist   $\check\iota^{N}=(\check \th^N, \check \eta^N)=\big((\check t^{N}, \check \xi^{N}), (\check s^{N}, \check  \zeta^{N})\big) \in \L_2(\hat{t},\hat{t})$ (omitting $^{n,\e}$) and $\Psi^N_{n,\e} \in C^0\big(\L_2(\hat{t},\hat{t})\big)$ such that
$\Psi^N_{n,\e} \le \psi^N_{n,\e}$ and
\bea
\label{PsineNmax}
 \Psi^N_{n,\varepsilon}(\ul{\check \iota}^N)
=\max_{\iota\in \L_2(\hat t, \hat t)} \,\Psi^N_{n,\varepsilon}(\ul\iota)
\;\ge\; \Big(\sup_{\iota\in \L_2(\hat t, \hat t)} \,\psi^N_{n,\varepsilon}(\ul\iota)-\frac1N\Big) \vee \Psi^{n,\e}(\hat \l).
\eea

\no{\rm (ii)} Moreover, $[\psi^N_{n,\e} - \Psi^N_{n,\e}](\ul\iota)= \pi^{n,\e,N}_1(\ul{\th}) + \pi^{n,\e,N}_2(\ul{\eta})$, for some nonnegative functions $\pi^{n,\e,N}_1\in C^{1,2}(\ol\cX^{\textcolor{black}{p}}_{[\check t^N, T]})$, $\pi^{n,\e,N}_2\in C^{1,2}(\ol\cX^{\textcolor{black}{p}}_{[\check s^N, T]})$, such that: for $j=1,2$ and at $\check\iota^N$,
\bea
\label{pineNest}
\dis 0\le \pi^{n,\e,N}_j\le {C\over N}; |\pa_t \pi^{n,\e,N}_j| \le {C\over \sqrt{N}};  \|\pa_X \pi^{n,\e,N}_j\|_{\textcolor{black}{p\over p-1}} \le {C\over N^{\textcolor{black}{p-1\over p}}}, \|\pa_{\bx X}  \pi^{n,\e,N}_j\|_{\textcolor{black}{p\over p-2}} \le {C\over N^{\textcolor{black}{p-2\over p}}}.
\eea

\no {\rm (iii)} $\sup_{n,N}\mathbb{E}\Big[\|\check \xi^{N}_{\cdot\wedge\check t^{N}}\|^{\textcolor{black}{p}}_{\textcolor{black}{p}}+\|\check \zeta^{N}_{\cdot\wedge\check s^{N}}\|^{\textcolor{black}{p}}_{\textcolor{black}{p}}\Big] \le C_\e,$ for some $C_\e$, which may depend on $\e$.
\end{prop}
\proof Recall the connection between $\Psi^{n,\e}$ and $\pi^{n,\e}_j$ in Proposition \ref{prop-psine} (ii) and that $\hat \xi -\hat \arm = \hat\arm - \hat \zeta = {1\over 2}[\hat \xi - \hat \zeta]$. Then, by \reff{psineN}, \reff{psine}, \reff{cIaxi}, and \reff{Upsilon0} we have
\bea
\label{psiNPsi}
&&\psi^N_{n,\e}\big(\ul{\hat \th}, \ul{\hat \eta}\big)- \Psi^{n,\e}(\ul{\hat \l})
=
\psi^N_{n,\e}\big(\ul{\hat \th}, \ul{\hat \eta}\big)- \psi^{n,\e}(\ul{\hat \l})
+ \pi^{n,\e}_0(\hat t) + \pi^{n,\e}_1(\ul{\hat \th}) + \pi^{n,\e}_2(\ul{\hat \eta})\nonumber
\\
&=&
n \Upsilon_0\big(\ul{\hat \th}, \ul{\hat \eta}\big) + \textcolor{black}{n^{q}}\dbE\Big[|\hat\xi_{\hat t} - \hat\zeta_{\hat t} |^{\textcolor{black}{p}}\Big] \nonumber\\
&&
- {\textcolor{black}{2^{p-1}}} n\Big[\Upsilon_0\big(\ul{\hat \th},(\hat t, \ul{\hat \arm})\big)+ \Upsilon_0\big(\ul{\hat \eta},(\hat t, \ul{\hat \arm})\big)\Big]- \textcolor{black}{2^{p-1} n^{q}}\dbE\Big[|\hat\xi_{\hat t} - \hat \arm_{\hat t}|^{\textcolor{black}{p}} + |\hat\zeta_{\hat t} - \hat \arm_{\hat t}|^{\textcolor{black}{p}}\Big]\nonumber
 \\
 &=&
 n\dbE\big[\Upsilon_{\hat t}(\hat\xi -\hat \eta)\big]
 -  {\textcolor{black}{2^{p-1}}} n\dbE\Big[\Upsilon_{\hat t}\Big({\hat\xi -\hat \zeta\over 2}\Big)
                            +\Upsilon_{\hat t}\Big({\hat\zeta -\hat \xi\over 2}\Big) \Big]
 =
 0,
\eea
where the last equality is due to the fact that $\Upsilon_t({\bx\over 2}) = \Upsilon_t(-{\bx\over 2})  ={1\over {\textcolor{black}{2^{p}}}} \Upsilon_t(\bx)$, which can be verified straightforwardly from \reff{Upsilon0}. Then one can prove all the statements following the same arguments as in  Proposition \ref{prop-psine}. In particular, we have $\pi^{n,\e, N}_1(\ul\th) := \sum_{i=0}^\infty 2^{-i} \ol \Upsilon_0(\ul\th, \ul{\check \th}^N_i)$ and $\pi^{n,\e, N}_2(\ul\eta) := \sum_{i=0}^\infty 2^{-i} \ol \Upsilon_0(\ul\eta, \ul{\check \eta}^N_i)$, for some appropriate $\check \iota^N_i = (\check \th^N_i, \check \eta^N_i)$. Note that the $\pi^{n,\e, N}_j$ here includes the time difference, while the $\pi^{n,\e}_j$ in \reff{pi12} does not.
\qed

The convergence of $\check \iota^{N}$ is more involved and will be conducted in two cases.

\subsection{ Proof of Proposition \ref{prop-comparison1}: Case 1.}
\label{sect-case1}

In this subsection we prove the proposition in the case that:
\bea
\label{case1}
\mbox{ there is a subsequence, still denoted as $(\check\theta^{N},\check\eta^{N})_N$, satisfying $\check s^{N} \le \check t^{N}$.}
\eea
Surprisingly in this case we do not need to invoke the viscosity subsolution property of $U_1$. The viscosity supersolution property of $U_2$ alone induces the desired estimate.

 \begin{lem}
 \label{lem-psineN2}
In the setting of Proposition \ref{prop-psineN1}, if $\check s^{N} \le \check t^{N}$, then
\bea
\label{psineN2est}
&\dis\big\|\check\xi^{N}_{\cdot\wedge\check t^{N}}
          -\hat \xi_{\cdot\wedge\hat t}\big\|_{\textcolor{black}{p}}^{\textcolor{black}{p}}
 +\big\|\check\zeta^{N}_{\cdot\wedge\check s^{N}}
            -\hat \zeta_{\cdot\wedge\hat t}\big\|_{\textcolor{black}{p}}^{\textcolor{black}{p}}
+ |\check t^{N} - \hat t|^2 + |\check s^{N} - \hat t|^2 + N|\check t^{N} - \check s^{N}|^2\nonumber\\
&\dis \le  \rho_{C_\e}\big(C_{n,\e} N^{-{1\over 2}}\big) + C_{n,\e} N^{-{1\over 4}},
\eea
 where $C_\eps$ is defined in Proposition \ref{prop-psine}, and $C_{\eps,n}$ is a constant depending on $n,\eps$ only.
Consequently, under \eqref{case1}, there exists  $N_{n,\e}$,  such that
  \bea
  \label{TeN}
   \check s^{N}\le \check t^{N} \le T'_{\e, m} := {T_{\e, m}+T\over 2}<T,
&\mbox{for all}&
n\ge n_{\e,m}, \  N\ge N_{n,\e}.
   \eea
\end{lem}

\proof Clearly, \reff{TeN} is a direct consequence of the claimed estimate \reff{psineN2est}, which we now focus on. Fix $(n,\e, N)$ such that $\check s^{N}\le \check t^{N}$. We omit further the subscripts/superscripts $(n,\e)$ in $\psi^{n,\e}, \Psi^{n,\e}$, $\psi^N_{n,\e}, \Psi^N_{n,\e}$, $\pi_j^{n,\e}$, $\pi_j^{n,\e,N}$. Note that $\check s^N \le \check t^N$  implies
\bea\label{s<t}
  \Upsilon_{\check t^{N}}(\check \zeta^{N}_{\cd\wedge \check s^{N}}) =  \Upsilon_{\check s^{N}}(\check \zeta^{N}),\
\Upsilon_0\big(\ul{\check \th}^N, (\check t^N, \ul{\check\zeta}^N_{\cd\wedge \check s^N})\big) = \Upsilon_0\big(\ul{\check \th}^N, \ul{\check\eta}^N\big), \ \pi_2(\check t^{N}, \ul{\check \zeta}^{N}_{\cd\wedge \check s^{N}}) = \pi_2(\ul{\check \eta}^{N}),
\eea
where we used \reff{Upsilon0} and \reff{pi12}. Since $\Psi^N \le \psi^N$, by \reff{PsineNmax} and \reff{Psinemax} we have
\beaa
 0\ge \Psi^N(\ul{\check\iota}^{N}) - \psi^N(\ul{\check\iota}^{N}) \geq \Psi(\ul{\hat \l})- \psi^N(\ul{\check\iota}^{N})  \ge \Psi(\check{t}^{N}, \ul{\check\xi}^{N}, \ul{\check\zeta}^{N}_{\cdot\wedge\check{s}^N}) - \psi^N(\ul{\check\iota}^{N}).
 \eeaa
 Thus, it follows from the definitions of $\psi$, $\psi^N$ in \reff{psine} and \reff{psineN}, as well as  Proposition \ref{prop-psine} (ii) and \reff{s<t} that
\bea
\label{psineN-est0}
0&\geq& \psi(\check{t}^{N}, \ul{\check\xi}^{N}, \ul{\check\zeta}^{N}_{\cdot\wedge\check{s}^N}) - \big[\pi_0(\check t^N) + \pi_1(\ul{\check\th}^{N}) + \pi_2(\ul{\check\eta}^{N})\big] - \psi^N(\ul{\check\iota}^{N})\nonumber\\
&=&\Big[U_2(\check s^N,\ul{\check \zeta}^N) -U_2(\check t^N, \ul{\check\zeta}^N_{\cd\wedge \check s^N})\Big]+ \Big[\bar\Upsilon_0\big(\ul{\check \th}^N, \ul{\hat \th}\big) + \bar\Upsilon_0\big(\ul{\check\eta}^N, \ul{\hat \eta}\big)\Big] \nonumber\\
&&\dis  +{\e\over 2T} (\check t^N-  \check s^N)\dbE\big[\Upsilon_{\check s^N}(\check \zeta^N) \big]+\textcolor{black}{2^{p-1}}n\Big[\Upsilon_0(\ul{\check\th}^N,(\hat t, \ul{\hat \arm})) + \Upsilon_0(\ul{\check \eta}^N,(\hat t, \ul{\hat \arm}))\Big]\ms\nonumber\\
&&\dis - n\Upsilon_0(\ul{\check \th}^N, \ul{\check\eta}^N)- \textcolor{black}{n^{q}}\mathbb{E}\big[|\check\xi^N_{\check t^N}-\check\zeta^N_{\check s^N}|^{\textcolor{black}{p}}\big] + N|\check s^N-\check t^N|^2\ms\\
&&\dis + \sup_{\a\in \cA_{[\hat t, T]}}
                      \cJ^{\hat\th,\hat\arm}_-(\check\th^N,\a)
                      + \inf_{\a\in \cA_{[\hat t, T]}}
                         \cJ_+^{\hat\th,\hat\arm}(\check\eta^N,\a).
 \nonumber
\eea
Recall \reff{cIaxi}. Note that, by \reff{Upsilon0} and \reff{triangle}, as well as the fact \textcolor{black}{$2^{p-1}[|x|^{p} + |y|^{p}] \ge |x-y|^{p}$},
\bea
\label{psineN-est1}
&&\dis \!\!\!\!\!\!\!\!\!\!\!\!\!\!\!\!\!\!  \textcolor{black}{2^{p-1}}\Big[\Upsilon_0(\ul{\check\th}^N,(\hat t, \ul{\hat \arm})) + \Upsilon_0(\ul{\check \eta}^N,(\hat t, \ul{\hat \arm}))\Big]
\ge
 \Upsilon_0(\ul{\check \th}^N, \ul{\check\eta}^N),
\nonumber\\
&&\dis\!\!\!\!\!\!\!\!\!\!\!\!\!\!\!\!\!\! \sup_{\a\in \cA_{[\hat t, T]}} \!\!
                      \cJ^{\hat\th,\hat\arm}_-(\check\th^N,\a)
                      + \!\!\inf_{\a\in \cA_{[\hat t, T]}} \!\!
                         \cJ_+^{\hat\th,\hat\arm}(\check\eta^N,\a)
\ge
\inf_{\a\in \cA_{[\hat t, T]}} \Big[
                      \cJ^{\hat\th,\hat\arm}_-(\check\th^N,\a)
                      +\cJ_+^{\hat\th,\hat\arm}(\check\eta^N,\a)\Big]
\\
&&\qq\qq\qq\qq\ge
\inf_{\a\in \cA_{[\hat t, T]}} \Big[\textcolor{black}{n^{q}}\dbE\big[|\hat\cI^{\a}_{\check t^N}(\check\xi^N) -\hat \cI^{\a}_{\check s^N}(\check\zeta^N) |^{\textcolor{black}{p}} \big]+ \hat F^{\a}_{\check s^N} - \hat F^{\a}_{\check t^N}\Big].
\nonumber
\eea
Plug these into \reff{psineN-est0} and note that  $\check t^N-  \check s^N\ge 0$, we have
\bea
\label{psineN-est2}
&&  N|\check s^N\!-\!\check t^N|^2
+
\Big[\bar\Upsilon_0\big(\ul{\check \th}^N,  \ul{\hat \th}\big)
        + \bar\Upsilon_0\big(\ul{\check\eta}^N,  \ul{\hat \eta}\big)\Big]
\le
 U_2(\check t^N\!, \ul{\check\zeta}^N_{\cd\wedge \check s^N})\!-\! U_2(\check s^N\!,\ul{\check \zeta}^N)   \ms
\\
&&\qq + \textcolor{black}{n^{q}}\mathbb{E}\big[|\check\xi^N_{\check t^N}
                                                      - \check\zeta^N_{\check s^N}|^{\textcolor{black}{p}}\big] -\inf_{\a\in \cA_{[\hat t, T]}} \Big[\textcolor{black}{n^{q}}\dbE\big[|\hat\cI^{\a}_{\check t^N}(\check\xi^N) - \hat \cI^{\a}_{\check s^N}(\check\zeta^N) |^{\textcolor{black}{p}} \big]+ \hat F^{\a}_{\check s^N} - \hat F^{\a}_{\check t^N}\Big].
\nonumber
\eea
Recall \reff{barXtxi}, \reff{fne}, and denote, for any $\a$,
\beaa
\D^\a_N := \big[\check\xi^N_{\check t^N}-\check\zeta^N_{\check s^N}\big] - \big[\hat\cI^{\a}_{\check t^N}(\check\xi^N) -\hat \cI^{\a}_{\check s^N}(\check\zeta^N)\big]= \int_{\check s^N}^{\check t^N} b^{\hat \th, \a}_r  dr + \int_{\check s^N}^{\check t^N} \si^{\hat \th, \a}_r dB_r.
\eeaa
By \reff{CTe} we have $\dbE[|\D^\a_N|^{\textcolor{black}{p}}] \le C_\e(\check t^N-\check s^N)^{\textcolor{black}{p\over 2}}$, and then, by Proposition \ref{prop-psineN1} (iii),
  \beaa
  &&\dis \mathbb{E}\Big[|\check\xi^N_{\check t^N}-\check\zeta^N_{\check s^N}|^{\textcolor{black}{p}} \Big] - \mathbb{E}\Big[|\hat\cI^{\a}_{\check t^N}(\check\xi^N) -\hat \cI^{\a}_{\check s^N}(\check\zeta^N)|^{\textcolor{black}{p}}\Big] \\
  &&\dis =\mathbb{E}\Big[|\check\xi^N_{\check t^N}-\check\zeta^N_{\check s^N}|^{\textcolor{black}{p}} - |\check\xi^N_{\check t^N}-\check\zeta^N_{\check s^N}- \D_N^\a|^{\textcolor{black}{p}}\Big]  \le C_\e \|\D_N^\a\|_{\textcolor{black}{p}} \le C_\e (\check t^N-\check s^N)^{1\over 2},
  \eeaa
Moreover, since $-U_2$ satisfies \reff{Ureg}, it is clear that
  \beaa
  U_2(\check t^N, \ul{\check\zeta}^N_{\cd\wedge \check s^N})-U_2(\check s^N,\ul{\check \zeta}^N)\le \rho_{C_\e}\big(\check t^N-\check s^N\big);\q \big| F^{\hat\th,\a}_{\check s^N} - F^{\hat\th,\a}_{\check t^N} \big|\le C_\e (\check t^N-\check s^N).
  \eeaa
  Then by    \reff{psineN-est2} we have
  \bea
\label{psineN-est3}
\dis N|\check s^N-\check t^N|^2 + \big[\bar\Upsilon_0\big(\ul{\check \th}^N,  \ul{\hat \th}\big) + \bar\Upsilon_0\big(\ul{\check\eta}^N,  \ul{\hat \eta}\big)\big]  \le  \rho_{C_\e}\big(\check t^N-\check s^N\big) +C_{n,\e}(\check t^N-\check s^N)^{1\over 2}.
\eea
 As $\rho_R(\d) \le 1+C_R\d$, without loss of generality, the above implies that $0\le \check t^N-\check s^N \le C_{n,\e} N^{-{1\over 2}}$. Then \reff{psineN2est} follows by plugging this into the right side of \reff{psineN-est3}.
\qed

\vspace{3mm}
\no{\bf Proof of Proposition \ref{prop-comparison1} under \eqref{case1}.}\q
Fix $\e\le \e_m$, $n\ge n_{\e,m}$, $N\ge N_{n,\e}$ such that $\check s^N \le \check t^N$.  As before we omit the subscripts/supscripts $(n, \e)$. Introduce
\beaa
\f^N_2(\ul\eta) := \f^{\hat\eta}_2(\ul\eta, \ul{\hat\arm}) + N|s - \check t^N|^2 + \pi^N_2(\ul\eta); ~ \phi^N_2(\ul\eta) := \inf_{\a\in \cA_{[\hat t, T]}} \Big[\textcolor{black}{2^{p-1}n^{q}}\mathbb{E}\big[|\hat \cI^{\a}_s( \zeta)|^{\textcolor{black}{p}}  \big]
                         + \hat F^{\a}_t\Big].
\eeaa
Here, in terms of the notation in \reff{cIaxi} and \reff{C+cXp}, we are setting $(\tilde t, \tilde \xi) = (t', \xi') = (\hat t, \hat \xi)$, $(\tilde b^\a, \tilde \si^\a, \tilde f^\a) = (b^{\hat\th, \a}, \si^{\hat\th, \a}, f^{\hat\th, \a})$, and \textcolor{black}{$k=2^{p-1}n^{q}$}. Then one can easily check that $\phi^N_2 \in C^+(\textcolor{black}{\ol\cX^{p}}_{[\check s^N,T]})$.  By \reff{TeN} we have $\check s^N < T$, and by Proposition \ref{prop-psine} (ii) and Proposition \ref{prop-psineN1} (ii) we see that $\f^N_2\in C^{1,2}(\textcolor{black}{\ol\cX^{p}}_{\textcolor{black}{[\check s^N, T]}})$. Then it follows from Proposition \ref{prop-psineN1} (i) that  $(-\f^N_2, -\phi^N_2)\in \textcolor{black}{\mathfrak{F}^-_{p}} U_2(\ul{\check \eta}^N)$.
Thus, by the viscosity $p$-supersolution property of $U_2$, and recalling \reff{barXtxi}, we have
\beaa
c_1 &\le&
\pa_t \f^N_2(\ul{\check \eta}^N) + \liminf_{\d\to 0}\sup_{\a\in \cA_{[\check s^N, T]}} {1\over \d} \int_{\check s^N}^{\check s^N+\d} \Big[  \dot \phi^N_2\Big(s, \ul{\bar X}^{{\check \eta}^N, \a}\Big)
\nonumber\\
&&\hspace{40mm}
-  H_s\Big(\ul{\check \zeta}^N_{\cd\wedge \check s^N}, -  \pa_X \f^N_2(\ul{\check \eta}^N), -  \pa_{\bx X} \f^N_2(\ul{\check \eta}^N), \ul\a_s\Big) \Big]ds.
\eeaa
Note that $t'=\hat{t}\leq \check s^N$. This allows us to apply Proposition \ref{prop-C+cXp} (i) and obtain
\bea
\label{paphi2-est1}
 \int_{\check s^N}^{\check s^N+\d}  \dot \phi^N_2(s, \ul{\bar X}^{\check \eta^N, \a})ds  &\le& \inf_{\tilde\a\in \cA_{[\check s^N, T]}} \int_{\check s^N}^{\check s^N+\d} \Big[C \textcolor{black}{n^{q}} I^{\check \eta^N}_s(\a, \tilde \a)
+\  f^{\hat \th,\tilde\a}_s\Big]ds\nonumber\\
&\le& \int_{\check s^N}^{\check s^N+\d} \Big[C \textcolor{black}{n^{q}}  I^{\check \eta^N}_s(\a, \a)+  f^{\hat \th, \a}_s\Big]ds,
\eea
where we denoted
\bea
\label{Itxia}
&& I^{\eta}_s(\a, \tilde \a)
:=
\|b^{\eta, \a}_s- b^{\hat \th, \tilde\a}_s\|_{\textcolor{black}{p}}
\sup_{ \a'\in \cA_{[\hat t, T]}}\| \cI^{\hat \th,\a'}_s( \bar X^{\eta, \a})
\|_{\textcolor{black}{p}}^{\textcolor{black}{p-1}}\nonumber\\
&&\dis\qq\qq\qq+\|\si^{\eta, \a}_s- \si^{\hat \th, \tilde\a}_s\|_{\textcolor{black}{p}}^2
  \sup_{ \a'\in \cA_{[\hat t, T]}}\| \cI^{\hat \th,\a'}_s( \bar X^{\eta, \a})
  \|_{\textcolor{black}{p}}^{\textcolor{black}{p-2}},\\
&&  \cI^{\hat \th,\a'}_s( \bar X^{\eta, \a}):=\cI^{b^{\hat\th},\si^{\hat\th},\hat t,\hat \arm,\a'}_s( \bar X^{\eta, \a}).\nonumber
\eea

Observe that it is crucial to have $\inf_{\tilde \a}$ in the right side of the first line of \eqref{paphi2-est1}, which helps to replace $\tilde \a$ with $\a$ in the second line.  Then, by \reff{HJB} we have
\bea
\label{paphi2-est2}
c_1 &\le& K_0^N + \liminf_{\d\to 0}\sup_{\a\in \cA_{[\check s^N, T]}} {1\over \d} \int_{\check s^N}^{\check s^N+\d} [K^{N}_1(s) + K^{N}_2(s) + K^{N}_3(s) + K^{N}_4(s)]ds,
\eea
where
\bea
\label{KN0-4}
\left.\ba{lll}
\dis K^N_0 := \pa_t \f^N_2(\ul{\check \eta}^N),\q
K^N_1(s) :=\dbE\big[ b^{{\check \eta}^N,\a}_s\cd \pa_X \f^N_2(\ul{\check \eta}^N)\big],\ms\\
\dis K^N_2(s) := {1\over 2} \dbE\big[ (\si\si^\top)^{{\check \eta}^N,\a}_s: \pa_{\bx X} \f^N_2(\ul{\check \eta}^N)\big],\ms
\\
K^N_3(s) := C \textcolor{black}{n^{q}} I^{{\check \eta}^N}_s(\a, \a),
\q
\dis K^N_4(s) := \dbE\big[f^{\hat \th, \a}_s - f^{{\check \eta}^N,\a}_s\big].
\ea\right.
\eea

We now estimate the $K^N_i$ separately. We shall send $N\to \infty$, along the subsequence such that $\check s^N \le \check t^N$, by applying the convergence in Proposition \ref{prop-psineN1} (i) and Lemma \ref{lem-psineN2} repeatedly. We shall also apply the estimates \eqref{lem-Upsilon}, Proposition \ref{prop-psine} (ii), (iii), Proposition \ref{prop-psineN1} (ii), (iii), as well as Assumption \ref{assum-standing} repeatedly. At below we let $o(1)$ denote generic terms which will vanish when $N\to \infty$. The convergence rate may depend on $n, \e$, but is uniform in $\a, \d, s$. Note that, by \reff{fne}, \reff{Upsilon0}, and \reff{statesop3333},
  \bea
  \label{paf2}
  \left.\ba{lll}
\dis  \partial_t\varphi^{\hat\eta}_2(\ul\eta, \ul{\hat\arm})=-\frac{\varepsilon}{ 2T}\mathbb{E}[\Upsilon(\eta)]+2({s}-{\hat{t}}),\ms\\
\dis \partial_X\varphi^{\hat\eta}_2(\ul\eta, \ul{\hat\arm}) =\partial_{\bx}\Upsilon_s(\zeta-\hat{\zeta}_{\cdot\wedge{\hat{t}}})
                     +\textcolor{black}{2^{p-1}}n\partial_{\bx}\Upsilon_s(\zeta-\hat{\arm}_{\cdot\wedge{\hat{t}}})+\varepsilon[1-{s\over 2T}]
                     \partial_{\bx}\Upsilon_s(\zeta)+\pa_X \pi_2(\ul\eta),\ms\\
\dis \partial_{\bx X}\varphi^{\hat\eta}_2(\ul\eta, \ul{\hat\arm})=\partial_{\bx\bx}\Upsilon_s(\zeta- {\hat \zeta}_{\cdot\wedge{\hat{t}}})
           +\textcolor{black}{2^{p-1}}n\partial_{\bx\bx}\Upsilon_s(\zeta-\hat{\arm}_{\cdot\wedge{\hat{t}}})
                      \\
     \dis\qq\qq\qq +\e[1-{s\over 2T}]\partial_{\bx\bx}\Upsilon_s(\zeta)+\pa_{\bx X}\pi_2(\ul\eta).
 \ea\right.
  \eea

First, thanks to the crucial assumption $\check s^N \le \check t^N$,
\bea
\label{KN0est}
K^N_0 &=& \partial_t\varphi^{\hat\eta}_2(\ul{\check \eta}^N, \ul{\hat\arm}) +2N(\check s^N-\check{t}^N) + \pa_t \pi^N_2(\ul{\check \eta}^N)\nonumber\\
&\le& \partial_t\varphi^{\hat\eta}_2(\ul{\check \eta}^N, \ul{\hat\arm}) + o(1)=-\frac{\varepsilon}{ 2T}\mathbb{E}[\Upsilon(\check \eta^N)]+2(\check s^N-{\hat{t}})+o(1) \nonumber\\
&=& -\frac{\varepsilon}{ 2T}\mathbb{E}[\Upsilon(\hat \eta)] + o(1) \le -\frac{\varepsilon}{ 2T}\mathbb{E}\big[|\hat \zeta_{\cd\wedge \hat t}|_\infty^{\textcolor{black}{p}}\big] + o(1).
\eea
Next, denote $\|\zeta\|_{2,t} := \big(\dbE^0_t[|\zeta_{\cd\wedge t}|_\infty^2]\big)^{1\over 2}$,
\bea
\label{KN1est}
&&\dis K^N_1(s)
=
\dbE\Big[ b^{{\check \eta}^N,\a}_s\cd \big[\pa_X \f^{\hat\eta}_2(\ul{\check \eta}^N)+ \pa_X\pi^N_2(\ul{\check \eta}^N)\big]\Big]
=
\dbE\Big[ b^{{\check \eta}^N,\a}_s\cd \pa_X \f^{\hat\eta}_2(\ul{\check \eta}^N)\Big]  + o(1)\nonumber\\
&&\le
C\dbE\Big[\Big(1+|\check{\zeta}^N_{\cdot\wedge\check{s}^N}|_\infty+\|\check{\zeta}^N\|_{2,\check{s}^N}\Big)\times\nonumber\\
&&\q \Big(|\partial_{\bx}\!\Upsilon_{\check s^N}( \check \zeta^N-\hat{\zeta}_{\cdot\wedge{\hat{t}}})|+n|\partial_{\bx}\!\Upsilon_{\check s^N}( \check \zeta^N\!-\!\hat{\arm}_{\cdot\wedge{\hat{t}}})|+\varepsilon|\partial_{\bx}\!\Upsilon_{\check s^N}(\check \zeta^N)| + |\pa_X\!\pi_2(\ul{\check \eta}^N)|\Big)\Big] \!+\! o(1)
\nonumber\\
 &&=
 C\dbE\Big[\Big(1+|\hat{\zeta}_{\cdot\wedge\hat{t}}|_\infty+\|\hat{\zeta}\|_{2, \hat t}\Big)
  \Big(n|\partial_{\bx}\Upsilon_{\hat t}(\hat \zeta-\hat{\arm})| +\varepsilon|\partial_{\bx}\Upsilon_{\hat t}(\hat \zeta)|  +  |\pa_X\pi_2(\ul{\hat \eta})|\Big)\Big] + o(1)\nonumber\\
&&\le
C \varepsilon \dbE\big[1+|\hat \zeta_{\cdot\wedge\hat t}|_\infty^{\textcolor{black}{p}}\big] +C_{\e}n\big(\dbE\big[|\hat \zeta_{\hat t} - \hat\xi_{\hat t}|^{\textcolor{black}{p}}\big]\big)^{\textcolor{black}{p-1\over p}}+C_{\e}\big(\dbE\big[|\pa_X\pi_2(\ul{\hat \eta}) |^{\textcolor{black}{p\over p-1}}]\big)^{{\textcolor{black}{p-1\over p}}} + o(1)\nonumber\\
&&\le
C \varepsilon \dbE\big[1+|\hat \zeta_{\cdot\wedge\hat t}|_\infty^{\textcolor{black}{p}}\big] +C_{\e}{n^{1-\left(q+{\beta\over{p-\beta}}\right){{p-1}\over{p}}}}
+C_{\e} n^{-{\textcolor{black}{p-1\over p}}}+ o(1),
\eea
where the {last} inequality thanks to \reff{xietane} and \reff{pineest}. Similarly,
\bea
\label{KN2est}
K^N_2(s)
\!\!\!&\le&\!\!\!
C\dbE\Big[\Big(1+|\hat{\zeta}_{\cdot\wedge\hat{t}}|^2_\infty+\|\hat{\zeta}\|_{2,\hat t}^2\Big)\times\nonumber\\
&&\dis\qq  \Big(n|\partial_{\bx\bx}\Upsilon_{\hat t}(\hat \zeta-\hat{\arm})| +\varepsilon|\partial_{\bx\bx}\Upsilon_{\hat t}(\hat \zeta)|  +  |\pa_{\bx X}\pi_2(\ul{\hat \eta})|\Big)\Big] + o(1)\nonumber\\
  \!\!\!&\le&\!\!\! C \varepsilon \dbE\big[1+|\hat \zeta_{\cdot\wedge\hat t}|_\infty^{\textcolor{black}{p}}\big] +C_{\e}n\big(\dbE\big[|\hat \zeta_{\hat t} - \hat\xi_{\hat t}|^{\textcolor{black}{p}}\big]\big)^{{\textcolor{black}{p-2\over p}}}+C_{\e}\big(\dbE\big[|\pa_X\pi_2(\ul{\hat \eta}) |^{\textcolor{black}{p\over p-2}}]\big)^{{\textcolor{black}{p-2\over p}}} + o(1)\nonumber\\
 \!\!\!&\leq&\!\!\! C \varepsilon \dbE\big[1+|\hat \zeta_{\cdot\wedge\hat t}|_\infty^{\textcolor{black}{p}}\big] +C_{\e}{n^{1-\left(q+{\beta\over{p-\beta}}\right){{p-2}\over{p}}}}
+C_{\e} n^{-{\textcolor{black}{p-2\over p}}}+ o(1).
\eea
To estimate $K^N_3$, we first recall \reff{cIaxi} and note that
\beaa
 \cI^{\hat \th,\a'}_s( \bar X^{\check\eta^N, \a})
= \check \zeta^N_{\check s_N} - \hat \arm_{\hat t}
+ \int_{\check s^N}^s b^{\check \eta^N, \a}_t dt
                                  +  \int_{\check s^N}^s\si^{\check \eta^N, \a}_t  dB_t  - \int_{\hat t}^s b^{\hat \th,\a'}_t dt
                                 -  \int_{\hat t}^s \si^{\hat \th,\a'}_t dB_t.
\eeaa
Then, for  $s\in[\check s^N, \check s^N +\d]$,
we have
 \beaa
 \dbE\big[| \cI^{\hat \th,\a'}_s( \bar X^{\check\eta^N, \a})|^{\textcolor{black}{p}}\big]
 &\le&
 C \dbE\big[| \check \zeta^N_{\check s_N} - \hat \arm_{\hat t}|^{\textcolor{black}{p}}\big]
 + C_\e(s- \hat t)^{\textcolor{black}{p\over 2}}\\
 &=&
 C \dbE\big[| \hat \zeta_{\hat t} - \hat \xi_{\hat t}|^{\textcolor{black}{p}}\big] + C_\e \d^{\textcolor{black}{p\over 2}} + o(1)
 \;\le\;
 C_\e \textcolor{black}{n^{-{\beta\over p-\beta}-q}} + C_\e \d^{\textcolor{black}{p\over 2}} + o(1).
 \eeaa
Moreover, by \reff{xietane} we also have
\beaa
&&\dbE\Big[|b^{\check \eta^N, \a}_s- b^{\hat \th, \a}_s|^{\textcolor{black}{p}}+|\si^{\check \eta^N, \a}_s - \si^{\hat \th, \a}_s|^{\textcolor{black}{p}}\Big]\\
&&\le C\dbE\big[|\check \zeta^N_{\cd\wedge \check s^N} - \hat\xi_{\cd\wedge \hat t}|_\infty^{\textcolor{black}{p}}\big]  =  C \dbE\big[| \hat \zeta_{\cd\wedge \hat t} - \hat \xi_{\cd\wedge \hat t}|_\infty^{\textcolor{black}{p}}\big] + o(1) \le C_\e \textcolor{black}{n^{-1-{\beta\over p-\beta}}} + o(1),\\
&&\textcolor{black}{|f^{\check \eta^N, \a}_s - f^{\hat \th, \a}_s|^{p}
\le C\dbE\big[|\check \zeta^N_{\cd\wedge \check s^N} - \hat\xi_{\cd\wedge \hat t}|_\infty^{p}\big]^{\beta}  \le C_\e n^{-\beta-{\beta^2\over p-\beta}} + o(1)}.
\eeaa
Plug these into \reff{Itxia}, we have
\bea
\label{KN3est}
\left.\ba{lll}
K^N_3(s) \le C_\e n^{q} \Big[\big(n^{-{\beta\over \textcolor{black}{p}-\beta}-q} + \d^{\textcolor{black}{p\over 2}}\big)^{\textcolor{black}{p-1\over p}} n^{(-1-{\beta\over \textcolor{black}{p}-\beta}){1\over \textcolor{black}{p}}} \ms\\
\qq\qq+ \big(n^{-{\beta\over \textcolor{black}{p}-\beta}-q} + \d^{\textcolor{black}{p}\over 2} \big)^{\textcolor{black}{p-2\over p}} n^{(-1-{\beta\over \textcolor{black}{p}-\beta}){2\over \textcolor{black}{p}}}\Big] + o(1)\ms\\
\qq\q \le C_\e  {n^{{2q\over {p}}- {2\over p} - {\b \over p-\b}}} + C_{n,\e} \d^{\textcolor{black}{p\over 2}-1} + o(1);\ms\\
K^N_4(s) \le C_\e n^{-(\beta+{\beta^2\over \textcolor{black}{p}-\beta}){1\over \textcolor{black}{p}}} + o(1)=C_\e n^{-{\beta\over p-\b}} + o(1). 
\ea\right.
\eea
Plug the estimates \reff{KN0est}-\reff{KN3est}  into \reff{paphi2-est2},  we obtain
\bea
\label{c1contradiction1}
c_1 &\le& -\frac{\varepsilon}{ 2T}\mathbb{E}[|\hat \zeta_{\cd\wedge \hat t}|_\infty^{\textcolor{black}{p}}]  +C \varepsilon \dbE\big[1+|\hat \zeta_{\cdot\wedge\hat t}|_\infty^{\textcolor{black}{p}}\big] \nonumber\\
&&+C_{\e}\Big[ {{n^{1-\left(q+{\beta\over{p-\beta}}\right){{p-2}\over{p}}}}
+ n^{-{{p-2\over p}}}+ n^{{2q\over {p}}- {2\over p} - {\b \over p-\b}} + n^{-{\beta\over p-\b}}}\Big]+o(1)\nonumber\\
&\le&  \e \big[C-{1\over 2T}\big] \mathbb{E}[|\hat \zeta_{\cd\wedge \hat t}|_\infty^{\textcolor{black}{p}}] + C\e + C_\e {n^{-\tilde{q}}} + o(1),
\eea
 {where
\bea
\label{tildebeta}
\left.\ba{c}
\dis\tilde q:= \tilde q_1  \wedge {p-2\over p} \wedge \tilde q_2 \wedge {\beta\over p-\b},\ms\\
\dis \tilde q_1 := \left(q+{\beta\over{p-\beta}}\right){{p-2}\over{p}} - 1,\q \tilde q_2 := {2\over p} +{\b \over p-\b}- {2q\over {p}}.
\ea\right.
\eea
Recall \reff{pbeta2}, we have
\bea
\label{pbeta3}
\left.\ba{lll}
\dis \tilde q_1 >  \left({p\over p-2}-{\beta\over {p-\beta}}+{\beta\over{p-\beta}}\right){{p-2}\over{p}} - 1 =0;\ms\\
\dis \tilde q_2 > {2\over p} + {\b \over p-\b} - {2\over p}\big(1+{\beta p\over {2p-2\beta}}\Big)=0.
\ea\right.
\eea
Then $\tilde q>0$.}  Sending $N\to\infty$ in \reff{c1contradiction1}, we obtain \reff{comparison1} immediately, in the case \reff{case1}.
\qed

\subsection{ Proof of Proposition  \ref{prop-comparison1}: Case 2}
We now turn to the case that
\bea
\label{case2}
 \check t^{N} < \check s^{N}\q\mbox{for all large $N$}.
 \eea
In this case,  we need to modify \reff{psineN} by adding another penalization. Fix $0<\e\le \e_m$, $n\ge n_{\e,m}$, and $N$ large enough such that \reff{case2} holds.  Consider the settings of Propositions \ref{prop-psine} and \ref{prop-psineN1}. Again we omit the subscripts/superscripts $(n,\e)$ when the contexts are clear. For each $M$, recall \reff{fne} and define:
\bea
\psi^{N,M}(\ul\iota)
&:=&
U_1({\ul\th})-U_2({\ul\eta})
- \big[\f^{\hat\th}_1(\ul\th,\ul{\hat\arm})
         + \f^{\hat\eta}_2(\ul\eta,\ul{\hat\arm}) \big] - \pi^{n,\e}_0(t) - N|s-t|^2
\nonumber\\
&&
- M\big|(s-t)-(\check s^{N}-\check t^{N}) \big|^2 -  \inf_{\a\in \cA_{[\hat t, T]}}\Big(
                 \cJ_-^{\hat \th,\hat\arm}(\theta,\alpha)
        + \cJ_+^{\hat \th,\hat\arm}(\eta,\alpha)
 \Big).~~~~
 \label{psineNM}
\eea

 \begin{prop}
 \label{prop-psineNM}
In the setting of Proposition \ref{prop-psineN1}, for any $0< \e\le \e_m$, $n\ge n_{\e, m}$, and $N$ large enough such that $\check t^{N} < \check s^{N}$, the following hold.

\ms
\no{\rm (i)} There exist   $\grave\iota^{N,M}=(\grave \th^{N,M}, \grave \eta^{N,M})=\big((\grave t^{N,M}, \grave \xi^{N,M}), (\grave s^{N,M}, \grave  \zeta^{N,M})\big) \in \L_2(\check t^{N},\check s^{N})$ (omitting $^{n,\e}$ again ) and $\Psi^{N,M} \in C^0\big(\L_2(\grave t^{N,M},\grave s^{N,M})\big)$ such that
$\Psi^{N,M} \le \psi^{N,M}$ and
\bea
\label{PsineNMmax}
 \Psi^{N,M}(\ul{\grave \iota}^{N,M})
=\max_{ \L_2(\grave t^{N,M},\grave s^{N,M})} \,\Psi^{N,M}
\;\ge\; \Big(\sup_{ \L_2(\grave t^{N,M},\grave s^{N,M})} \,\psi^{N,M}-\frac1N\Big) \vee \psi^{N,M}(\ul{\check\iota}^N).
\eea

\no{\rm (ii)} Moreover, $[\psi^{N,M} - \Psi^{N,M}](\ul\iota)= \pi^{N,M}_1(\ul{\th}) + \pi^{N,M}_2(\ul{\eta})$, for some nonnegative functions $\pi^{N,M}_1\in C^{1,2}(\ol\cX^{\textcolor{black}{p}}_{[\grave t^{N,M}, T]})$, $\pi^{N,M}_2\in C^{1,2}(\ol\cX^{\textcolor{black}{p}}_{[\grave s^{N,M}, T]})$, such that: for $j=1,2$ and at $\grave\iota^{N,M}$,
\bea
\label{pineNMest}
\dis 0\le \pi^{N,M}_j\le {C\over N}; |\pa_t \pi^{N,M}_j| \le {C\over \sqrt{N}};  \|\pa_X \pi^{N,M}_j\|_{\textcolor{black}{p\over p-1}} \le {C\over N^{\textcolor{black}{p-1\over p}}}, \|\pa_{\bx X}  \pi^{N,M}_j\|_{\textcolor{black}{p\over p-2}} \le {C\over N^{\textcolor{black}{p-2\over p}}}.
\eea

\no {\rm (iii)} $\sup_{n,N,M}\mathbb{E}\Big[\|\grave \xi^{N,M}_{\cdot\wedge\grave t^{N,M}}\|^{\textcolor{black}{p}}+\|\grave \zeta^{N,M}_{\cdot\wedge\grave s^{N,M}}\|^{\textcolor{black}{p}}\Big] \le C_\e,$ for some $C_\e$, which may depend on $\e$.

\ms
\no {\rm (iv)}  $M\big|(\grave s^{N, M}- \grave t^{N, M})
                         -(\check s^{N}-\check t^{N})\big|^2 \le C_{n,\e}$, for some $C_{n,\e}$ which may depend on $n, \e$, but not on $N,M$; Consequently, there exists $M_N:=M_{n,\e,N}$ large enough  such that $\grave t^{N, M_N} < \grave s^{N, M_N}$.

\ms
\no {\rm (v)} For the $M_N$ in (iv), and abbreviating the notations further that $\grave \iota^N := \grave \iota^{N, M_N}$, we have
\beaa
&&
\mathbb{E}\Big[|\grave\xi^{N}_{\cdot\wedge\grave t^{N}}
-\hat \xi_{\cdot\wedge\hat t}|_\infty^{\textcolor{black}{p}} +|\grave\zeta^{N}_{\cdot\wedge\grave s^{N}}
-\hat \zeta_{\cdot\wedge\hat t}|_\infty^{\textcolor{black}{p}}\Big]
+ |\grave t^{N}- \hat t|^2 + |\grave s^{N} - \hat t|^2
\\ &&
\qq \qq \qq\qq
+ N|\grave t^{N} - \grave s^{N}|^2+M_N\Big|(\grave s^{N}- \grave t^{N})-(\check s^{N}-\check t^{N})\Big|^2
\;\le\;
C_{n,\e} N^{-\textcolor{black}{\beta\over 4-\beta}},
\eeaa
for some $C_{n,\e}$, which may depend on $n, \e$, but not on $N$. Consequently, there exists $N_{n,\e}$ such that $\grave t^{N} < \grave s^{N}<T$ for all $N\geq N_{n,\e}$.
\end{prop}

\proof  (i)-(iii) follow similar arguments as in Proposition \ref{prop-psine}. We now fix $(n,\e, N, M)$ and as in the previous proofs,  we omit the subscripts/superscripts $(n, \e)$, but keep $N, M$.

\ms

 \no (iv) First, since $\Psi^{N,M}\le \psi^{N,M}$ and $\Psi^N \le \psi^N$, we have
 \bea
 \label{psiNMcompare}
\psi^{N,M}(\ul{\grave \iota}^{N,M}) \ge \Psi^{N,M}(\ul{\grave \iota}^{N,M}) \ge \psi^{N,M}(\ul{\check \iota}^N) \ge  \psi^N(\ul{\check \iota}^N) \ge \Psi^N(\ul{\check \iota}^N) \ge \Psi^N(\ul{\grave \iota}^{N,M}),
\eea
where the second inequality is due to \reff{PsineNMmax},  the third one is by comparing \reff{psineN} and \reff{psineNM} directly, and the last one is due to \reff{PsineNmax}.
Then, by \reff{psineN}, \reff{psineNM}, and Proposition \ref{prop-psineN1} (ii),
\beaa
0 &\ge&  \Psi^N(\ul{\grave \iota}^{N,M}) -  \psi^{N,M}(\ul{\grave \iota}^{N,M})\\
&=&\inf_{\a\in \cA_{[\hat t, T]}}\Big(
                 \cJ_-^{\hat \th,\hat\arm}(\grave \th^{N,M},\alpha)
        + \cJ_+^{\hat \th,\hat\arm}(\grave \eta^{N,M},\alpha)
 \Big)+ M\big|(\grave s^{N,M} - \grave t^{N,M})-(\check s^{N}-\check t^{N}) \big|^2\\
 && - \Big[\sup_{\a\in \cA_{[\hat t, T]}}
                 \cJ_-^{\hat \th,\hat\arm}(\grave \th^{N,M},\alpha) + \inf_{\a\in \cA_{[\hat t, T]}} \cJ_+^{\hat \th,\hat\arm}(\grave \eta^{N,M},\alpha)
  + \pi^N_1(\grave \th^{N,M}) + \pi^N_2(\grave \eta^{N,M})\Big].
 \eeaa
 This, together with (iii) and \reff{pineNest}, implies the estimate in Item (iv) immediately.

\ms
\no (v) Recall again that $\grave \iota^N = \grave \iota^{N, M_N}$. By \reff{psiNMcompare}, \reff{PsineNmax}, and then \reff{Psinemax}, we have
\beaa
&&\psi^{N, M_N}(\ul{\grave \iota}^{N}) \ge \psi^{N}(\ul{\check \iota}^N) \ge \Psi(\ul{\hat\l}) \ge \Psi(\grave s^{N}, \ul{\grave \xi}^{N}_{\cd\wedge \grave t^{N}}, \ul{\grave \zeta}^{N}),
\eeaa
where $\Psi = \Psi^{n,\e}$. Since $\grave t^{N} < \grave s^{N}$, similarly to \reff{s<t} we have
\beaa
\dis\Upsilon_{\grave s^{N}}(\grave \xi^{N}_{\cd\wedge \grave t^{N}}) =  \Upsilon_{\grave t^{N}}(\grave \xi^{N}),~
\Upsilon_0\big(\ul{\grave \eta}^N, (\grave s^{N}, \ul{\grave \xi}^{N}_{\cd\wedge \grave t^{N}})\big) = \Upsilon_0\big(\ul{\grave \eta}^N, \ul{\grave\th}^N\big),~
\pi_1(\grave s^{N}, \ul{\grave \xi}^{N}_{\cd\wedge \grave t^{N}}) =\pi_1(\ul{\grave \th}^{N}).
\eeaa
Then, by comparing  \reff{psine} and \reff{psineNM}, and recalling Proposition  \ref{prop-psine} (ii) and \reff{fne}, we have
\beaa
0 &\ge& \psi(\grave s^{N}, \ul{\grave \xi}^{N}_{\cd\wedge \grave t^{N}}, \ul{\grave \zeta}^{N}) - \pi_0(\grave s^{N}) - \pi_1(\grave s^{N}, \ul{\grave \xi}^{N}_{\cd\wedge \grave t^{N}}) - \pi_2(\ul{\grave \eta}^{N}) - \psi^{N, M_N}(\ul{\grave \iota}^{N}) \\
&=& \big[U_1(\grave s^{N}, \ul{\grave \xi}^{N}_{\cd\wedge \grave t^{N}}) - U_1(\ul{\grave \th}^N)\big]  + \big[\bar\Upsilon_0\big(\ul{\grave \th}^N, \ul{\hat\th})\big) +\bar\Upsilon_0\big(\ul{\grave \eta}^N, \ul{\hat\eta})\big) \big]\\
&& +{\textcolor{black}{2^{p-1}}} n\big[\Upsilon_0(\ul{\grave \th}^N,(\hat t, \ul{\hat \arm})) +\Upsilon_0(\ul{\grave \eta}^N,(\hat t, \ul{\hat \arm}))
\big] - n\Upsilon_0(\ul{\grave \th}^{N},\ul{\grave \eta}^{N})- \textcolor{black}{n^{q}}\mathbb{E}\big[|\grave \xi^{N}_{\grave t^{N}}-\grave \zeta^{N}_{\grave s^{N}}|^{\textcolor{black}{p}}\big]\\
 &&+ \pi_0(\grave t^N) - \pi_0(\grave s^{N})+ N|\grave s^N-\grave t^N|^2+ M_N\big|(\grave s^N-\grave t^N)-(\check s^{N}-\check t^{N}) \big|^2\\
&&  + \inf_{\a\in \cA_{[\hat t, T]}}\Big(
                 \cJ_-^{\hat \th,\hat\arm}(\ul{\grave \th}^N,\alpha)
         + \cJ_+^{\hat \th,\hat\arm}(\ul{\grave \eta}^N,\alpha)\Big)\\
        &\ge& \big[U_1(\grave s^{N}, \ul{\grave \xi}^{N}_{\cd\wedge \grave t^{N}}) - U_1(\ul{\grave \th}^N)\big]  + \big[\bar\Upsilon_0\big(\ul{\grave \th}^N, \ul{\hat\th}\big) +\bar\Upsilon_0\big(\ul{\grave \eta}^N, \ul{\hat\eta}\big) \big]- \textcolor{black}{n^{q}}\mathbb{E}\big[|\grave \xi^{N}_{\grave t^{N}}-\grave \zeta^{N}_{\grave s^{N}}|^{\textcolor{black}{p}}\big]\\
 &&+ \pi_0(\grave t^N) - \pi_0(\grave s^{N})+ N|\grave s^N-\grave t^N|^2+ M_N\big|(\grave s^N-\grave t^N)-(\check s^{N}-\check t^{N}) \big|^2\\
&&  + \inf_{\a\in \cA_{[\hat t, T]}} \Big[\textcolor{black}{n^{q}}\dbE\big[|\hat\cI^{\a}_{\grave t^N}(\grave\xi^N) -\hat \cI^{\a}_{\grave s^N}(\grave\zeta^N) |^{\textcolor{black}{p}} \big]+ \big[\hat F^{\a}_{\grave s^N} - \hat F^{\a}_{\grave t^N}\big]\Big],
\eeaa
where the last inequality is due to estimates similar to \reff{psineN-est1}. Then,
\beaa
&&N|\grave s^N-\grave t^N|^2 + M_N\big|(\grave s^N-\grave t^N)-(\check s^{N}-\check t^{N}) \big|^2+ \bar\Upsilon_0\big(\ul{\grave \th}^N, \ul{\hat\th}\big) +\bar\Upsilon_0\big(\ul{\grave \eta}^N, \ul{\hat\eta}\big)\nonumber\\
&&
\hspace{25mm}
\le \big|U_1(\grave s^N, \ul{\grave \xi}^N_{\cd\wedge \grave t^N})-U_1(\grave t^N,\ul{\grave\xi}^N)\big|+  |\pi_0(\grave t^N) - \pi_0(\grave s^{N})| \nonumber\\
&&\dis
\hspace{29mm}-\inf_{\a\in \cA_{[\hat t, T]}} \Big[\textcolor{black}{n^{q}}\mathbb{E}\big[|\hat\cI^{\a}_{\grave t^N}(\grave\xi^N) -\hat \cI^{\a}_{\grave s^N}(\grave\zeta^N)|^{\textcolor{black}{p}} - |\grave\xi^N_{\grave t^N} -\grave\zeta^N_{\grave s^N}|^{\textcolor{black}{p}}\big] + \big[\hat F^{\a}_{\grave s^N} - \hat F^{\a}_{\grave t^N}\big]\Big]\\
&&
\hspace{25mm}
\le C_{n,\e}|\grave s^N-\grave t^N|^{\textcolor{black}{\beta\over 2}},
\eeaa
where the last inequality is due to \reff{pineest} and similar arguments as for \reff{psineN-est3}.
In particular, $N|\grave s^N-\grave t^N|^2 \le C_{n,\e}|\grave s^N-\grave t^N|^{\textcolor{black}{\beta\over2}}$, and thus $|\grave s^N-\grave t^N|\le C_{n,\e}N^{-\textcolor{black}{2\over 4-\beta}}$. Plugging this into the right side of the last inequality provides the required estimate in Item (v).
\qed

\bs
\no{\bf Proof of Proposition  \ref{prop-comparison1} in Case \reff{case2}.}
Fix $0<\e\le \e_m$, $n\ge n_{\e, m}$, $N\ge N_{n,\e}$, and recall the $M_N := M_{n,\e, N}$ in Proposition \ref{prop-psineNM} (iv) and (v).  As before we omit the subscripts/supscripts $(n, \e)$, and abbreviate $\grave \iota^N:= \grave \iota^{N,M_N}$. In particular, by otherwise considering a larger $N_{n,\e}$, we have $\grave t^N < \grave s^N < T$.
Introduce:
\beaa
\f^N_1(\ul\th) &:=& \f^{\hat\th}_1(\ul\th,\ul{\hat\arm}) + \pi_0(t) + N|t- \grave s^{N}|^2+M_N\big| t-\check t^{N} - \grave s^{N} + \check s^{N}\big|^2 + \pi^{N, M_N}_1(\ul\th);\\
\phi^N_1(\ul\th) &:=& \inf_{\a\in \cA_{[\hat t, T]}} \Big( \cJ_-^{\hat \th,\hat\arm}(\theta,\alpha)+\k_1(\a)\Big),\q \k_1(\a):=  \cJ_+^{\hat \th,\hat\arm}(\grave\eta^N,\alpha);\\
\f^N_2(\ul\eta) &:=& \f^{\hat\eta}_2(\ul\eta,\ul{\hat\arm}) + N|s- \grave t^{N}|^2+M_N\big|s-\check s^{N} - \grave t^{N} + \check t^{N}\big|^2 + \pi^{N, M_N}_2(\ul\eta);\\
 \phi^N_2(\ul\eta) &:=& \inf_{\a\in \cA_{[\hat t, T]}}  \Big(\cJ_+^{\hat \th,\hat\arm}(\eta,\alpha) +\k_2(\a)\Big),\q \k_2(\a):= \cJ_-^{\hat \th,\hat\arm}(\grave\theta^N,\alpha).
 \eeaa
  One can easily see that
$
 (\f^N_1, \phi^N_1)\in \textcolor{black}{\mathfrak{F}^+_{p}} U_1(\ul{\grave \th}^N)$ and $(-\f^N_2, -\phi^N_2)\in \textcolor{black}{\mathfrak{F}^-_{p}} U_2( \ul{\grave \eta}^N).
$
  Then, by the $p$-viscosity properties of $U_1$ and  $U_2$, we have
\bea
\label{I12}
\left.\ba{c}
\dis c_1 \le  K^N_0 + \liminf_{\d\to 0}\inf_{\a\in \cA_{[\grave t^N, T]}} {1\over \d} \int_{\grave t^N}^{\grave t^N+\d}  \big[K^N_1(\d, \a, t) + K^N_2(\d, \a, t)\big] dt \\
\dis +  \liminf_{\d\to 0}\sup_{\a\in \cA_{[\grave s^N, T]}} {1\over \d}  \int_{\grave s^N}^{\grave s^N+\d}   \big[K^N_3(\d, \a, s) + K^N_4(\d, \a, s)\big] ds,
\ea\right.
\eea
where
\bea
\label{KNd}
\left.\ba{lll}
\dis K^N_0 := \pa_t \f^N_1(\ul{\grave \th}^N) +\pa_t \f^N_2(\ul{\grave \eta}^N);\ms\\
\dis K^N_1(\d, \a, t) := \dbE\Big[b^{\grave \th^N,\a}_t\cd \pa_X \f^N_1(\ul{\grave \th}^N) + {1\over 2}(\si\si^\top)^{\grave \th^N,\a}_t: \pa_{\bx X} \f^N_1(\ul{\grave \th}^N)\Big];\ms\\
\dis K^N_2(\d, \a, t) := \dot \phi^N_1(t, \bar X^{\grave \th^N, \a}) + f^{\grave \th^N,\a}_t;\ms\\
\dis K^N_3(\d, \a, s) := \dbE\Big[b^{\grave \eta^N,\a}_s\cd \pa_X \f^N_2(\ul{\grave \eta}^N) + {1\over 2}(\si\si^\top)^{ \grave \eta^N,\a}_s: \pa_{\bx X} \f^N_2(\ul{\grave \eta}^N)\Big];\ms\\
\dis K^N_4(\d, \a, s) := \dot \phi^N_2(s, \bar X^{\grave \eta^N, \a}) - f^{\grave \eta^N,\a}_s.
\ea\right.
\eea

We next estimate $K^N_i$ separately. Similarly to the previous subsection, we shall send $N\to \infty$, under \reff{case2}, by applying the convergence in Proposition \ref{prop-psineNM} (ii), (v) repeatedly. We shall also apply the estimates \eqref{lem-Upsilon}, Proposition \ref{prop-psine} (ii), (iii), Proposition \ref{prop-psineN1} (ii), (iii), Proposition \ref{prop-psineNM} (v), as well as Assumption \ref{assum-standing} repeatedly. At below again we let $o(1)$ denote generic terms which will vanish when $N\to \infty$. The convergence rate may depend on $n, \e$, but is uniform in $\a, \d, t, s$. Recall \reff{paf2}, and we have similar expressions for the derivatives of $\f^{\hat\th}_1(\ul\th, \ul{\hat\arm})$.

First, noting the obvious but crucial cancellations, we have
\bea
\label{KN0est2}
K^N_0 &=& \pa_t \f^{\hat\th}_1(\ul{\grave \th}^N,\ul{\hat\arm}) + \pa_t\pi_0(\grave t^N)  + \pa_t \pi^{N, M_N}_1(\ul{\grave \th}^N)+\pa_t \f^{\hat\eta}_2(\ul{\grave \eta}^N,\ul{\hat\arm}) +  \pa_t \pi^{N, M_N}_2(\ul{\grave \eta}^N)\nonumber\\
&=&-\frac{\varepsilon}{2T}\mathbb{E}\big[\Upsilon(\grave \eta^N)\big] +2\big[(\grave t^N-{\hat{t}})+(\grave s^N-{\hat{t}})\big] + \pa_t\pi_0(\grave t^N) + o(1)\nonumber\\
&=&-\frac{\varepsilon}{2T}\mathbb{E}\big[ \Upsilon(\hat\zeta)\big]  + \pa_t\pi_0(\hat t) + o(1)\nonumber\\
&\le& -\frac{\varepsilon}{2T}\mathbb{E}\big[  |\hat \zeta_{\cd\wedge \hat t}|_\infty^{\textcolor{black}{p}}\big]  + Cn^{-{1\over 2}} + o(1).
\eea
Next, similarly to \reff{KN1est} and \reff{KN2est}, we have
\bea
\label{KN1est2}
\dis  K^N_1(\d, \a, t) \le C_{\e}\textcolor{black}{n^{-\tilde{q}}}+ o(1);\q
\dis K^N_3(\d, \a, s) \le C \varepsilon \dbE\big[1+|\hat \zeta_{\cdot\wedge\hat t}|_\infty^{\textcolor{black}{p}}\big] +C_{\e}\textcolor{black}{n^{-\tilde{q}}}+ o(1).
\eea
Moreover, note that, in terms of the notations in Proposition \ref{prop-C+cXp}, for $\phi^N_2$ we have $t' = \grave t^N$ and $t = \grave s^N$. Since $\grave t^N< \grave s^N$, we may apply Proposition \ref{prop-C+cXp} (i) again to estimate $K^N_4$: for any $\a$ and recalling \reff{Itxia},
\beaa
\int_{\grave s^N}^{\grave s^N+\d}   K^N_4(\d, \a, s) ds &\le& \inf_{\tilde\a\in \cA_{[t, T]}}
\int_{\grave s^N}^{\grave s^N+\d} \big[C\textcolor{black}{n^{q}}I^{\grave \eta^N}_s(\a, \tilde \a)+ f^{\hat\th,\tilde\a}_s\big]ds - \int_{\grave s^N}^{\grave s^N+\d}  f^{\grave \eta^N,\a}_sds\\
&\le& \int_{\grave s^N}^{\grave s^N+\d}  \big[C\textcolor{black}{n^{q}}I^{\grave \eta^N}_s(\a, \a)+ f^{\hat\th,\a}_s- f^{\grave \eta^N,\a}_s\big]ds.
\eeaa
Then it follows from the same arguments as in \reff{KN3est} that
$
{1\over \d}\int_{\grave s^N}^{\grave s^N+\d}  K^N_4(\d, \a, s) ds \le C_\e \textcolor{black}{n^{-\tilde{q}}} + C_n \textcolor{black}{\d^{{p\over 2}-1}} + o(1),
$
and thus
\bea
\label{KN4est2}
\sup_{\a\in \cA_{[\grave s^N, T]}} {1\over \d}  \int_{\grave s^N}^{\grave s^N+\d}  K^N_4(\d, \a, s)ds \le C \textcolor{black}{n^{-\tilde{q}}} + C_n \textcolor{black}{\d^{{p\over 2}-1}} + o(1).
\eea

Finally we estimate $K^N_2$. Recall the notations $K^N_2$ in \reff{KNd}, $\bar X^{\grave \th^N, \a}$ in \reff{barXtxi}, and the setting in \reff{cIaxi}-\reff{C+cXp}. Note that for $\phi^N_1$ we have $t' = \grave s^N > \grave t^N=t$, then we can only apply Proposition \ref{prop-C+cXp} (ii) here.  That is, there exists $\a^\d\in \cA_{[\grave t^N, T]}$, which may depend on $\ul X^{\grave \th^N,\a}_{\cd\wedge \grave t^N}= \grave \xi^N_{\cd\wedge \grave t^N}$, but not on $(b^{\grave \th^N,\a}_s, \si^{\grave \th^N,\a}_s)_{s\ge \grave t^N}$, and thus is independent of $\a$, such that
\beaa
\int_{\grave t^N}^{\grave t^N+\d}  K^N_2(\d, \a, t) dt \le
\int_{\grave t^N}^{\grave t^N+\d}  \Big[C\textcolor{black}{n^{q}}I^{\grave \th^N}_t(\a, \a^\d)- f^{\hat\th,\a^\d}_t+  f^{\grave \th^N,\a}_t\Big]dt+\textcolor{black}{\d^2},
\eeaa
and thus
\beaa
\inf_{\a\in \cA_{[\grave t^N, T]}}\int_{\grave t^N}^{\grave t^N+\d}  K^N_2(\d, \a, t) dt \le \int_{\grave t^N}^{\grave t^N+\d}  \Big[C\textcolor{black}{n^{q}}I^{\grave \th^N}_t(\a^\d, \a^\d)- f^{\hat\th,\a^\d}_t+  f^{\grave \th^N,\a^\d}_t\Big]dt+\textcolor{black}{\d^2}.
\eeaa
Then again by \reff{KN3est}  we have
\bea
\label{KN2est2}
\inf_{\a\in \cA_{[\grave t^N, T]}}{1\over \d}\int_{\grave t^N}^{\grave t^N+\d}  K^N_2(\d, \a, t) dt \le C \textcolor{black}{n^{-\tilde{q}}}+\textcolor{black}{\d} + C_n \textcolor{black}{\d^{{p\over 2}-1}} + o(1).
\eea

Plug \reff{KN0est2}, \reff{KN1est2}, \reff{KN4est2}, and \reff{KN2est2} into \reff{I12}, we obtain
\beaa
c_1 \le \e[C-{1\over 2T}]\mathbb{E}\big[ |\hat \zeta_{\cd\wedge \hat t}|_\infty^{\textcolor{black}{p}} \big] + C\e+C_{\e}\textcolor{black}{n^{-\tilde{q}}} + o(1).
\eeaa
Send $N\to \infty$, we obtain \reff{comparison1} immediately, in the case \reff{case2}. This, together with the proof in Subsection \ref{sect-case1},  completes the proof of Proposition  \ref{prop-comparison1}.
\qed

\begin{rem}
\label{rem-comparison2}
We remark that the arguments in this subsection do not work when $\check s^N < \check t^N$, and thus we cannot eliminate Subsection \ref{sect-case1}. Indeed, in this case   we can only apply Proposition \ref{prop-C+cXp} (ii), instead of (i), to estimate $\dot\phi^N_2$ and $K^N_4$. That is, there exists $\a^\d$, independent of $\a$, s.t.
\beaa
\int_{\grave s^N}^{\grave s^N+\d} K^N_4(\d,\a,s)ds\le  \int_{\grave s^N}^{\grave s^N+\d} \big[C \textcolor{black}{n^{q}} I^{\grave \eta^N}_s(\a, \a^\d) +  f^{\hat \th,\a^\d}_s - f^{\grave \eta^N,\a}_s\big]ds+\d^2.
\eeaa
However, for the term $\int_{\grave s^N}^{\grave s^N+\d}  K^N_4(\d,\a,s)ds$, we need to take supremum over $\a$:
\beaa
\sup_{\a\in \cA_{[\grave s^N, T]}}\int_{\grave s^N}^{\grave s^N+\d} \!\!\!\!\! K^N_4(\d,\a,s)ds\le \sup_{\a\in \cA_{[\grave s^N, T]}}  \int_{\grave s^N}^{\grave s^N+\d} \!\!\!\!\!\big[C \textcolor{black}{n^{q}} I^{\grave \eta^N}_s(\a, \a^\d) +  f^{\hat \th,\a^\d}_s - f^{\grave \eta^N,\a}_s\big]ds+\d^2.
\eeaa
Recall \reff{Itxia}, the right side of the above  won't be small anymore, and then we won't be able to derive the desired contradiction.
\end{rem}

\section{The comparison result: Proof of Theorem \ref{thm-comparison}}
\label{sect-comparison3}


\no{\bf (i) We first prove the case that $U_1$ satisfies \reff{Vreg1}.} Fix $\d_0\in(0, {1\over 2C})$ for the $C$ in  \reff{comparison1}. We proceed in three steps.

\ms
\no {\it Step 1.} We first assume \reff{U2c} holds true and $T\le \d_0$. Assume by contradiction that $m>0$.
Then, since $C<{1\over 2\d_0} \le {1\over 2T}$, by \reff{comparison1} we have
\beaa
c_1 \le C\e+C_\e \textcolor{black}{n^{-\tilde{q}}}, \q \mbox{for} \ n\geq n_{\e,m}.
\eeaa
By first sending $n\to\infty$ and then $\e\to 0$, we obtain the desired contradiction: $c_1\le 0$.

\ms
\no {\it Step 2.} We next assume \reff{U2c} holds true but $T$ can be arbitrarily large. Consider a partition $0=T_0<\cds<T_k =T$ such that $T_{i+1}-T_i\le \d_0$ for all $i$. By Step 1 we should have $U_1(t, \ul\xi)\le U_2(t, \ul\xi)$ for all $(t, \xi) \in \ol\cX^2_{[T_{k-1}, T_k]}$. In particular, $U_1(T_{k-1},\cd) \le U_2(T_{k-1}, \cd)$. Now consider the equation on $[T_{k-2}, T_{k-1}]$, by Step 1 again we obtain $U_1(t, \ul\xi)\le U_2(t, \ul\xi)$ for all $(t, \xi) \in \ol\cX^2_{[T_{k-2}, T_{k-1}]}$. Repeat the arguments backwardly in time, we prove the result on the whole interval $[0, T]$.

\ms
\no {\it Step 3.} In the general case, let $c_1>0$ be an arbitrary number. Set $U_2^{c_1} := U_2 + {(1+T)^2\over 1+t} c_1$. It is clear that $-U_2^{c_1}$ also satisfies \reff{Ureg} and, providing $U_2$ is smooth,
$$
\pa_t U_2^{c_1} = \pa_t U_2 - {(1+T)^2\over (1+t)^2} c_1\le  \pa_t U_2 - c_1,\q \pa_X U^{c_1}_2 = \pa_X U_2,\q \pa_{\bx X} U^{c_1}_2 = \pa_{\bx X} U_2.
$$
Then one can easily see that $U_2^{c_1}$ is a viscosity supersolution of $\cL U^{c_1}_2 \le -c_1$. Since $U_1(T,\cd) \le U_2(T, \cd) < U^{c_1}_2(T,\cd)$, then by Step 2 we have $U_1(t,\cd) \le U^{c_1}_2(t,\cd)$ for all $t\in [0, T]$. Now since $c_1>0$ is arbitrary, we see that $U_1(t,\cd) \le U_2(t,\cd)$ for all $t\in [0, T]$.

\ms
\no{\bf (ii) We next prove the case that $U_2$ satisfies \reff{Vreg1}.} The proof is almost the same as  in (i), with the following modifications.
First, \reff{psine} is modified by changing the last term:
\beaa
 \psi^{n,\e}(\ul\l):=U_1(t,\ul\xi)-U_2(t,\ul\zeta)- n\dbE\big[\Upsilon_t(\xi-\zeta)\big]- \textcolor{black}{n^{q}}\mathbb{E}\big[|\xi_{t}-\zeta_{t}|^{\textcolor{black}{p}}\big]-\e\Big(1- {t\over 2T}\Big)\mathbb{E}\big[ \Upsilon_t(\xi)\big].
 \eeaa
 This is used to derive the counterpart of \reff{U1difference}, based on $\Psi(\ul{\hat\l}) \ge \Psi(\hat t, \ul{\hat\xi}, \ul{\hat\xi})$:
 \beaa
 \left.\ba{c}
\dis 0\le
 U_2(\hat t, \ul{\hat\xi})- U_2(\hat t, \ul{\hat\zeta})
 - n\dbE\big[\Upsilon_{\hat t}(\ul{\hat\xi}- \ul{\hat\zeta})\big]- \textcolor{black}{n^{q}}\mathbb{E}\big[|\hat\xi_{\hat t}-\hat\zeta_{\hat t}|^{\textcolor{black}{p}}\big]
 \ms\\
 \dis 
 + \sum_{i\ge 0}2^{-i}
 \dbE\Big[\Upsilon_{\hat t}(\hat \xi - (\zeta_i)_{\cdot\wedge t_i})
               -\Upsilon_{\hat t}(\hat \zeta- (\zeta_i)_{\cdot\wedge t_i})
        \Big].
        \ea\right.
 \eeaa
 Then the estimates {\reff{Vreg1}} for $U_2$ will lead to \reff{xietane}.

 Next, we modify \reff{fne} by moving the $\e$-term from $ \f^{\hat\eta}_2(\ul\eta,\ul{\hat\arm})$ to $ \f^{\hat\th}_1(\ul\th,\ul{\hat\arm})$:
\beaa
\left.\ba{c}
\dis \f^{\hat\th}_1(\ul\th,\ul{\hat\arm})
:=
\bar\Upsilon_0\big(\ul\th, \hat\th)\big)
+\textcolor{black}{2^{p-1}} n\Upsilon_0(\ul\th,(\hat t, \ul{\hat \arm}))
+\e\Big(1- {t\over 2T}\Big)\dbE[\Upsilon_t(\xi)]+\pi^{n,\e}_1(\ul\th);
\ms\\
\dis  \f^{\hat\eta}_2(\ul\eta,\ul{\hat\arm})
:=
\bar\Upsilon_0\big(\ul\eta, \hat\eta)\big)
+\textcolor{black}{2^{p-1}} n\Upsilon_0(\ul\eta,(\hat t, \ul{\hat \arm}))
+\pi^{n,\e}_2(\ul\eta).
\ea\right.
\eeaa
This ensures that the calculation in \reff{psiNPsi} remains true. Moreover, in this case the derivatives of $\f^{\hat\eta}_2$ is simplified slightly. Then \reff{comparison1} will become: for all  $\e\le \e_{m}$, and $n\ge n_{\e, m}$,
\beaa
c_1
\le C\e
+C_\e \textcolor{black}{n^{-\tilde{q}}} + \left\{\ba{lll} \qq\qq 0,\qq\qq\q  \mbox{under \reff{s<t}};\\ \e (C -{1\over 2T})
    \big\|\hat\xi^{n,\e}_{\cdot\wedge\hat t^{n,\e}}\big\|_{\textcolor{black}{p}}^{\textcolor{black}{p}},\q  \mbox{under \reff{case2}}.
    \ea\right.
\eeaa
This leads to the desired contradiction as in Case 1.
\qed


\section{Appendix}
\label{sect-appendix}

In this Appendix we prove some technical results. { We first prove the  DPP \reff{DPP}, following the arguments in \cite{cosso3}.}

\ms
\noindent {\no{\bf Proof of the DPP \reff{DPP}.} }   {Denote the right side of \reff{DPP}:
$$
   \Pi_t(\underline{\xi}):=\inf_{\a\in \cA_{[t, T]}}
\Big\{V_{t+\d}(\ul X^{t,\xi, \a})
         \!+\! \int^{t+\delta}_{t}\!\!\!\!f _s( \ul X^{t,\xi,\a},\ul\a_s)ds\Big\}.
$$
We first prove that $V_t(\underline{\xi})\geq \Pi_t(\underline{\xi})$. For every fixed $\alpha\in \cA_{[t, T]}$, by the definition of $V$,
 $$
   V_{t+\delta}(\underline{X}^{t,\xi,\alpha})=\inf_{\a'\in \cA_{[t+\delta, T]}}\left[\int^T_{t+\delta}f_s(\underline{X}^{t+\delta,\underline{X}^{t,\xi,\alpha},\a'},{\underline{\a}}'_s)ds+g(\underline{X}^{t+\delta,\underline{X}^{t,\xi,\alpha},\a'})\right].
 $$
By the uniqueness of equation \reff{state1}, we have the flow property
$
   X^{t,\xi,\a}=X^{t+\delta,X^{t,\xi,\a}, \a}.
$
Hence, by choosing $\a'=\a$,
 $$
   V_{t+\delta}(\underline{X}^{t,\xi,\alpha})+\int^{t+\delta}_tf_s(\underline{X}^{t,\xi,\alpha},{\underline{\a}}_s)ds
   \leq\int^T_{t}f_s(\underline{X}^{t,\xi,\alpha},{\underline{\a}}_s)ds+g(\underline{X}^{t,\xi,\alpha})=J_t(\underline{\xi},\underline{\a}).
 $$
Taking the infimum in $\cA_{[t, T]}$, we obtain that $V_t(\underline{\xi})\geq \Pi_t(\underline{\xi})$.}

{ Next, we prove $V_t(\underline{\xi})\leq \Pi_t(\underline{\xi})$. For every $\varepsilon>0$, let $\a^{\e}\in \cA_{[t, T]}$ be such that
\begin{eqnarray}\label{0704add}
     \Pi_t(\underline{\xi})\geq V_{t+\d}(\ul X^{t,\xi, \a^{\e}})
         \!+\! \int^{t+\delta}_{t}\!\!\!\!f _s( \ul X^{t,\xi,\a^{\e}},\ul\a^{\e}_s)ds-\e.
\end{eqnarray}
By the definition of $V_{t+\d}(\underline{X}^{t,\xi,\a^{\e}})$,  there exists ${\a'}^{\e}\in \cA_{[t+\d, T]}$ such that
\begin{eqnarray}\label{0704add1}
            V_{t+\d}(\underline{X}^{t,\xi,\a^{\e}})\geq\int^T_{t+\delta}f_s(\underline{X}^{t+\d,{X}^{t,\xi,\a^{\e}},{\a'}^{\e}},{\underline{\a'}^{\e}}_s)ds
            +g(\underline{X}^{t+\d,{X}^{t,\xi,\a^{\e}},{\a'}^{\e}})-\e.
\end{eqnarray}
Notice that $
           {\a''}^{\e}:=\a^{\e}1_{[t,t+\delta]}+{\a'}^{\e}1_{(t+\delta,T]}\in \cA_{[t,T]}$, and that $
                {X}^{t+\d,{X}^{t,\xi,\a^{\e}},{\a'}^{\e}}=X^{t,\xi,{\a''}^{\e}}
$, again by the uniqueness of equation  \reff{state1}. Hence, by \reff{0704add} and \reff{0704add1}, we get
\begin{eqnarray*}
     \Pi_t(\underline{\xi})\geq  \int^{T}_{t}\!\!\!\!f _s( \ul X^{t,\xi,{\a''}^{\e}},\ul{\a''}^{\e}_s)ds+g(\underline{X}^{t,\xi,{\a''}^{\e}})-2\e\geq V_t(\xi)-2\e.
\end{eqnarray*}
Letting $\e\rightarrow0$, we have $V_t(\underline{\xi})\leq \Pi_t(\underline{\xi})$.} \qed

\bs
 {We next prove Lemma \ref{lem-dupire} concerning the uniqueness of the path derivatives  for  $\f\in C^{1,2}(\ol\cX^p_{[\hat t, T]})$.}

\no{\bf Proof of Lemma \ref{lem-dupire}.}   We first prove the Lemma in the case that $t\in [\hat t, T]$, $\xi\in  \cM^p_{[t, T]}$.

First, let $X\in \cM^p_{[\hat t, T]}$ be such that ${X_{\cd\wedge t}} = \xi_{\cd\wedge t}$ and $\b_s ={\bf 0}, \g_s={\bf 0}$, $s\ge t$, in \reff{cMp}. Then ${X_{\cd\wedge s}} =  {X_{\cd\wedge t}} = \xi_{\cd\wedge t}$ for all $s\ge t$, and thus by \reff{Ito} we obtain,
\beaa
\f_{t+\d}(   \ul \xi_{\cd\wedge t}) = \f_t( \ul\xi_{\cd\wedge t}) + \int_t^{t+\d} \pa_t \f_s(  \ul \xi_{\cd\wedge t}) ds,\q  0<\d\le T-t.
\eeaa
Since $\pa_t \f$ is continuous, we obtain a  representation of $\pa_t \f$ which  implies its uniqueness:
\beaa
\pa_t \f_t(  \ul\xi_{\cd\wedge t}) = \lim_{\d\to 0} {\f_{t+\d}(  \ul\xi_{\cd\wedge t}) - \f_t( \ul\xi_{\cd\wedge t}) \over \d}.
\eeaa

Next, replace the above $X$ with $\b_s \equiv \b_t$, $s\in [t, T]$ for some arbitrary $\b_t\in \dbL^p(\cF_t; \dbR^d)$. By \reff{Ito} we have, for any $0<\d\le T-t$,
\beaa
\f_{t+\d}(   \ul X) = \f_t( \ul X) + \int_t^{t+\d} \pa_t \f_s( \ul X) ds + \dbE\Big[\int_t^{t+\d} \partial_X\f_s(\ul X) ds \cd \b_t\Big].
\eeaa
By the uniqueness of $\pa_t \f$ and the continuity of $\pa_t\f, \pa_X \f$, we see that
\beaa
\dbE\big[\partial_X\f_t(\ul\xi) \cd \b_t\big] = \lim_{\d\to 0}\Big[{\f_{t+\d}(  \ol X) - \f_t( \ol X)\over \d} - \pa_t \f_t( \ul\xi)\Big]
\eeaa
is unique. Since $\b_t$ is arbitrary, we see that $\partial_X\f_t(\ul\xi)$ is unique, $\dbP$-a.s.
Similarly, by considering $X$ with $\b_s \equiv 0$ but $\g_s \equiv \g_t$, $s\in [t, T]$, we obtain the uniqueness of $\dbE\big[\partial_{\bx X}\f_t(\ul\xi) : \g_t\g_t^\top\big]$. Since $\pa_{\bx X} \f$ is symmetric and $\g_t$ is arbitrary, we see that $\partial_{\bx X}\f_t(\ul\xi)$ is unique, a.s.

Finally, for arbitrary $(t, \xi)\in \ol\cX^p_{[\hat t, T]}$, clearly there exist a sequence $\xi^n\in  \cM^p_{[\hat t, T]}$ such that $\dis\lim_{n\to\infty}\|\xi^n_{\cd\wedge t} - \xi_{\cd\wedge t}\|_p =0$. By the above arguments $\pa_t \f, \pa_X\f, \pa_{xX}\f$ are unique at $(t, \xi^n)$. Then it follows from their continuity that they are unique at $(t, \xi)$ as well.
\qed

\bs
 {The next result concerns the absolute continuity of the singular part $\phi$ of test functions.}

\no{\bf Proof of Proposition \ref{prop-C+cXp}.}
For any $\a\in \cA_{[\tilde t, T]}$, by the standard It\^{o} formula we have:
\bea
\label{DcI}
&&\dis\Big|\dbE\big[\big|\cI^\a_{t+\d}(X)\big|^p - \big|\cI^\a_t( X) \big|^p\big]\Big| \nonumber\\
&&\dis= \Big| \dbE\Big[ \int_{t}^{t+\d} \big[p|\cI^\a_s(X)|^{p-2} (\cI^\a_s(X)) \cd (\b_s - \tilde{b}^\a_s)+  {p\over 2} |\cI^\a_s(X)|^{p-2}|\g_s - \tilde{\si}^\a_s|^2\nonumber\\
 &&\dis\qq +  {p(p-2)\over 2} |\cI^\a_s(X)|^{p-4}|\cI^\a_s(X)(\g_s - \tilde{\si}^\a_s)|^2\big]ds\Big]\Big|\nonumber\\
 &&\dis \le C_p \dbE\Big[ \int_{t}^{t+\d} \big[|\cI^\a_s(X)|^{p-1}|\b_s - \tilde{b}^\a_s|+  |\cI^\a_s(X)|^{p-2} |\g_s - \tilde{\si}^\a_s|^2\big]ds\Big].
 \eea
 Then, recalling \reff{tildebsifbound} with $\G_*$ and \reff{C+cXp2}, and noting that $\k(\a)$ does not depend on $t$,
 \beaa
&&\dis \Big|\phi_{t+\d}(\ul X) - \phi_t(\ul X)\Big| \le k \sup_{\a\in \cA_{[\tilde t, T]}} \Big|\dbE\big[\big|\cI^\a_{t+\d}(X)\big|^p - \big|\cI^\a_t(X)\big|^p\big]\Big|  +  \sup_{\a\in \cA_{[\tilde t, T]}}\Big[\int_t^{t+\d}|\tilde{f}^\a_s|ds\Big]\\
\dis &&\q \le C_{p,k}\sup_{\a\in \cA_{[\tilde t, T]}}\dbE\Big[ \int_{t}^{t+\d} \big[|\cI^\a_s(X)|^p+  |\b_s|^p  +  |\g_s|^p  +|\tilde b^\a_s|^p+|\tilde \si^\a_s|^p +1\big]ds\Big]\\
\dis &&\q \le C_{p,k}\int_{t}^{t+\d} \dbE\Big[ |X_s|^p + |\tilde \xi_{\tilde t}|^p + |\b_s|^p + |\g_s|^p + |\G_*|^p+ 1 \Big]ds.
\eeaa
This implies that $\phi_t( \ul X)$ is absolutely continuous in $t$.

We next prove (i). Note that $\cA_{[\tilde t, T]} = \{\a\oplus_t \tilde \a: \a\in \cA_{[\tilde t, T]}, \tilde \a\in \cA_{[t, T]}\}$, where the concatenation is defined as: $\a\oplus_t\tilde \a := \a \1_{[\tilde t, t)} + \tilde \a\1_{[t, T]}$.  Recall the $\k(\a)$ in \reff{C+cXp2} and note that, since $t\ge t'$, we have $\k(\a\oplus_t \tilde \a) = \k(\a)$.  Thus
\beaa
\dis \phi_{t+\d}( \ul X) - \phi_t( \ul X) &=& \inf_{\a\in \cA_{[\tilde  t, T]}} \inf_{\tilde \a\in \cA_{[t, T]}}\Big[ k \dbE\big[ \big|\cI^{\a\oplus_t \tilde \a}_{t+\d}(X)\big|^p\big] + F_{t+\d}^{\a\oplus_t \tilde \a} + \k(\a)\Big] \\
&&\dis  -  \inf_{\a\in \cA_{[\tilde  t, T]}} \Big[ k  \dbE\big[\big|\cI^{\a}_t( X)\big|^p \big]+ F^\a_t+ \k(\a)\Big]\\
&\le&  \inf_{\tilde \a\in \cA_{[t, T]}} \sup_{\a\in \cA_{[\tilde  t, T]}} \Big\{k\dbE\Big[ \big|\cI^{\a\oplus_t \tilde \a}_{t+\d}(X)\big|^p  -  \big|\cI^{\a}_t(X) \big|^p\Big] + \int_t^{t+\d} \tilde{f}^{\tilde \a}_s ds\Big\}.
\eeaa
We remark that the above involves $\a$ only on $[\tilde t, t)$. By \reff{DcI} we prove  \reff{pa+phi-est2} immediately.

It remains to prove (ii). For any $\d>0$, by \reff{C+cXp2} there exists $\a^\d$ s.t.
\beaa
\phi_t( \ul X) \ge  k \dbE\big[\big|\cI^{\a^\d}_t( X)\big|^p\big] +F^{\a^\d}_t + \k(\a^\d) - \d^2.
\eeaa
It is clear that the above property, and hence $\a^\d$, does not involve $(\b_s, \g_s)_{s\ge t}$. Then,
\beaa
&&\dis \phi_{t+\d}( \ul X) - \phi_t( \ul X) \\
&&\dis \le \Big[ k \dbE\big[\big|\cI^{\a^\d}_{t+\d}(X)\big|^p\big] +F^{\a^\d}_{t+\d} + \k(\a^\d)\Big] - \Big[ k \dbE\big[\big|\cI^{\a^\d}_t(X)\big|^p\big] +F^{\a^\d}_t + \k(\a^\d)\Big] + \d^2\\
&&\dis = k\dbE\Big[  \big|\cI^{\a^\d}_{t+\d}(X)\big|^p  -  \big|\cI^{\a^\d}_t(X)\big|^p \Big] + \int_t^{t+\d} \tilde{f}^{\a^\d}_s ds +\d^2.
\eeaa
By  \reff{DcI} again we can easily prove  \reff{pa+phi-est1}.
\qed

 In order to prove Proposition \ref{prop-consistent},  {establishing the consistency between viscosity solutions and classical solutions,} we need the following slight generalization of \cite[Lemma F.2]{cosso3}.
  \begin{lem}
\label{lem-xax0630}
Let $t\in [0, T]$ and  $F: \cA_{T} \to \dbR$ be continuous.  Then,
\bea
\label{xzx}
\lim_{s\downarrow t}\inf_{\a_s\in \cA_{s}}F(\ul\a_s) = \inf_{\a_t\in \cA_{t}}F(\ul\a_t).
\eea
\end{lem}
\proof Since $ \cA_{t}\subset  \cA_{s}$ for $s\geq t$, we immediately have
 $
\dis\lim_{s\downarrow t}\inf_{\a_s\in \cA_{s}}F(\ul\a_s) \le \inf_{\a_t\in \cA_{t}}F(\ul\a_t).
$
Assume by contradiction that \reff{xzx} does not hold. Then there exist $\varepsilon>0$ and $N\in\dbN$ such that $\cB_n \neq \emptyset$ for all $n\geq N$, where, denoting $t_n:= t+{1\over n}$,
$$
       \cB_n:=\Big\{\a_{t_n}\in \cA_{t_n}: F(\ul\a_{t_n})\leq \inf_{\a_t\in \cA_{t}}F(\ul\a_t)-\e\Big\}.
$$
Note that $\cA_{t_{n+1}}\subset\cA_{t_n}$ and thus $\cB_{n+1}\subset\cB_n$. Moreover, by the continuity of $F$ we see that $\cB_n$ is closed. Since $\cA_T$ is complete, by Cantor's intersection theorem, there exists
$\hat{\a}\in \cA_T$ such that $\hat{\a} \in \bigcap_{n\geq N}\cB_n\subset \bigcap_{n\geq N}\cA_{t_n}$.  Note that $\bigcap_{n\geq N}\cA_{t_n}$ is the set of $\cF_{t+}$-measurable $A$-valued random variables. Since the augmented filtration of $\dbF$ is right continuous, there exists $\hat\a' \in \cA_t$ such that $\hat\a' = \hat\a$, $\dbP$-a.s., see e.g. \cite[Proposition 1.2.1]{zhang}. Then  $
F(\ul{\hat\a'}) = F(\ul{\hat\a})\leq \inf_{\a_t\in \cA_{t}}F(\ul\a_t)-\e.
$
This is a desired contradiction.
\qed

\ms
\no{\bf Proof of Proposition \ref{prop-consistent}.} We shall only prove the equivalence of subsolution properties. The supersolution case follows essentially the same argument.

 First, if $U$ is a viscosity subsolution, {then $U$ is a viscosity $p$-subsolution for some $p\geq 2$}. Note that $C^{1,2}(\ol\cX^2_{[0, T]}) \subset C^{1,2}(\textcolor{black}{\ol\cX^p_{[0, T]}})$. For any $(t,  \xi)\in \textcolor{black}{\ol\cX^p_{[0, T)}}$, it is clear that $(U, \mathbf{0})\in \textcolor{black}{\mathfrak{F}^+_{p}}U( t, \ul\xi)$, then it is straightforward to derive from \reff{sub} the classical subsolution property of $U$ at $(t, \xi)$. Moreover, by the arguments in Remark \ref{rem-viscosity2} (ii), we obtain the classical subsolution property at all $(t,  \xi)\in \ol\cX^2_{[0, T)}$.

On the other hand, assume $U$ is a classical subsolution. {For every $p\geq 2$ (namely with $p_0=2$ in Definition \ref{defn-viscosity} (iii)),} fix  $(t, \xi) \in \textcolor{black}{\ol\cX^p_{[0, T)}}$ and $(\varphi, \phi)\in \textcolor{black}{\mathfrak{F}^+_{p}}U(t, \ul\xi)$. For any $\a\in \cA_{[t, T]}$,
 recall \reff{barXtxi} and denote $\bar X^\a:= \bar X^{t,\xi,\a}$. By the It\^o formula \reff{Ito} and \reff{cAtest+} we have: for $\forall\d>0$,
 \beaa
 0 &\ge&\Big[\big[U-(\varphi+\phi)\big]_{t+\d}( \ul{\bar X}^{\a}) -  \big[U-(\varphi+\phi)\big]_t(\ul\xi)\Big]\\
 &=&  \int_t^{t+\d} \dbE\Big[\pa_t (U-\f)_s(  \ul{\bar X}^{\a}) + \pa_X (U-\f)_s(  \ul{\bar X}^{\a}) \cd b^{t, \xi, \a}_{s} \\
 && \hspace{15mm}
 + {1\over 2} \pa_{\bx X}(U-\f)_s( \ul{\bar X}^{\a}) :  (\si\si^\top)^{t, \xi, \a}_{s} - \dot\phi_s(  \ul{\bar X}^{\a})\Big]ds.
 \eeaa
 By the smoothness of $U-\f$ and the regularity of  $b, \si, f$, including Assumption \ref{assum-regt}, we have
  \beaa
 && \int_t^{t+\d}\Big[\pa_t U_t( \ul\xi) +H_t( \ul\xi,\partial_X U_t(\ul\xi),\partial_{\bx X}U_t(\ul\xi), \ul\a_s)\Big]ds\\
 &\le&  \int_t^{t+\d} \Big[\pa_t \f_t( \ul\xi) + H_s( \ul\xi_{\cd\wedge t},\partial_X\f_t(\ul\xi),\partial_{\bx X}\f_t(\ul\xi), \ul\a_s) + \dot\phi_s(  \ul{\bar X}^{\a})\Big]ds + o(\d),
 \eeaa
 where $o(\d)$ is uniform in $\a$. Then
 \beaa
 &&\inf_{\a\in \cA_{[t, T]}}{1\over \d}\int_t^{t+\d}\Big[\pa_t U_t( \ul\xi) +H_t( \ul\xi,\partial_X U_t(\ul\xi),\partial_{\bx X}U_t(\ul\xi), \ul\a_s)\Big]ds\\
 &\le& \inf_{\a\in \cA_{[t, T]}}{1\over \d} \int_t^{t+\d} \Big[\pa_t \f_t( \ul\xi) + H_s( \ul\xi_{\cd\wedge t},\partial_X\f_t(\ul\xi),\partial_{\bx X}\f_t(\ul\xi), \ul\a_s) + \dot\phi_s(  \ul{\bar X}^{\a})\Big]ds + o(1),
 \eeaa
 Note that,  by Assumption \ref{assum-regt}, we have $H_t( \ul\xi,\partial_X U_t(\ul\xi),\partial_{\bx X}U_t(\ul\xi), \cdot)$ is continuous.
 Send $\d\to 0$,  by  Lemma \ref{lem-xax0630}   we obtain
  \beaa
 \cL U_t( \ul\xi) \le \partial_t \f_t(\ul\xi)+\liminf_{\d\to 0}\inf_{\a \in \cA_{[t, T]}} {1\over \d}\int_t^{t+\d}\!\!\! \big[ H_s( \ul\xi_{\cd\wedge t},\partial_X\f_t(\ul\xi),\partial_{\bx X}\f_t(\ul\xi), \ul\a_s) + \dot\phi_s(  \ul{\bar X}^{\a})\big]ds.
 \eeaa
Then \reff{sub} follows from  the classical subsolution property of $U$.
 \qed

 \bs
 We next prove Proposition \ref{prop-MFC2},  {which establishes the basic properties of the value function in the mean control setting.}

 \no{\bf Proof of Proposition \ref{prop-MFC2}.} (i) We shall prove only $\check V_t(\ul\xi) \ge \check V_t(\ul \xi')$.  It suffices to prove $\check J_t(\ul\xi, \ul\a) \ge \check V_t(\ul \xi')$ for an arbitrarily fixed $\a\in \cA_{[t, T]}$. By standard approximation arguments, we may assume without loss of generality that, for some $t=t_0<\cds<t_n = T$:
\bea
\label{discretea}
\a_s = \sum_{i=0}^{n-1} \a_i(\o, B^1_{\cd\wedge t_i}, B^0_{\cd\wedge t_i}) \1_{[t_i, t_{i+1})}(s),
\eea
 where $\a_i: \O \times \dbX\to A$ is $\cF_0$ measurable in $\o$ and uniformly continuous in $\bx$. Moreover,   since $\cF_0$ is rich enough, we can generate $\cF_0$-measurable random variables $\{{\bf U_i}\}_{i\ge 0}$ which are i.i.d. with distribution Uniform$([0,1])$ and are independent of $(\xi, \xi', \a)$, and hence are independent of $X^{t,\xi, \a}$. Denote $X^{\tilde \a} := X^{t, \xi, \tilde \a}$, $X^{'\tilde \a}:= X^{t,\xi', \tilde\a}$ for arbitrary $\tilde \a$.

First, since $\dbP_{\xi_{\cd\wedge t} |\cF^0_t} = \dbP_{\xi'_{\cd\wedge t} |\cF^0_t}$, $\dbP$-a.s., then $\dbP_{(\xi_{\cd\wedge t_0}, B^0)} = \dbP_{(\xi'_{\cd\wedge t_0}, B^0)}$. By using ${\bf U}_0$, one can easily construct $\a'_{t_0} \in \cA_{t_0}$ such that $\dbP_{(\a_{t_0}, \xi_{\cd\wedge t_0}, B^0)} = \dbP_{(\a'_{t_0}, \xi'_{\cd\wedge t_0}, B^0)}$. This clearly implies that $\dbP_{(\a_{t_0}, X^\a_{\cd\wedge t_1},  B^0)} = \dbP_{(\a'_{t_0}, X^{'\a'}_{\cd\wedge t_1}, B^0)}$. Repeat the arguments by using the ${\bf U}_i$, we may construct $\a'=\sum_{i=0}^{n-1} \a'_{t_i}\1_{[t_i, t_{i+1})} \in \cA_{[t, T]}$ such that  $\dbP_{(\a, X^\a,  B^0)} = \dbP_{(\a', X^{'\a'}, B^0)}$. Then clearly $\check J_t(\ul\xi, \ul\a)  = \check J_t(\ul\xi', \ul\a') \ge \check V_t(\ul \xi')$. By the arbitrariness of $\a$ we have $\check V_t(\ul\xi) \ge \check V_t(\ul \xi')$. Similarly $\check V_t(\ul\xi') \ge \check V_t(\ul \xi)$, then equality holds and hence \reff{checkV} is well defined.

Moreover, by \reff{Vreg} and Proposition \ref{prop-MFC1} we obtain part of \reff{Vreg55} immediately:
\bea
\label{checkVreg2}
  |\check V_t( \mu)| \le C\Big(1+ W^{\beta}_2(\mu_{\cd\wedge t}, \d_{\bf 0})\Big), \q |\check V_t(\mu)- \check V_{t}(\mu')|\leq C W^{\beta}_2(\mu_{\cd\wedge t}, \mu'_{\cd\wedge t}).
\eea

 (ii) We proceed in several steps. For $(t,\xi)\in \ol\cX^2_{[0,T]}$ and $\a\in \cA_{[t, T]}$,  denote
 \beaa
 X^\a := X^{t,\xi,\a},\qq Y^\a_t:= \dbE_{\cF^0_t}\Big[\check g\big(X^\a, \dbP_{ X^\a|\cF^0_T}\big)
                +
                \int_t^T \check f_s\big(X^\a,\a_s, \dbP_{ (X^\a, \a_s)|\cF^0_s}\big)ds
        \Big].
        \eeaa

 {\it Step 1.} We first assume $\xi\in \cX_2^\perp(\dbF^0_t)$. Then $\dbP_{\xi_{\cd\wedge t}|\cF^0_t} = \dbP_{\xi_{\cd\wedge t}}$, and thus  the claimed result follows directly from \reff{checkV}.

 Moreover, assume $\a\in \cA_{[t, T]}$ takes the form \reff{discretea}. For each $\bx^0 \in C^0([0, T]; \dbR^{d_0})$ with $\bx^0_0=0$, where $d_0$ is the dimension of $B^0$, denote
 \beaa
 &\dis \a^{\bx^0} := \sum_{i=0}^{n-1} \a_i(\o, B^1_{\cd\wedge t_i}, (\bx^0\otimes_t B^0)_{\cd\wedge t_i}) \1_{[t_i, t_{i+1})}\in \cA_{[t, T]},~\mbox{and}~ \check J'_t(\ul\xi, \ul\a, \bx^0) := \check V_t(\ul\xi, \ul \a^{\bx^0}),\\
 &\dis \mbox{where}\q (\bx^0\otimes_t B^0)_s := \bx^0_s \1_{[0, t)}(s) + [\bx^0_t + B^0_s-B^0_t]\1_{[t, T]}(s).
 \eeaa
We note that this is in the spirit of the regular conditional probability distribution, but we avoid it by considering the special form \reff{discretea}. Since $\xi\in \cX_2^\perp(\dbF^0_{t})$, we can easily show that
 \bea
 \label{YaJ'}
Y^\a_t =  J'_t(\ul\xi, \ul\a, B^0), \q\dbP\mbox{-a.s.}
        \eea
 In particular, this implies that $Y^\a_t \ge \check V_t(\ul\xi)$, $\dbP$-a.s.

{\it Step 2.} We next consider the case $\xi = \sum_{i=1}^n \xi_i \1_{E_i}$, where $\xi_i\in \cX_2^\perp(\dbF^0_t)$ and $\{E_i\}_{1\le i\le n} \subset \cF^0_t$ form a partition of $\O$. Note that $\dbP_{\xi|\cF^0_t} = \sum_{i=1}^n \dbP_{\xi_i} \1_{E_i}$, $\dbP$-a.s. Then
\beaa
V'_t(\ul\xi) := \dbE\big[\check V_t(\dbP_{\xi|\cF^0_t})\big] = \sum_{i=1}^n \check V_t(\dbP_{\xi_i}) \dbP(E_i).
\eeaa
We first prove $V_t(\ul\xi) \ge V'_t(\ul\xi)$. For any $\a\in \cA_{[t, T]}$ taking the form \reff{discretea}, by \reff{YaJ'} we have
\beaa
J_t(\ul\xi,\ul\a) = \dbE[Y^\a_t]  = \dbE\Big[\sum_{i=1}^n Y^{t, \xi_i, \a}_t\1_{E_i}\Big] \ge \dbE\Big[\sum_{i=1}^n \check V_t(\dbP_{\xi_i})\1_{E_i}\Big] = V'_t(\ul\xi).
\eeaa
Now by standard approximation arguments, we have $J_t(\ul\xi,\ul\a) \ge V'_t(\ul\xi)$ for all $\a\in \cA_{[t, T]}$, and thus $V_t(\ul\xi) \ge V'_t(\ul\xi)$.

To see the opposite inequality, fix an arbitrary $\e>0$. For each $i=1,\cds, n$, there exists $\a^i$ taking the form \reff{discretea} (with different $n$ there) such that $Y^{t, \xi_i, \a^i}_t \le \check V_t(\dbP_{\xi_i}) + \e$. Construct
$
\a^\e := \sum_{i=1}^n \a^i \1_{E_i}.
$
One can easily see that $Y^{t,\xi, \a^\e}_t = \sum_{i=1}^n Y^{t, \xi_i, \a^i}_t \1_{E_i}$. Then
\beaa
V_t(\ul\xi) &\le& \check J_t(\ul\xi, \ul \a^\e) = \dbE\big[Y^{t,\xi, \a^\e}_t\big] = \dbE\Big[\sum_{i=1}^n Y^{t, \xi_i, \a^i}_t \1_{E_i}\Big] \\
&\le& \dbE\Big[\sum_{i=1}^n \big[\check V_t(\dbP_{\xi_i}) + \e\big] \1_{E_i}\Big]  =  V'_t(\ul\xi) + \e.
\eeaa
Since $\e>0$ is arbitrary, we obtain $V_t(\ul\xi) \ge V'_t(\ul\xi)$, and hence the equality.

{\it Step 3.} By the regularity of $V$ in \reff{Vreg} and that of $\check V$ in \reff{checkVreg2}, and by standard approximation arguments, similarly to \reff{discretea} we may assume without loss of generality that $\xi$ takes the following form: for some $0=t_0<\cds<t_n=t$,
\beaa
&\dis\xi_s \equiv \xi'_{t_0}, ~ t\in [t_0, t_1];\q \xi_s = {t_{i+1}-t\over t_{i+1}-t_i} \xi'_{t_{i-1}} + {t-t_i\over t_{i+1}-t_i} \xi'_{t_i} ,~ t\in [t_i, t_{i+1}], ~i=1,\cds, n-1;\\
&\dis\mbox{where}\q \xi'_{t_i} = h_i\Big(\o, (B^1_{t_1}, \cds, B^1_{ t_i}), (B^0_{t_1}, \cds, B^0_{ t_i})\Big),
\eeaa
and  $h_i: \O \times \dbR^{di}\to \dbR^d$ is $\cF_0$ measurable in $\o$ and uniformly continuous in $y\in \dbR^{di}$. Fix $\e>0$, and let $\{O^{\e}_j\}_{1\le j\le n_{\e}}$ be a partition of $\big\{y\in \dbR^{d_0n}: |y|\le {1\over \e}\big\}$ such that each $O^\e_j$ has diameter less than $\e$. Fix an arbitrary $y^\e_j\in O^{\e}_j$ for each $j$, and denote $O^{\e,i}_j:= \{Proj_i(y): y\in O^\e_j\}$, where $Proj_i(y)\in \dbR^{di}$ denotes the first $i$ components of $y\in \dbR^{dn}$.  Define
\beaa
&\dis\xi^\e_s \equiv \xi'_{t_0}, ~ t\in [t_0, t_1];\q \xi^\e_s = {t_{i+1}-t\over t_{i+1}-t_i} \xi^{'\e}_{t_{i-1}} + {t-t_i\over t_{i+1}-t_i} \xi^{'\e}_{t_i} ,~ t\in [t_i, t_{i+1}], i=1,\cds, n-1;\\
&\dis\mbox{where}~ \xi^{'\e}_{t_i} = \sum_{j=1}^{n_\e}\xi^{'\e,i}_{t_i} \1_{\{(B^0_{t_1}, \cds, B^0_{t_i}) \in O^{\e,i}_j\}},~\xi^{'\e,i}_{t_i}:= h_i\Big(\o, (B^1_{t_1\wedge t_i}, \cds, B^1_{t_n\wedge t_i}), Proj_i(y^\e_j)\Big).
\eeaa
Note that $\xi^{'\e,i}_{t_i}$ is independent of $\cF^0_t$, and, on $\{(B^0_{t_1}, \cds, B^0_{t_i}) \notin \cup_j O^{\e,i}_j\}$,  $\xi^{'\e}_{t_i} =0$, which is also independent of $\cF^0_t$. Then one can easily see that $\xi^\e$ takes the form as in Step 2, and thus it follows from Step 2 that $V_t(\ul\xi^\e) = V'_t(\ul\xi^\e)$. Send $\e\to 0$, by the desired regularities we obtain $V_t(\ul\xi) =V'_t(\ul\xi)$.

(iii) follows from (ii) and DPP \reff{DPP} for $V$. Finally, for $0\le t< t+\d\le T$, $\mu\in \cP_2(\dbX)$, and $\xi \in  \cX_2^\perp(\dbF^0; t, \mu)$, by (iii) we have, again denoting $X^\a := X^{t,\xi,\a}$,
\beaa
\check V_t(\mu) - \check V_{t+\d}(\mu_{\cd\wedge t})
&=& \inf_{\a\in \cA_{[t, T]}} \dbE\Big[\check V_{t+\d}( \dbP_{X^{\a}|\cF^0_{t+\d}}) - \check V_{t+\d}(\dbP_{X^{\a}_{\cd\wedge t}})
\\
&&\qq
+ \int^{t+\delta}_{t}\check f _s(X^{\a}, \a_s, \mathbb{P}_{(X^{\a},\a_s)|\cF^0_s})ds\Big].
\eeaa
Apply \reff{checkVreg2} on $\check V_{t+\d}( \dbP_{X^{\a}|\cF^0_{t+\d}}) - \check V_{t+\d}(\dbP_{X^{\a}_{\cd\wedge t}})$, we can easily obtain
\beaa
\big|\check V_t(\mu) - \check V_{t+\d}(\mu_{\cd\wedge t}) \big| \le C\textcolor{black}{\big(1+ W^{\beta}_2(\mu_{\cd\wedge t}, \d_{\bf 0})\big) {\d}^{\beta\over 2}}.
\eeaa
This, together with \reff{checkVreg2}  again, implies \reff{Vreg55} immediately.
 \qed

 \ms

 {We now turn to the proof of Proposition \ref{prop-MFCIto} on the functional It\^{o} formula in the mean field control setting.} For that purpose we first introduce briefly the path derivatives and specify the space $C_b^{1,2,2}([0, T]\times \dbX \times \cP_2(\dbX))$, which actually serves as the technical conditions for the proposition. We refer to \cite{CD2, WZ} for more details. We first need to extend $U$ to the c\`adl\'ag space. Let $\dbD$ denote the space of $d$-dimensional c\`adl\'ag paths on $[0, T]$, equipped with the uniform norm, and $\cP_2(\dbD)$ the space of square integrable probability measures on $\dbD$, equipped with the $2$-Wasserstein distance. Given an adapted function $U: [0, T]\times \dbX \times \cP_2(\dbX) \to \dbR$, let $\hat U: [0, T]\times \dbD \times \cP_2(\dbD) \to \dbR$ be an extension. The path derivatives of $\hat U$ are functions $\pa_t \hat U: [0, T]\times \dbD \times \cP_2(\dbD) \to \dbR$, $\pa_\bx \hat U: [0, T]\times \dbD \times \cP_2(\dbD) \to \dbR^d$ and $\pa_\mu \hat U: [0, T]\times \dbD \times \cP_2(\dbD) \times \dbR^d \to \dbR^d$ determined by:
\beaa
&\dis \pa_t \hat U_t(\bx, \mu) := \lim_{\d\to 0} {1\over \d}\Big[\hat U_{t+\d}(\bx_{\cd\wedge t}, \mu_{\cd\wedge t}) - \hat U_t(\bx_{\cd\wedge t}, \mu_{\cd\wedge t})\Big];\\
&\dis \hat U_t(\bx + \D x \1_{[t, T]}, \mu) - \hat U_t(\bx, \mu) = \pa_\bx \hat U_t(\bx, \mu)\cd \D x + o(|\D x|),\q\forall \D x\in \dbR^d;\\
&\dis \hat U_t(\bx, \dbP_{\xi + \eta_t \1_{[t, T]}}) -\hat U_t(\bx, \dbP_\xi) = \dbE\Big[\pa_\mu \hat U_t(\bx, \dbP_\xi, \xi_t)\cd \eta_t\Big] + o(\|\eta_t\|_2),\q \forall \eta_t \in \dbL^2(\cF_t; \dbR^d).
\eeaa
In particular, it can be proved that the Lions' derivative $\pa_\mu \hat U$ takes the above specific structure. We can then define the higher order derivatives in the same manner. In particular, $\pa_{\tilde x\mu} \hat U_t(\bx, \mu, \tilde x)$ is the standard (finite-dimensional) derivative of $\pa_\mu \hat U$ with respect to $\tilde x$, and  $\pa_{\mu\mu} \hat U: [0, T]\times \dbD \times \cP_2(\dbD) \times \dbR^d \times \dbR^d \to \dbS^d$ involves two additional variables $(\tilde x, \bar x)\in \dbR^d\times \dbR^d$. Let $C^{1,2,2}_b([0, T]\times \dbD\times \cP_2(\dbD))$ denote the space of adapted and continuous functions $\hat U: [0, T]\times \dbD \times \cP_2(\dbD) \times \dbR^d \to \dbR$ such that the derivatives $\pa_t \hat U, \pa_\bx\hat U, \pa_\mu\hat U, \pa_{\bx\bx} \hat U, \pa_{\bx\mu} \hat U, \pa_{\tilde x,\mu} \hat U, \pa_{\mu\mu} \hat U$ exist, and are continuous and bounded.

We then define the path derivatives of $U$ by restricting those of $\hat U$ to the continuous paths, for example, $\pa_\mu U$ is defined by restricting $\pa_\mu \hat U$ on $[0, T]\times \dbX\times \cP_2(\dbX)\times \dbR^d$. It can be shown that the path derivatives of $U$ are independent of the choice of the extension $\hat U$ and thus are intrinsic to $U$. Moreover, we say $U\in C_b^{1,2,2}([0, T]\times \dbX \times \cP_2(\dbX))$ if there exists an extension $\hat U\in C_b^{1,2,2}([0, T]\times \dbD \times \cP_2(\dbD))$.

\bs

\no{\bf Proof of  Proposition \ref{prop-MFCIto}.}  We shall only prove \reff{Ito-common} in the integral form on $[0, T]$ for an arbitrary extension  $\hat U\in C_b^{1,2,2}([0, T]\times \dbD \times \cP_2(\dbD))$. For notational simplicity we denote $U=\hat U$ at below. Set $0=t_0<\cds<t_n=T$ with $t_i := {i\over n}T$. For $j=1,2$, denote
\beaa
&\dis X^{j,n}_t := \sum_{i=0}^{n-1} X^j_{t_i}\1_{[t_i, t_{i+1})}(t) + X^j_T \1_{\{t_n\}}(t),\q \mu^n := \dbP\circ (X^{2,n})^{-1};\\
&\dis X^{j,n, i,  l}_t := X^{j,n}_{\cd \wedge t_i} +  l X^j_{t_i, t_{i+1}} \1_{[t_{i+1}, T]},\q \mu^{n,i,  l} := \dbP\circ (X^{2,n, i,  l})^{-1},\q  l \in [0,1],
\eeaa
where $X^j_{t_i, t_{i+1}} := X^j_{t_{i+1}}-X^j_{t_{i+1}}$. Note that  $X^{j,n}_{\cd\wedge t_{i+1}} = X^{j,n}_{\cd\wedge t_{i}} + X^j_{t_i, t_{i+1}} \1_{[t_{i+1}, T]}$.  Then,
\beaa
&&\dis U_T(X^{1,n}, \mu^n) - U_0(X^{1,n}, \mu^n)  =\sum_{i=0}^{n-1} \Big[U_{t_{i+1}}(X^{1,n}, \mu^n) - U_{t_i}(X^{1,n}, \mu^n)\Big]\\
&&\dis= \sum_{i=0}^{n-1} \Big[\big[ U_{t_{i+1}}(X^{1,n}_{\cd\wedge t_i}, \mu^n_{\cd\wedge t_i})- U_{t_i}(X^{1,n}, \mu^n)\big] + \big[U_{t_{i+1}}(X^{1,n}, \mu^n) - U_{t_{i+1}}(X^{1,n}_{\cd\wedge t_i}, \mu^n_{\cd\wedge t_i})\big]\Big]\\
&&\dis = \sum_{i=0}^{n-1} \int_{t_i}^{t_{i+1}} \pa_t U_t(X^{1,n}_{\cd\wedge t_i}, \mu^n_{\cd\wedge t_i}) dt + \sum_{i=0}^{n-1}\int_0^1 \big[I^{n,i}_1( l) + I^{n,i}_2( l)\big] d l,
\eeaa
where, denoting by $(\tilde X^{2,n, i,  l}, \tilde X^2_{t_i, t_{i+1}})$ and $(\bar X^{2,n, i,  l}, \bar X^2_{t_i, t_{i+1}})$  conditionally independent copies of $(X^{2,n, i,  l}, X^2_{t_i, t_{i+1}})$, conditional on $\dbF^0$,
\beaa
I^{n,i}_1( l) &:=& \pa_\bx U_{t_{i+1}}(X^{1,n, i,  l}, \mu^{n,i,  l}) \cd X^1_{t_i, t_{i+1}},\\
 I^{n,i}_2( l) &:=& \dbE_{\cF_{t_{i+1}}}\big[\pa_\mu U_{t_{i+1}}(X^{1,n, i,  l}, \mu^{n,i,  l}, \tilde X^{2,n, i,  l}) \cd \tilde X^2_{t_i, t_{i+1}}\big].
\eeaa
Let $o({1\over n})$ denote  a generic term whose $\dbL^2$-norm vanishes faster than ${1\over n}$ when $n\to \infty$.
By the desired regularities, one can easily check that
\beaa
I^{n,i}_1( l) &=& \pa_\bx U_{t_{i+1}}(X^{1,n}_{\cd\wedge t_i}, \mu^{n}_{\cd\wedge t_i}) \cd X^1_{t_i, t_{i+1}} \\
&&+ \Big[\pa_\bx U_{t_{i+1}}(X^{1,n, i,  l}, \mu^{n,i,  l}) - \pa_\bx U_{t_{i+1}}(X^{1,n, i,0}, \mu^{n,i,0})\Big] \cd X^1_{t_i, t_{i+1}}\\
&=&\pa_\bx U_{t_{i+1}}(X^{1,n}_{\cd\wedge t_i}, \mu^{n}_{\cd\wedge t_i}) \cd X^1_{t_i, t_{i+1}} + \int_0^ l \Big[\pa_{\bx\bx} U_{t_{i+1}}(X^{1,n, i,  l'}, \mu^{n,i,  l'})  X^1_{t_i, t_{i+1}} \\
&&+ \dbE_{\cF_{t_{i+1}}}\big[\pa_{\bx\mu} U_{t_{i+1}}(X^{1,n, i,  l'}, \mu^{n,i,  l'}, \tilde X^{2,n, i,  l'})  \tilde X^2_{t_i, t_{i+1}}\big]\Big] d l'\cd X^1_{t_i, t_{i+1}}\\
&=& \int_{t_i}^{t_{i+1}}\pa_\bx U_t(X^1, \mu)\cd dX^1_t  +   l \pa_{\bx\bx} U_{t_{i+1}}(X^{1,n}_{\cd\wedge t_i}, \mu^{n}_{\cd\wedge t_i}) : (X^1_{t_i, t_{i+1}})(X^1_{t_i, t_{i+1}})^\top  \\
&& +  l\dbE_{\cF_{t_{i+1}}}\big[\pa_{\bx\mu} U_{t_{i+1}}(X^{1,n}_{\cd\wedge t_i}, \mu^{n}_{\cd\wedge t_i}, \tilde X^{2,n}_{\cd\wedge t_i}) : (\tilde X^2_{t_i, t_{i+1}}) (X^1_{t_i, t_{i+1}})^\top \big]\Big]  + o({1\over n})\\
&=& \int_{t_i}^{t_{i+1}}\pa_\bx U_t(X^1, \mu)\cd dX^1_t  +   l \int_{t_i}^{t_{i+1}}\pa_{\bx\bx} U_t(X^{1}, \mu) : \si^1_t(\si^1_t)^\top dt \\
&& +  l \int_{t_i}^{t_{i+1}}\dbE_{\cF_{t_{i+1}}}\big[\pa_{\bx\mu} U_t(X^1, \mu, \tilde X^2) : \tilde \si^{2,0}_t (\si^{1,0})^\top \big]dt  + o({1\over n});\\
I^{n,i}_2( l) &=&  \dbE_{\cF_{t_{i+1}}}\Big[\pa_\mu U_{t_{i+1}}(X^{1,n}_{\cd\wedge t_i}, \mu^{n}_{\cd\wedge t_i}, \tilde X^{2,n}_{\cd\wedge t_i}) \cd \tilde X^2_{t_i, t_{i+1}}\\
&& +  \big[\pa_\mu U_{t_{i+1}}(X^{1,n, i,  l}, \mu^{n,i,  l}, \tilde X^{2,n, i,  l}) - \pa_\mu U_{t_{i+1}}(X^{1,n, i, 0}, \mu^{n,i, 0}, \tilde X^{2,n, i, 0})\big]\cd \tilde X^2_{t_i, t_{i+1}}\Big]\\
&=& \dbE_{\cF_{t_{i+1}}}\Big[\pa_\mu U_{t_{i+1}}(X^{1,n}_{\cd\wedge t_i}, \mu^{n}_{\cd\wedge t_i}, \tilde X^{2,n}_{\cd\wedge t_i}) \cd \tilde X^2_{t_i, t_{i+1}}\Big]\\
&& + \int_0^ l \dbE_{\cF_{t_{i+1}}}\Big[\big[\pa_{\bx\mu} U_{t_{i+1}}(X^{1,n, i,  l'}, \mu^{n,i,  l'}, \tilde X^{2,n, i,  l'}) X^1_{t_i, t_{i+1}} \\
&&\q + \pa_{\mu\mu} U_{t_{i+1}}(X^{1,n, i,  l'}, \mu^{n,i,  l'}, \bar X^{2,n, i,  l'}, \tilde X^{2,n, i,  l'}) \bar X^2_{t_i, t_{i+1}}\\
&& \q +  \pa_{\tilde x\mu} U_{t_{i+1}}(X^{1,n, i,  l'}, \mu^{n,i,  l'}, \tilde X^{2,n, i,  l'}) \tilde X^2_{t_i, t_{i+1}}\big]\cd \tilde X^2_{t_i, t_{i+1}}\Big]d l'\\
&=& \dbE_{\cF_{t_{i+1}}}\Big[\int_{t_i}^{t_{i+1}}\pa_\mu U_t(X^1, \mu, \tilde X^2) \cd d\tilde X^2_t\Big] \\
&& +  l \dbE_{\cF_{t_{i+1}}}\Big[\int_{t_i}^{t_{i+1}}\big[\pa_{\bx\mu} U_t(X^1, \mu, \tilde X^2) : \tilde \si^{2,0}_t (\si^{1,0})^\top \\
&&\q + \pa_{\mu\mu} U_t(X^1, \mu, \bar X^2, \tilde X^2):  \tilde \si^{2,0}_t (\bar \si^{2,0})^\top +  \pa_{\tilde x\mu} U_t(X^1, \mu, \tilde X^2) : \tilde \si^{2}_t (\tilde\si^{2})^\top\Big] + o({1\over n}).
\eeaa
Then, noting that $\int_0^1  l d l = {1\over 2}$,
\beaa
&&U_T(X^1, \mu) - U_0(X^1, \mu) =U_T(X^{1,n}, \mu^n) - U_0(X^{1,n}, \mu^n)  + o(1)\\
&& = \int_0^T \Big[\pa_t U_t(X^1, \mu) dt + \pa_\bx U_t(X^1, \mu)\cd dX^1_t  +  {1\over 2} \pa_{\bx\bx} U_t(X^{1}, \mu) : \si^1_t(\si^1_t)^\top\Big] dt \\
&& +  \int_0^T\dbE_{\cF_t}\Big[\pa_{\bx\mu} U_t(X^1, \mu, \tilde X^2) : \tilde \si^{2,0}_t (\si^{1,0})^\top + {1\over 2}\pa_{\mu\mu} U_t(X^1, \mu, \bar X^2, \tilde X^2):  \tilde \si^{2,0}_t (\bar \si^{2,0})^\top \\
&&\q +  {1\over 2}\pa_{\tilde x\mu} U_t(X^1, \mu, \tilde X^2) : \tilde \si^{2}_t (\tilde\si^{2})^\top \Big]dt + o(1).
\eeaa
By sending $n\to \infty$, this is exactly \reff{Ito-common}.
\qed

\par

\end{document}